\documentclass[10pt]{article}
\expandafter\let\csname equation*\endcsname\relax
\expandafter\let\csname endequation*\endcsname\relax
\usepackage{bm}
\usepackage{enumerate}
\usepackage[top=1in, bottom=1in, left=1in, right=1in]{geometry}
\usepackage[dvipsnames]{xcolor}
\usepackage{mathrsfs}
\usepackage{amsfonts}\usepackage{multirow}
\usepackage{amssymb,dsfont}
\usepackage{caption,comment}
\usepackage{blkarray}
\usepackage{amsmath}
\usepackage{float}
\usepackage{footnote}
\usepackage{amssymb,amsthm}
\usepackage[toc,page]{appendix}
\usepackage{graphicx,afterpage}
\usepackage{epstopdf}
\usepackage{tabularx} 
\usepackage[normalem]{ulem}
\usepackage[pdftex,colorlinks=true]{hyperref}
\usepackage{authblk}
\usepackage{soul}
\hypersetup{CJKbookmarks,%
	bookmarksnumbered,%
	colorlinks,%
	linkcolor=black,%
	citecolor=black,%
	plainpages=false,%
	pdfstartview=FitH}
\numberwithin{equation}{section}
\numberwithin{figure}{section}

\newcommand\tabcaption{\def\@captype{table}\caption}

\newtheorem{thm}{Theorem}[section]

\newtheorem{aspt}[thm]{Assumption}

\newtheorem*{thm*}{Theorem}
\newtheorem{example}[thm]{Example}

\newtheorem{theorem}{Theorem}[section]
\newtheorem{lemma}[thm]{Lemma}

\newcommand{\cP}{\mathcal P}

\newcommand{\cT}{\mathcal T}

\renewcommand{\d}{\,\mathrm{d}}
\newcommand{\R}{\mathbb R}

\newcommand{\eps}{\varepsilon}

\newcommand{\supp}{\mathrm{supp}}

\usepackage[english]{babel}
\usepackage[
backend=biber,
style=alphabetic,
]{biblatex}

\usepackage{afterpage}

\newcommand\blankpage{%
	\null
	\thispagestyle{empty}%
	\addtocounter{page}{-1}%
	\newpage}

\title{Connections between sequential Bayesian inference and evolutionary dynamics}
\date{}

\author[1]{Sahani Pathiraja}
\author[2]{Philipp Wacker}

\affil[1]{School of Mathematics \& Statistics, UNSW Sydney, Australia, \url{https://orcid.org/0000-0002-0114-3164}}
\affil[2]{Faculty of Engineering, University of Canterbury, New Zealand, \url{https://orcid.org/0000-0001-8718-4313}}

\begin{document}
	
	\maketitle
	\begin{abstract}
		It has long been posited that there is a connection between the dynamical equations describing evolutionary processes in biology and sequential Bayesian learning methods. This manuscript describes new research in which this precise connection is rigorously established in the continuous time setting.  Here we focus on a partial differential equation known as the Kushner-Stratonovich equation describing the evolution of the posterior density in time. Of particular importance is a piecewise smooth approximation of the observation path from which the discrete time filtering equations, which are shown to converge to a Stratonovich interpretation of the Kushner-Stratonovich equation. This smooth formulation will then be used to draw precise connections between non-linear stochastic filtering and replicator-mutator dynamics.  Additionally, gradient flow formulations will be investigated as well as a form of replicator-mutator dynamics which is shown to be beneficial for the misspecified model filtering problem.  It is hoped this work will spur further research into exchanges between sequential learning and evolutionary biology and to inspire new algorithms in filtering and sampling. 
	\end{abstract}
	\tableofcontents

	\section{Introduction}

	It has been posited that there is a connection between sequential Bayesian inference and dynamical models describing evolutionary biological processes.  Understanding and studying this connection has the potential to provide valuable insights on improved algorithms for complex Bayesian inference and sampling tasks arising in a wide range of fields in the data science, engineering and machine learning.  Specifically, the key connection to sequential Bayesian estimation is via the so-called replicator-mutator partial differential equations \cite{kimura_diffusion_2024,hofbauer1985selection}, describing the time evolution of a large population of individuals with certain traits or attributes due to mutation and reproduction (or selection).
	Broadly speaking, sequential Bayesian estimation procedures bear striking similarity to the way species respond to evolutionary pressure moderated by a fitness landscape. The correspondence is as follows: 
	\begin{itemize}
		\item states or parameters $\leftrightarrow$ traits
		\item prior distribution $\leftrightarrow$ current population
		\item prediction (in the case of filtering or hidden markov models) $\leftrightarrow$ mutation 
		\item likelihood function $\leftrightarrow$ fitness landscape governing selection or birth-death
	\end{itemize}
	This connection has been discussed most notably in \cite{DelMoral1997,Shalizi2009,harper2009replicator,akyildiz2017probabilistic,Czegel2022}, primarily in the context of discrete time and discrete trait space problems. 
	Some of the earliest connections between discrete time particle based Bayesian updating and genetic mutation-selection models seem to have been in e.g. \cite{DelMoral1997,DelMoral2004,DelMoral2005}. \cite{Shalizi2009} raised awareness to the similarity to replicator equations specifically; they show how Bayesian updating corresponds to one step of a discrete-time and continuous-trait replicator equation without mutation. Around the same time, a similar point was made in \cite{harper2009replicator}. This connection to replicator equations with mutation was further extended to the setting of sequential in time inference with hidden Markov models (also known as sequential filtering, discrete time data assimilation)) in \cite{akyildiz2017probabilistic,Czegel2022}.  They showed that discrete time replicator-mutator dynamics consists of a sequence of (discrete in time) alternating mutation and updating steps, as in sequential filtering. The PhD thesis \cite{zhang2017controlled} draws some interesting connections to optimisation and interacting particle approaches to sequential filtering, e.g. the Feedback Particle Filter \cite{yang_multivariable_2012, laugesen_poissons_2014}.  In this manuscript we focus on making this connection precise in the continuous time and continuous trait space case, which to the best of our knowledge appears to be missing in the literature.  
	More recently, there has been interest in incorporating replicator or ``birth-death'' dynamics into sampling algorithms for optimisation and inversion tasks, see e.g. \cite{lu_accelerating_2019, lu_birthdeath_2023, chen_efficient_2024}. 
	
	\subsection{Replicator-mutator equations}
	
	The replicator-mutator equation is a broad class of dynamical systems modelling the response of a distribution of traits to evolutionary adaptation to an external fitness landscape.  Early influential work in establishing these models can be attributed to \cite{crowkimura,akin1979geometry,schuster1983replicator}; this equation is also sometimes referred to as the ``Crow-Kimura'' equation.  These models are popular in the mathematical evolution, biology, game theory and ecology literature, primarily in modelling the discrete-time evolution of a distribution of discrete-valued traits such as gene loci \cite{kopp_rapid_2014,kopp2018phenotypic,Matuszewski2014}. Their continuous time counterparts are typically the focus of analysis where their stability and geometric properties have been established for various fitness functionals, see e.g. \cite{Chalub2009,chalub_gradient_2020}.

	

	Our focus is on the continuous trait space setting which has received comparatively less attention, with the exception of \cite{kimura_stochastic_1965,cressman2006stability,gil_mathematical_2017,vlasic_markovian_2020}.     The continuous-time continuous-trait replicator-mutator partial differential equation (PDE) is given by  
	\begin{align}
		\partial_t \rho_t(x) &= \underbrace{\mathcal L^\star \rho_t(x)}_{\text{mutation}} + \underbrace{\rho_t(x) \left(\mathbb{E}_{z \sim \rho_t}[f_t(x,z)] - \mathbb{E}_{\rho_t}[f_t] \right)}_{\text{replication}} ,\label{eq:norm_repmut}
	\end{align}
	where $\rho_t(x) \geq 0$ 
	denotes the probability density function describing the distribution of traits $x \in \mathbb{R}^n$ in the population.  The mutation term is typically modelled after a diffusion process, where $\mathcal L^\star$ is the adjoint generator of a diffusion process, e.g., $\mathcal L^\star  = -\Delta$ for Brownian dynamics. If $\mathcal L^\star = 0$, i.e., no mutation exists, we call this the (pure) replicator equation.  
	The (potentialy non-local \& time-dependent) selection or fitness function is denoted by $f_{t}(x,z)$ and the net birth-death rate for a given trait $x$ at time $t$ is given by $\mathbb{E}_{z \sim \rho_t}[f_t(x,z)]$, where the subscript indicates that the expectation is taken over the $z$ variable only.  The expectation here is used to model ``mean-field'' type interactions within the population that affect an individual's fitness.  A simplified form of the replicator-mutator often appearing in the literature \cite{Tsimring1996,alfaro_explicit_2014,Alfaro2017} is the local PDE 
	\begin{align}
		\label{eq:norm_repmut_local}
		\partial_t \rho_t(x) &= \mathcal L^\star \rho_t(x)+ \rho_t(x) \left(f_t(x) - \mathbb{E}_{\rho_t}[f_t] \right),
	\end{align}
	where the non-local fitness function in \eqref{eq:norm_repmut} has been replaced by a fitness function that no longer depends on the current distribution of traits, only on the value of the trait itself.  An ubiquitous example that is connected to least squares estimation is the quadratic fitness, $f_t(x) = -\frac{1}{2}\|Hx - y_t\|_\Xi^2$.  Such a selection function penalises traits $x$ that have a large misfit to the (potentially) time varying data $y_t$ via a mapping $H$ and  preconditioned with a covariance matrix $\Xi$.  From an evolutionary biology perspective, $y_t$ may play the role of an ``optimal feature'' at time $t$ that is best adapted to the current environment and $x_t$ the trait that generates this physical feature in an individual.   

	\subsection{Sequential Bayesian inference or stochastic filtering}
    \label{sec:backfilt}
    Sequential Bayesian inference in continuous time (also stochastic filtering or non-linear filtering) is most often formulated in the following form.  Consider Euclidean spaces $X=\R^m$ and $Y=\R^n$, covariance matrices $\Sigma\in \R^{m\times m}, \Xi\in \R^{n\times n}$, sufficiently regular mapping $g: X\to X$ and $h:X\to Y$, and the following signal-observation pair, 
	\begin{align}
		\label{eq:sig}
		dX_t &= g(X_t)dt + \Sigma^{1/2} \, dW_t; \quad \text{(signal)}\\
		\label{eq:obs}
		dZ_t &= h(X_t)dt + \Xi^{1/2} \, dB_t \quad \text{(observation)}  .
	\end{align}
	The goal of stochastic filtering is to reconstruct the signal $X_t$ by means of the noisy observation path $\{Z_s\}_{s\leq t}$. Since $X_t$ cannot be uniquely identified from this data, the correct object to study is the conditional density of the hidden state $X_t$ at time $t$, given the observation filtration $\mathcal{Z}_t := \sigma(Z_s: s\leq t)$, which we denote by $p_t(x)$.  Throughout the article we will assume all distributions have densities with respect to the lebesgue measure.  It is known that $p_t$ evolves in time according to the Kushner-Stratonovich equation,  
	\begin{align}
		\label{eq:introkseqn}
		d p_t(x) = \mathcal{L}^\ast p_t(x) + p_t(x) \left(h(x) -  \mathbb{E}_{p_t}[h] \right)^\top \Xi^{-1} (dZ_t(\omega)  - \mathbb{E}_{p_t}[h] dt) 
	\end{align}
	where $\omega$ is used to denote a fixed realisation and $\mathcal{L}^\ast$ denotes the adjoint of the infinitesimal generator of \eqref{eq:sig} (i.e. the Kolmogorov forward operator), 
	\begin{align}
		\label{eq:generator}
		\mathcal{L}^\ast p_t(x) = -\nabla \cdot (p_t(x) \; g(x)) + \frac{1}{2}\nabla \cdot \left(\Sigma \;  \nabla p_t(x) \right). 
	\end{align}
	A direct comparison of terms between \eqref{eq:introkseqn} and \eqref{eq:norm_repmut_local} hints at a possible connection between these two PDEs, although \eqref{eq:norm_repmut_local} is a PDE driven by a finite dimensional ``rough'' path (in that it is not differentiable in time) whereas the fitness function $f_t(x)$ in \eqref{eq:norm_repmut_local} does not obviously coincide with such a term.  As we will demonstrate in Section \ref{sec:filt}, one way to establish this connection is to consider a piecewise linear approximation to the observation path, $Z_t^d$, and the corresponding evolution PDE of the conditional distribution of $X_t$ given the approximate data path $\{Z_s^d\}_{s \leq t}$.\\
	\\
	As solving \eqref{eq:introkseqn} is generally infeasible whenever $m,n$ is large, a popular numerical approach is to construct interacting diffusion processes whose law coincides with the solution of \eqref{eq:introkseqn}.  A notable example here is the ensemble Kalman-Bucy filter \cite{evensen2003ensemble, evensen94, bishop_mathematical_2023},  
	\begin{align*}
		d\bar{X}_t = g(X_t)dt + \Sigma^{1/2}dW_t  \, + \, \mathcal{Q}_t \Xi^{-1} (dZ_t - h(\bar{X}_t) dt \, + \, \sqrt{2} \Xi^{1/2} dW_t ) 
	\end{align*}
    where $\mathcal{Q}_t$ denotes the cross-covariance between $\bar{X}_t$ and $h(\bar{X}_t)$.  This approximate filter is implemented with a set of $N$ diffusion processes with $\mathcal{Q}_t$ replaced by its empirical $N$-particle version, so that the empirical distribution from the $N$ particles approximates $p_t$. 
	In practice, it is often necessary to adopt heuristic strategies such as covariance inflation and localisation to account for finite sample size effects and model errors \cite{anderson_adaptive_2007, Duc2020, hamill_accounting_2005, Mitchell2000, bishop_perturbations_2018, bishop_mathematical_2023, scheffler_dynamical_2022}.  We will demonstrate in Section 4 that the Crow-Kimura replicator-mutator equation \eqref{eq:norm_repmut} with a particular fitness function is directly related to covariance inflated ensemble Kalman-Bucy filtering. Our main contributions are summarised next.

	\subsection{Main contributions}

	Our first main contribution is in rigorously connecting replicator-mutator equations to stochastic filtering.  In doing so, we obtain a connection to continuous time Bayesian inversion with static data as a special case, which we describe in Section 2. We additionally investigate gradient flow aspects of this particular replicator-mutator with a non-local fitness function in Lemma \ref{lem:fisherrao}, where results for related forms have appeared recently \cite{chen_sampling_2024, maurais_sampling_2024,wang_measure_2024,lu_accelerating_2019, lu_birthdeath_2023, belhadji_weighted_2025, zhu_kernel_2024}.\\   
    \\
    More specifically, consider a particular Crow-Kimura replicator-mutator equation with \textit{non-local} fitness function, 
	\begin{equation}
		\label{eq:repmutlin}
		\begin{split}
			\partial_t \rho_t^d(x) &= \underbrace{\mathcal{L}^\ast \rho_t^d}_{\text{mutation}}  + \underbrace{\rho_t(x) \left(\mathbb{E}_{z \sim \rho_t}[f_t(x,z)] - \mathbb{E}_{\rho_t}[f_t] \right)}_{\text{replication}}\\ 
			f_t(x,z) &:= -\frac{r}{2}\|h(x)- \frac{dZ_t^d}{dt}\|_\Xi^2  + s \left \langle h(x) - \frac{dZ_t^d}{dt}, h(z) - \frac{dZ_t^d}{dt} \right \rangle_{\Xi}, \quad r > 0, s < r \\
			\mathcal{L}^\ast \rho_t^d &:=- \nabla \cdot (\rho_t^d(x) g(x)) + \frac{1}{2}  \nabla \cdot (\Sigma \nabla \rho_t^d(x)) 
		\end{split}
		\tag{$\textbf{CK}_d$}
	\end{equation}
	where the mutation term corresponds to evolution under the signal process \eqref{eq:sig} and $Z_t^d$ corresponds to a piecewise smooth approximation of an observation path generated by \eqref{eq:obs} (further details in Section \ref{sec:filt}).  A form of the non-local fitness function $f_t(x,z)$ with $H = I, r=1, \frac{dZ_t^d}{dt} = 0$ appears in the evolutionary biology literature (see e.g. \cite{cressman2006stability}).  In that context, the fitness function takes into account both a given trait's fitness on its own, but also beneficial or adversary effects of the group fitness through the $s \langle Hx-\frac{dZ_t^d}{dt}, Hz- \frac{dZ_t^d}{dt} \rangle_\Xi$ term.  More specifically, 
	\begin{itemize}
		\item $s =r$: the fitness function takes the form $f_t(x) = \|Hx-Hm\|_\Xi^2$ and traits $x$ with the property $Hx = Hm$ are most fit. Here ``conformity'' or ``collaboration'' (in the sense of being close to the population's feature average) is prioritised entirely over fit to the data term. 
		\item $s \in(0,r)$ corresponds to a fitness function with a trade-off between maximising data utility and conformity, respectively. 
		\item $s < 0$ corresponds to a fitness function that benefits from diversity or non-conformity: individual fitness can be improved by moving slightly from the population mean, even at the cost of increasing the misfit term. 
	\end{itemize}

	As will be demonstrated in Sections \ref{sec:filt} and \ref{sec:inflrepmut}, the local fitness function (with $s=0$) has an interpretation as a time dependent log-likelihood function.  The non-local form $(s \neq 0)$ can, in the linear-Gaussian setting, be directly related an ensemble Kalman-Bucy filter with both additive (e.g. \cite{hamill_accounting_2005}) and multiplicative  \cite{bishop_mathematical_2023} covariance inflation.    Our main result is Theorem \ref{theo:limitcrowkimura}, which we restate here informally, 
	\begin{thm*}
		\textbf{(Informal Theorem \ref{theo:limitcrowkimura}:) Limiting process of \eqref{eq:repmutlin}} Consider a time discretisation of the interval $[0,T]$ with step size $\delta_d$ such that $d \times \delta_d = T$.  Denote by $\{Z_t^d\}_{0 < t \leq T}$ a piecewise linear approximation with step size $\delta_d$ of an observation path generated by \eqref{eq:obs}.  Finally, let $\mu_t^d \propto \rho_t^d$ denote the unnormalised form of $\rho_t^d$ satisfying \eqref{eq:repmutlin}.  Then under some conditions, 
		\begin{align*}
			\mu_t^d \rightarrow q_t \quad \text{as} \quad d \rightarrow \infty, \enskip \forall \enskip  t \in [0, T],
		\end{align*}
		where $q_t$ satisfies the modified Zakai equation,
		\begin{align}
			\label{eq:modifzakai}
			d q_t(x) = \mathcal{L}^\ast q_t(x) -\frac{s}{2}h(x)^\top \Xi^{-1} h(x) q_t(x)  + (r-s)q_t(x) h(x)^\top \Xi^{-1} dZ_t, 
			\tag{$\textbf{ZK}$}
		\end{align}
		with $\mathcal{L}^\ast$ is as defined in \eqref{eq:repmutlin}.  
	\end{thm*}
	The required conditions and specific mode of convergence are detailed in Section \ref{sec:filt}.  This theorem allows us to bridge between a piecewise smooth replicator-mutator PDE and a ``generalised'' form of the fundamental equation of stochastic filtering which is a (stochastic)-PDE.  Additionally, for the special case $r=1, s = 0$, \eqref{eq:repmutlin} then corresponds to a replicator-mutator PDE with \textit{local} fitness function $f_t(x)$,
	\begin{align*}
		f_t(x) =  -\frac{1}{2}\|h(x)- \frac{dZ_t^d}{dt}\|_\Xi^2
	\end{align*}
	and the corresponding limiting process is given by the classical Zakai equation from non-linear filtering,
	\begin{align}
		\label{eq:zakaiintro}
		d q_t(x) = \mathcal{L}^\ast q_t(x) + q_t(x) h(x)^\top \Xi^{-1} dZ_t.
	\end{align}
	Interpretations of the limiting process \eqref{eq:modifzakai} and \eqref{eq:repmutlin} for $r > 0, s \neq 0$ are investigated further in Section \ref{sec:misspec} in the context of linear-Gaussian filtering.    
	It is evident that the solution of \eqref{eq:modifzakai} for $s \neq 0, r \neq 1$ no longer coincides with the filtering density.  Our final set of contributions (as summarised in the below table) gives an inference-based interpretation of the following non-local replicator-mutator equation,  
	\begin{equation}
		\label{eq:repmutLG}
		\begin{split}
			\partial_t \rho_t^d(x) &= \underbrace{\mathcal{L}^\ast \rho_t^d}_{\text{mutation}}  +  \underbrace{\rho_t(x) \left(\mathbb{E}_{z \sim \rho_t}[f_t(x,z)] - \mathbb{E}_{\rho_t}[f_t] \right)}_{\text{replication}} \\ 
			f_t(x,z) &:= -\frac{r}{2}\| Hx - \frac{dZ_t^d}{dt}\|_\Xi^2  + s \left \langle Hx - \frac{dZ_t^d}{dt}, Hz - \frac{dZ_t^d}{dt} \right \rangle_{\Xi}, \quad r > 0, s < r \\
			\mathcal{L}^\ast \rho_t^d &:=- \nabla \cdot (\rho_t^d(x) Gx + \frac{1}{2}  \nabla \cdot (\Sigma \nabla \rho_t^d(x)) 
		\end{split}
		\tag{$\textbf{LG-CK}_d$}
	\end{equation}
	where $Z_t^d$ is a piecewise linear approximation as described before, but for an observation path generated by \eqref{eq:sig}-\eqref{eq:obs} with $g(x):= Gx, \, G \in \mathbb{R}^{m \times m}$ and $h(x):=Hx, \, H \in \mathbb{R}^{n \times m}$.  In particular, when $s = 0$ but $r$ is any arbitrary positive value, \eqref{eq:repmutLG} describes the probability density evolution of multiplicative covariance inflated ensemble Kalman-Bucy filtering \eqref{eq:inflateenkbf}.  The non-local form (i.e. with $s \neq 0)$ can be seen as form of inflated ensemble Kalman-Bucy filtering involving both additive and multiplicative covariance inflation, as detailed in Section \ref{sec:inflrepmut}.  The key (and important difference) to standard implementations of inflation is in the specification of the weighting matrix $T$, see discussion below lemma \ref{lem:meanfieldproc} in section \ref{sec:inflrepmut}.

	\begin{center}
		\begin{tabularx}{1\textwidth} { 
				| >{\centering\arraybackslash}X 
				| >{\centering\arraybackslash}X 
			| >{\centering\arraybackslash}X | }
		\hline
		\textbf{Replicator-Mutator} & \textbf{Filtering}  & \textbf{Equivalence result (as $d \rightarrow \infty$)}  \\
		\hline
		\eqref{eq:repmutlin} with $r =1,s=0$  & Zakai equation \eqref{eq:zakaiintro}  & Theorem \ref{theo:limitcrowkimura}  \\
		\hline
		\eqref{eq:repmutLG} with $r>0,s=0$ & Multiplicative covariance inflated linear-Gaussian EnKBF \eqref{eq:inflateenkbf}  & Lemma \ref{lem:meanfieldproc}  \\
		\hline
		\eqref{eq:repmutLG} with $r>0,s \neq 0$  & Additive +  Multiplicative covariance inflated linear-Gaussian EnKBF & Lemma \ref{lem:meanfieldproc}  \\
		\hline
	\end{tabularx}
	\end{center}
	
	 Finally, In section \ref{sec:biasvar} we further investigate the performance of \eqref{eq:repmutLG} with various $r,s$ values for linear-Gaussian filtering in the presence of model misspecification due to an unknown bias term in the signal dynamics. It is demonstrated both analytically (section \ref{sec:mseanalysis}, Lemma \ref{lem:optimalrs}) and numerically (section \ref{Sec:numericsmissp}) that there are infinitely many $(r,s)$ pairs minimising the time asymptotic mean squared error (MSE). Most importantly, we show in Lemma  \ref{lem:covimpact} that there is a unique $(r,s)$ pair that simultaneously minimises asymptotic MSE and gives rise to an asymptotic covariance that equals MSE (as is achieved with the Kalman-Bucy filter in the perfect model setting).  We give exact expressions for this unique $(r,s)$ pair in the scalar setting in terms of the system parameters, demonstrating that in most cases this corresponds to the \textit{non-local} form of \eqref{eq:repmutLG} with $s \neq 0$.  It remains to be seen if such a non-local fitness function has similar benefits for misspecified filtering in the non-linear case; we leave a detailed analysis of this to future work.  We hope that the connections established in this paper spur further exchanges between the field of evolutionary dynamics and sequential learning and inference.

\subsection{Notation}

We define some notation that will be used throughout the manuscript. 
\begin{itemize}
	\item [] Set $X = \mathbb R^m$ the trait space, $H: X\to Y$ be a mapping, and $Y = \R^n$. 
	\item [] $m$ = dimension of state/trait space 
	\item [] $n$ = dimension of observation/fitness space
	\item [] $\| x\|_\Xi^2 = x^\top \Xi^{-1}x$
	\item [] $\dot{m}_t$ is used to denote $\frac{\d m_t}{\d t}$ for any vector $m_t \in \mathbb{R}^d$ depending only on $t$. 
	\item [] $\partial_t p_t(x)$ is used to denote $\frac{\partial p_t(x)}{\partial t}$ for $\rho_t(x): [0,\infty) \times \mathbb{R}^m \rightarrow \mathbb{R}$
	\item [] $\mathbb{E}_{p_t}[f] = \int f(x) \rho_t(x)dx$ and where necessary, we specify the variable in the subscript to indicate the variable to which the integration operation applies, i.e. $\mathbb{E}_{z \sim \rho_t}[f_t(x,z)] = \int f_t(x,z) \rho_t(z) dz$.  
	\item [] $f_t(x,z): \mathbb{R}^m \times \mathbb{R}^m \rightarrow \mathbb{R}$ is the non-local fitness landscape at time $t$ (omission of the $t$ subscript is used to indicate a time-independent fitness landscape).  Frequency dependent fitness of trait $x \in \mathbb{R}^m$ is indicated by the term $\mathbb{E}_{z \sim \rho_t}[f_t(x,z)] \equiv \int f_t(x,z) \rho_t(z) dz$.  
	\item [] $p_t$ and $q_t$ is used to refer to the  normalised and unnormalised filtering density satisfying the Kushner-Stratonovich and Zakai equations  respectively. 
	\item [] $\rho_t$ and $\mu_t$ is used to refer to the normalised and unnormalised density function from Crow-Kimura replicator-mutator respectively. 
\end{itemize}

\section{The replicator equation \& continuous time Bayesian inversion}

Before presenting the rigorous connection between sequential filtering and replicator-mutator equations with time dependent fitness, we present the case of static fitness functions and continuous time Bayesian inversion.  Inversion \cite{engl1996regularization} refers to the task of inferring an unknown parameter $x\in X =\R^m$ from a noisy measurement $y\in Y=\R^n$ of form
\begin{equation}
	\label{eq:bayesian_inverse}
	y = h(x) + \varepsilon,
\end{equation}
where $H: X \to Y$ is the so-called forward operator and $\varepsilon$ is measurement noise, commonly assumed to be zero-mean Gaussian, $\varepsilon\sim \mathcal N(0,\Xi)$. The operator $H$ is analogous to the observation operator $h(x)$ in the filtering setting. Unfavorable properties of the operator and the noise makes a naive inversion procedure ill-defined, so additional regularisation of the inversion procedure is necessary. The Bayesian approach to inversion (\cite{stuart2010inverse}) requires a prior distribution $q_0$ of the unknown parameter $x$ and produces the Bayesian posterior measure $\nu^y$ via
\begin{equation}
	\label{eq:posterior}
	\frac{\d \nu^y}{\d \nu_0}(x) \propto \exp\left(-\frac{1}{2}\|h(x)-y\|^2_\Xi\right).
\end{equation}
This inversion procedure constitutes a one-step method from the prior $\nu_0$ to the posterior $\nu^y$. In the following, we assume that all measures have a Lebesgue density, which we denote $p_0$ and $p^y$.  A commonly adopted prior $p_0$ is the multivariate Gaussian $\mathcal N(m_0, C_0)$.  

When $m,n$ are large and the prior and posterior are significantly ``different'' it can be advantageous to gradually transform the prior to the posterior, either over a finite or infinite time horizon.  When this is done over a finite time interval, this procedure is known as tempering or annealing in the sequential monte carlo literature \cite{Gelman1998,Neal2001,chopin_connection_2024} and also the homotopy approach to Bayesian inversion (e.g. \cite{reich2011dynamical,blomker2022continuous, chada_tikhonov_2019}).  More specifically, this approach modifies the single step from prior $\mu_0$ to posterior $\mu^y$ into a smooth transition by introducing a pseudotime $t\in[0,1]$, and intermediate measures with density $p_t$ via
\[p_t(x)\propto \exp\left(-\frac{t}{2}\|h(x)-y\|^2_\Xi\right)\cdot p_0(x) \]
such that $p_1 = p^y$. 
As outlined in e.g. \cite{pidstrigach2023affine} it can be seen that a family of probability densities defined via $p_t(x) \propto \exp(tf(x)) p_0(x)$ is the solution of the infinite-dimensional system of ODEs\
\begin{equation}
	\label{eq:repeq_pointwise} \frac{\d p_t(x)}{\d t}= \frac{\d}{\d t} \left(\frac{\exp(t f(x)) p_0(x)}{\mathbb{E}_{p_t}[f]} \right) = (f(x) - \mathbb{E}_{p_t}[f])  p_t(x),
\end{equation}
This is identical to the pure replicator equation, i.e., the replicator-mutator equation without any mutation, with the fitness function being the log-likelihood,
\begin{align*}
	f(x) = -\frac{1}{2} \|h(x)-y\|^2_\Xi.
\end{align*}
In the simplest setting, $y$ is a time-independent feature, which then means that the fitness function is also static.  From a sampling point of view, it is then possible to construct both deterministic and stochastic schemes which have the replicator equation as their density evolution equation.  To do so, we switch to the so-called \textit{Eulerian} perspective on this problem, as presented e.g. in \cite{reich2011dynamical}.  The goal is to construct a vector field $(t,x)\mapsto v(t,x)$ such that the family of diffeomorphisms $T_t: X\to X$ with $T_t(x_0) := x(t)$ on trait space defined by solutions of the ODE
\begin{equation}
	\dot x(t) = v(t, x(t)),\, x(0)=x_0 \label{eq:ODE_homotopy}
\end{equation}
smoothly pushes forward the initial population distribution to the distribution at a later time via
\begin{equation}
	(T_t)_\# p_0 = p_t.
\end{equation}
The advantage of this perspective is that if we were able to find such a vector field $v$, then we can approximate the solution of \eqref{eq:repeq_pointwise} by sampling $J$ particles $\{x_0^{(i)}\}_{i=1}^J\sim p_0$, evolve them according to \eqref{eq:ODE_homotopy}, and the resulting ensemble $\{x^{(i)}(t)\}_{i=1}^J$ then constitutes valid samples from $p_t$. \cite[Theorem 5.34]{villani2021topics} states that such a velocity field $v$ satisfies the continuity equation
\begin{equation}
	\label{eq:homotopy_conteq} \partial_t p_t(x) = - \nabla_x \cdot (v(t,x) p_t(x)).
\end{equation}
On the other hand, comparing the right hand sides of \eqref{eq:repeq_pointwise} and \eqref{eq:homotopy_conteq} means that the vector field needs to be a solution of the Poisson equation
\begin{equation}
	\label{eq:homotopy_poisson} -\nabla_x \cdot (v(t,x) p_t(x)) = (f(x) - \mathbb{E}_{p_t}[f])  p_t(x).
\end{equation}
which leads us back to the familiar pure replicator dynamics.  It is worthwhile noting that this equation also arises in the construction of interacting particle filtering algorithms, where the fitness function is time-dependent due to the time-varying data term \cite{crisan_approximate_2010,laugesen_poissons_2014, pathiraja_mckean-vlasov_2021}.  Finally, it is possible to include a mutation or ``exploration term'' to aid in generating samples when $m$ is large (i.e. the underlying trait space is high dimensional), and is even necessary for particle-based implementations of the above (e.g. \cite{DelMoral1997,DelMoral2006,Moral2014, chen_efficient_2024, lu_accelerating_2019,lu_birthdeath_2023,zhang_mean-field_2019}).  In this setting, one arrives at a connection to the Crow-Kimura replicator-mutator equation, albeit with a static (time independent) fitness functional. 

In the remainder of this section we discuss some geometric properties of the replicator equation, showing that the continuous trait-space replicator equation follows a gradient flow of the mean fitness energy functional with respect to the Fisher-Rao metric.  This extends the well-known result that the discrete trait space replicator equation is a gradient flow with respect to the Shahshahani metric (i.e. the finite dimensional version of the Fisher-Rao metric), see for example Theorem 7.8.3 in \cite{hofbauer1998evolutionary}, \cite{fujiwara1995gradient} (for the entropy rather than fitness functional), \cite{Harper2009}, \cite{Chalub2009} and \cite{Chalub2021} (in the case of two traits/species) or \cite{baez2021fundamental} and \cite{baez2016relative} (including some interesting remarks about speed and acceleration in this metric) and \cite{raju_lie_2020, raju_variational_2018} for a control-theoretic perspective.  There is comparatively much less work on gradient flow interpretations for continuous trait spaces, which we now focus on.


We start by reminding ourselves of the basics of information geometry. We define $\cP$ as the manifold of absolutely continuous probability measures on $\R^n$. Every $p\in \cP$ will be identified with its (Lebesgue) density $p(x)$. At a $p\in \cP$ such that $p(x) > 0$ everywhere, the tangent space of $\cP$ is given by $\cT_p \cP = \{ \sigma\in C^\infty(\R^n): \int \sigma(x)\d x = 0\}$. If $p$ has support $\supp(p)\subsetneq \R^n$, then the tangent space is $\cT_p \cP = \{ \sigma\in C^\infty(\R^n): \sigma|_{\supp(p)^C} = 0, \int \sigma(x)\d x = 0\}$. The associated cotangent space is given by $\cT_p^\star \cP = \{\phi\in C^\infty(\R^n)\} /\sim$, where the equivalence relation $\sim$ is defined as $\phi_1\sim\phi_2$ if and only if $(\phi_1-\phi_2)|_{\supp(p)} \equiv \text{const}$. The dual pairing between $\phi\in \cT_p^\star\cP$ and $\sigma\in \cT_p\cP$ is then canonically defined as $\langle \sigma, \phi\rangle = \int \sigma(x)\phi(x)\d x$. The Rao-Fisher metric on the tangent space is given by $g_p(\sigma_1,\sigma_2) = \int \frac{\sigma_1(x)}{p(x)}\frac{\sigma_2(x)}{p(x)}\d p(x)$. This still makes sense on the boundary, i.e. if $p(x) = 0$ for some $x$, because then  $\sigma(x) = 0$ for $\sigma\in\cT_p\cP$ and we interpret the resulting expression $\frac{0}{0} = 0$. 
The metric defines an invertible metric tensor $G(p) : \cT_p\cP \to \cT_p^\star \cP$ via $g_p(\sigma_1,\sigma_2) = \langle \sigma_1, G(p)[\sigma_2]\rangle$ (or equivalently, $G(p)[\sigma] = g_p(\cdot, \sigma)$), which in this case is given explicitly by
\[ (G(p)[\sigma])(x) = \begin{cases}
	\frac{\sigma(x)}{p(x)} &\text{ if } p(x) > 0\\
	0 &\text{ else, }
\end{cases}\]
which is consistent since $\sigma \in \cT_p\cP$ is required to vanish on the set $\supp(p)^C$, anyway. In fact, any map identical to $G(p)$ up to a global constant on $\supp(p)$, and with arbitrary values on $\supp(p)^C$, would give a valid representative due to the quotient structure on $\cT_p^\star \cP$. 

The identity $g_p(\sigma_1,\sigma_2) = \langle \sigma_1, G(p)[\sigma_2]\rangle$ then holds true since
\begin{align*}
	\langle \sigma_1, G(p)[\sigma_2]\rangle &= \int \sigma_1(x)  G(p)[\sigma_2](x) \d x = \int \sigma_1(x)\frac{\sigma_2(x)}{p(x)}\d x = \int \frac{\sigma_1(x)}{p(x)}\frac{\sigma_2(x)}{p(x)}\d p(x)
\end{align*}

The metric tensor associated to the Fisher-Rao metric has inverse $G(p)^{-1}:\cT_p^\star\cP \to \cT_p \cP$
\[
(G(p)^{-1}[\phi])(x) = \left(\phi(x) - \int \phi(y) \d p(y)\right) p(x)
\]

which is well-defined since $G(p)\left[G(p)^{-1}[\phi] \right]= \phi - \int\phi(y)\d p(y)\sim \phi$.  We now recall some facts on gradient flows.  Let $\mathcal{P}$ denote a linear space; then the time evolution of $p_t \in \mathcal{P}$ is a gradient flow if it can be written as 
\begin{align*}
	\partial_t p_t = -\mathcal{K}(p_t)\mathcal{F}'(p_t)
\end{align*}
where $\mathcal{F}: \mathcal{P} \rightarrow \mathbb{R}$ is an energy functional,  $\mathcal{F}': \mathcal P \to T_p^\ast \mathcal P$ is the Frechet derivative and $\mathcal{K}(p): T_p^\ast \mathcal{P} \rightarrow T_p\mathcal{P}$ is a linear operator characterising the dissipation mechanism (loosely, giving meaning to how quickly $\mathcal{F}$ increases/decreases).  In this context we are most interested in the case where $\mathcal{K}$ is related to the metric tensor: The dissipation mechanism $\mathcal{K}(p)$ is then taken to be $\mathcal{G}(p)^{-1}$. We are now ready to state our main result of this section, the proof of which can be found in section \ref{sec:prooffisherrao}.
\begin{lemma}
	\label{lem:fisherrao}
	The pure replicator equation \eqref{eq:norm_repmut} with $\mathcal L^\star = 0$ and with frequency dependent but time-independent fitness $f_t(x,z) = f(x,z)$ performs a gradient flow of the average fitness functional $ \mathcal F(p): \mathcal{P}(\mathbb{R}^d) \rightarrow \mathbb{R}$ 
	\begin{align}
		\label{eq:expecfit}
		\mathcal F(p) = -\frac{1}{2} \iint f(x,z) p(z) p(x) d z d x
	\end{align}
	with respect to the Fisher-Rao metric.  
\end{lemma}

Similar results for frequency-dependent fitness functionals that are closely related to but not special cases of \eqref{eq:expecfit} have recently appeared in the sampling and machine learning literature (e.g. for the Kullback-leilber divergence \cite{chen_sampling_2024, maurais_sampling_2024,wang_measure_2024,lu_accelerating_2019, lu_birthdeath_2023}).  More specifically, \cite{lu_accelerating_2019} shows that a dynamical system similar to the Crow-Kimura replicator-mutator equation is a gradient flow of the Kullback-Leibler divergence with respect to the Wasserstein-Fisher-Rao metric.  Around the same time as this work, \cite{zhu_kernel_2024} provided an extensive rigorous study of kernelised Fisher-Rao gradient flows of maximum mean discrepancy (MMD) and \cite{belhadji_weighted_2025} proposed MMD flows in the Wasserstein-Fisher-Rao geometry.  The corresponding MMD$^2$ functional evaluated on kernel mean embeddings of the current and target distribution can be seen as a special case of the non-local functional \eqref{eq:expecfit}.  

\section{Time varying Crow-Kimura replicator-mutator and Stochastic filtering}
\label{sec:filt}

Recall the stochastic filtering problem as described in section \ref{sec:backfilt} and that the main quantity of interest in filtering is the conditional density $p_t(x)$.
The unnormalised density $q_t$ evolves according to the Zakai equation 
\begin{align}
	\label{eq:zakai}
	dq_t = \mathcal{L}^\ast q_t(x) + q_t h(x)^T \Xi^{-1} dZ_t 
\end{align}
which can be equivalently expressed in Stratonovich form as (see e.g. \cite{hu_approximation_2002}, \cite{pathiraja_mckean-vlasov_2021})
\begin{align}
	\label{eq:zakaistrat}
	d q_t(x) = \mathcal{L}^\ast q_t(x) -\frac{1}{2} h(x)^\top \Xi^{-1} h(x) q_t(x) +  q_t(x) h(x)^\top \Xi^{-1} \circ dZ_t 
\end{align}
In order to connect to discrete measurement processes more commonly encountered in practice, as well as to the Crow-Kimura equation replicator-mutator equation, consider a piecewise smooth observation process.  The piecewise smooth approach to approximating rough signals has a long history in robust filtering \cite{crisan_robust_2013} and has been well-studied more generally by the so-called Wong-Zakai style theorems and in the context of rough path theory \cite{kelly_smooth_2016, pathiraja_l_2024}.  In this section, we will consider a very specific form of piecewise smooth approximation using piecewise linear interpolations of Brownian noise.  More precisely, consider a partition of the time interval $[0,T]$, $0 < t_1 < t_2 \dots < t_d = T$ with time-step $t_{i+1} - t_i = \delta_d$ for all $i$ and such that $\delta_d \rightarrow 0$ as $d \rightarrow \infty$ (i.e. $\delta_d = \frac{T}{d}$).  Define the piecewise linear approximation to a Brownian path as
\begin{align*}
	B_t^d = B_{t_i} + \frac{t - t_i}{\delta_d} (B_{t_{i+1}} - B_{t_i}), \quad t \in [t_i, t_{i+1}),
\end{align*}
which is piecewise differentiable with time derivative 
\begin{align*}
	\frac{dB_t^d}{dt} = \frac{1}{\delta_d} (B_{t_{i+1}} - B_{t_i}), \quad t \in [t_i, t_{i+1}).
\end{align*}
Then a piecewise smooth version of \eqref{eq:obs} can be constructed, for all $i = 1, 2, \dots d$ at 
\begin{align}
	\label{eq:smoothobsnew}
	\frac{dZ_t^d}{dt} = h(x_{t_i}^\ast) + \Xi^{1/2} \frac{dB_t^d}{dt}, \quad t  \in [t_i, t_{i+1})
\end{align}
where the notation $x_t^\ast$ is used to denote the true hidden state trajectory or reference trajectory that generates the observed measurement.  This form allows us to connect more easily to observation models more commonly encountered in practice, i.e.,  
\begin{align*}
	y_t(x) := \frac{dZ_t^d}{dt} = h(x_t) + \epsilon_t, \quad \epsilon_t \sim \mathcal{N}(0, R) 
\end{align*}
and as we will see in the remainder of the section, to also build a bridge between the replicator-mutator equation and the Kushner-Stratonovich equation.  
It is worth noting that we consider a piecewise smooth approximation of the observation noise term only, rather than a smooth approximation of the entire observation trajectory as is done more traditionally in the stochastic filtering literature \cite{hu_approximation_2002, crisan_approximate_2010}.  Specifically, \cite{hu_approximation_2002} consider the following approximation 
\begin{align*}
	Z_t^\Pi = Z_{t_i} = \frac{Z_{t_{i+1}} - Z_{t_i}}{t_{i+1} - t_i}(t-t_i), \enskip t \in [t_i, t_{i+1}) 
\end{align*}
where $Z_t$ is a fixed realisation of \eqref{eq:obs} and $Z_t^\Pi$ denotes the corresponding piecewise approximation.  This then yields (assuming $t_{i+1} - t_i = \delta_d$),
\begin{align*}
	\frac{dZ_t^\Pi}{dt} = \frac{Z_{t_{i+1}} - Z_{t_i}}{t_{i+1} - t_i} = \frac{1}{\delta_d} \int_{t_i}^{t_{i+1}} dZ_s  = \int_{t_i}^{t_{i+1}} h(x_s^\ast)ds + \frac{dB_t^d}{dt}, \quad t \in [t_i, t_{i+1})
\end{align*}
so that the observation involves a time integrated version of the hidden state, $\int_{t_i}^{t_{i+1}} h(x_s^\ast)ds$ rather than $h(x_{t_i}^\ast)$ as in \eqref{eq:smoothobsnew}, which may be more practically relevant particularly when the time between observations is large.  This distinction is primarily for comparison to measurement models encountered in practice, both approximate forms can be shown to have valid limiting forms.   

In regards to the Crow-Kimura replicator-mutator equation, for the remainder of this section, we focus on the following time-varying quadratic fitness landscape with $s < r$, $r > 0$,  
\begin{align}
	\label{eq:genquadfit}
	f_t(x,z) = -\frac{r}{2}\|h(x) - \frac{dZ_t^d}{dt}\|_{\Xi}^2 + s \left \langle h(x) - \frac{dZ_t^d}{dt}, h(z) - \frac{dZ_t^d}{dt} \right \rangle_{\Xi} 
\end{align}
and time-varying optimal feature $\frac{dZ_t^d}{dt}$ given by \eqref{eq:smoothobsnew}.  Notice that with the definition of $\frac{dZ_t^d}{dt}$ in \eqref{eq:smoothobsnew}, $f_t(x,z)$ is a piecewise constant in time functional in $x$.  In the remainder of this section, we will establish equivalence (in a sense to be made precise), between the the C-K replicator-mutator with \eqref{eq:genquadfit} and a modified form of the Zakai equation.  For the special case of $r=1,s=0$ in \eqref{eq:genquadfit}, we will see that C-K replicator-mutator converges to the standard Zakai equation \eqref{eq:zakai} as $d \rightarrow \infty$.   The following reformulation of the Crow-Kimura replicator-mutator equation with \eqref{eq:genquadfit} will be a useful aid.  Its proof can be found in section \ref{sec:proofckreform}.   

\begin{lemma}
	\label{lem:ckreform}
	Consider the Crow-Kimura replicator-mutator equation 
	\begin{align}
		\label{eq:repmutnonlin}
		\partial_t \rho_t(x) = - \nabla \cdot (\rho_t(x) g(x)) + \frac{1}{2} \nabla \cdot (\Sigma \rho_t(x))  +  \rho_t(x) (\mathbb{E}_{z \sim \rho_t}[f_t(x,z)] -  \mathbb{E}_{\rho_t}[f_t]) 
	\end{align}
	with fitness landscape given by \eqref{eq:genquadfit}.  This equation can be expressed in the form 
	\begin{align}
		\nonumber 
		\partial_t \rho_t(x) = \mathcal{L}^\ast \rho_t(x) +  (r-s) &\left( - \frac12\left( h(x)^\top\Xi^{-1} h(x)   - \mathbb{E}_{\rho_t}[h^\top\Xi^{-1} h]\right) + \left(h(x) - \mathbb{E}_{\rho_t}[h] \right)^\top\Xi^{-1}\frac{dZ_t^d}{dt} \right)\rho_t(x) \\
		\label{eq:repasKS}
		-\frac{s}{2} &\left( \| h(x) - \mathbb{E}_{\rho_t}[h] \|_\Xi^2  - \mathbb{E}_{\rho_t}\left[\| h(x) - \mathbb{E}_{\rho_t}[h] \|_\Xi^2 \right] \right) \rho_t(x). 
	\end{align}
	where $\mathcal{L}^\ast$ denotes the adjoint of the generator of the diffusion process \eqref{eq:sig}.  Additionally, the unnormalised form is given by 
	\begin{align}
		\label{eq:unnormCK}
		\partial_t \mu_t(x) = \mathcal{L}^\ast \mu_t(x) +  \left(   - \frac{r}{2} h(x)^\top\Xi^{-1} h(x) +  (r-s) h(x)^\top \Xi^{-1}\frac{dZ_t^d}{dt}  \right) \mu_t(x) 
	\end{align}
	
\end{lemma}

Before presenting the main theorem of this section, consider the following simple motivating example to demonstrate why as $d \rightarrow \infty$, the Crow-Kimura replicator-mutator with fitness landscape \eqref{eq:genquadfit} converges to a pde driven by Stratonovich rather than Ito noise.  It should be noted that although the pde is driven by Stratonovich noise, it can be transformed to an Ito version from which the familiar Kushner-Stratonovich equation can be recovered (for $r=1, s=0$).    

\begin{example}
	\textbf{Simplified one dimensional replicator-mutator.}  Consider the 1d linear-Gaussian filtering problem with trivial signal dynamics, i.e. $g(x) = 0, \Sigma = 0, h(x) = Hx, \Xi = 1$. In this case, the  unnormalised Crow-Kimura replicator-mutator equation \eqref{eq:unnormCK} with $r=1,s=0$ takes the form,
	\begin{align}
		\label{eq:cklingauss}
		\partial_t \mu_t^d(x) =  - \frac12 (Hx)^2 \mu_t^d(x)   + \mu_t^d(x) Hx  \frac{dZ_t^d}{dt} 
	\end{align}
	The Zakai equation (unnormalised filtering pde) \eqref{eq:zakai} has the following (Stratonovich) representation,
	\begin{align}
		\label{eq:zakaistratlingauss}
		\partial_t q_t(x) &=  - \frac{1}{2 }(Hx)^2  q_t(x) +  q_t(x) Hx \circ dZ_t 
	\end{align}
	To help demonstrate why the limit of the smooth approximated noise in  \eqref{eq:cklingauss} must indeed be of Stratonovich type, consider the following pde driven by finite dimensional Ito noise,
	\begin{align}
		\label{eq:zakaiitolingauss}
		\partial_t \rho_t(x) = - \frac{1}{2}(Hx)^2  \rho_t(x) +  \rho_t(x) Hx  dZ_t 
	\end{align}
	The following numerical experiment demonstrates empirically the convergence of \eqref{eq:cklingauss} to \eqref{eq:zakaistratlingauss} rather than \eqref{eq:zakaiitolingauss}.  A sequence of smooth observations over $[0,T]$ with step size $\delta_d$ is generated as 
	\begin{align}
		\label{eq:piecelin}
		B_t^d = B_{t_n} + \left( \frac{t - t_n}{t_{n+1} - t_n} \right) (B_{t_{n+1}} - B_{t_{n}}), \quad t \in [t_{n}, t_{n+1}) 
	\end{align}
	where $B_t$ corresponds to a Brownian motion so that $(B_{t_{n+1}} - B_{t_{n}}) \sim \mathcal{N}(0, \delta_d)$.  The obs increment $\frac{dZ_t^d}{dt}$ is then defined as in \eqref{eq:smoothobsnew}.  With the approximation \eqref{eq:piecelin}, $\frac{dZ_t^d}{dt}$ is constant in the time interval $[t_n, t_{n+1})$ for every $n = 0,1,2, \dots$, and to emphasise the lack of dependence on $t$, we denote it by a random variable $\xi_n$ where
	\begin{align*}
		\xi_n \sim  \mathcal{N} \left(Hx_0^\ast, \frac{1}{\delta_d} \right) 
	\end{align*}
	for some fixed $x_0^\ast$ denoting the true hidden state at time $0$.
	By interpreting \eqref{eq:cklingauss} as a linear pde of the form 
	\begin{align*}
		\partial_t \mu_t^d(x) = A_t(x)\mu_t^d(x)
	\end{align*}
	with piecewise smooth in time coefficient 
	$A_t(x):=  - \frac{1}{2}(Hx)^2 + Hx \xi_n$
	for $t \in [t_n, t_{n+1})$ , it can be discretised via the usual euler scheme over $[0,T]$ with time step $\Delta$ and $t_{i+1} - t_i = \Delta t < \delta_d$ for all $i=0,1,2 \dots $. Let $\tilde{\mu}_i(x)$ denote the approximation to $\mu_t^d(x)$ at $t = t_i$,
	\begin{align*}
		\tilde{\mu}_{i+1}(x) = \tilde{\mu}_{i}(x) + \Delta t A_{i}(x)\tilde{\mu}_{i}(x)
	\end{align*}
	with 
	$A_{i}(x) = - \frac{1}{2}(Hx)^2 + Hx \xi_{n_i}$ 
	for $n_i$ such that $t_{i} \in [t_{n_i}, t_{n_i+1})$.  Similarly,  \eqref{eq:zakaistratlingauss} and \eqref{eq:zakaiitolingauss}  can be simulated with Euler-Maruyama schemes with the same time step $\Delta t$.  There the observation path is a solution of $dZ_t = Hx_0^\ast dt + dB_t$, simulated at a fine time interval. Figure \ref{fig:itostratexample} shows the results for a single time instant.  Importantly, the Crow-Kimura replicator equation \eqref{eq:repasKS} (cyan line) coincides closely with \eqref{eq:zakaistratlingauss} (black line) (up to normalisation), while the Ito version \eqref{eq:zakaiitolingauss} (pink line) is significantly different. 
	
	\begin{figure}[h]
		\centering
		\includegraphics[width=0.7\linewidth]{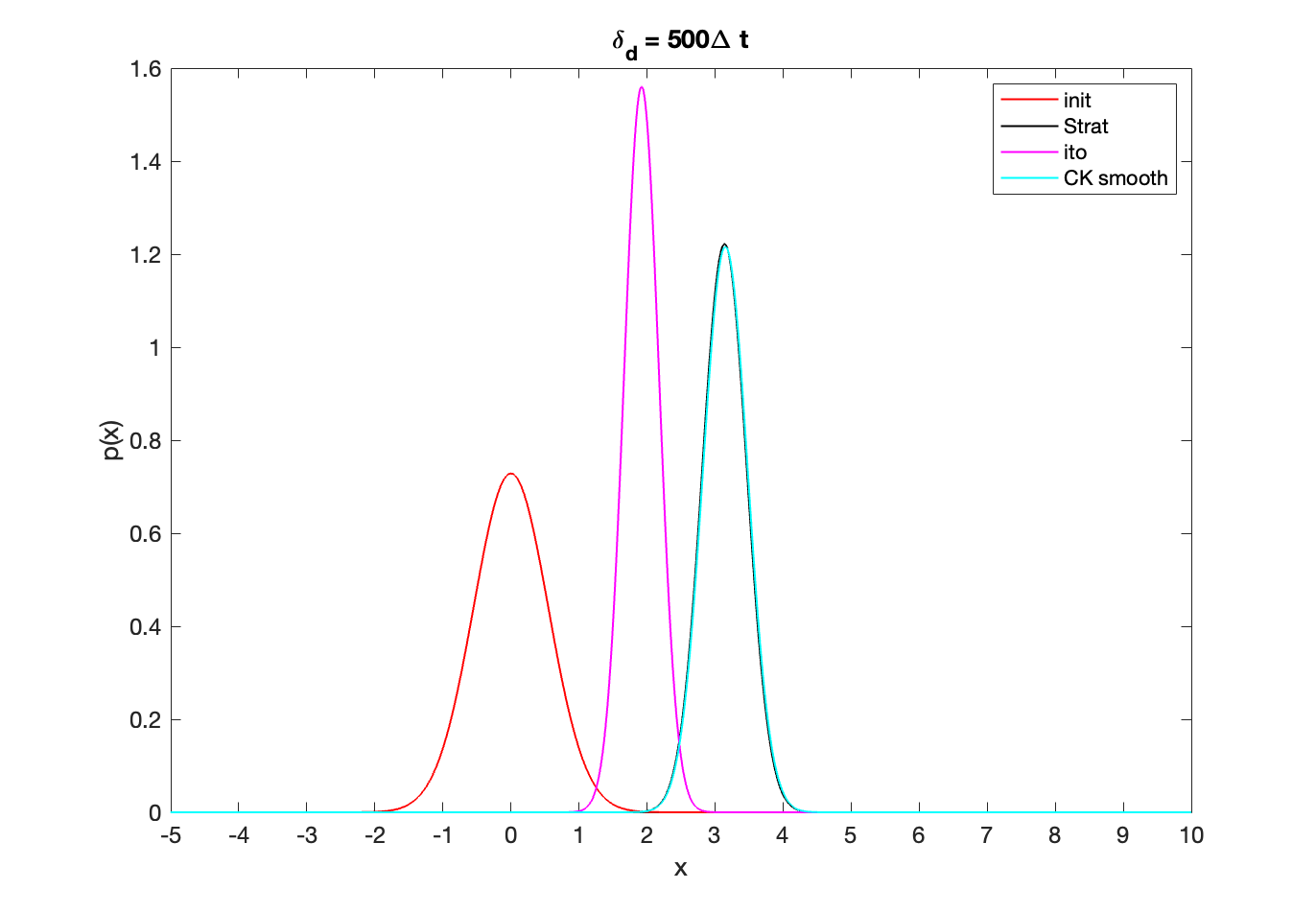}
		\caption{Snapshot in time for the above filtering problem with $H=2, m_0=0, P_0 = 0.3, x_0^\ast = 5, \Xi = 1$ and the smooth observations are constructed with $\delta_d = 500 \Delta t$.  Clearly the Ito interpretation (pink line) is not the correct limit for the Crow-Kimura with smooth approximation.  The correct Stratonovich interpretation (black line) aligns closely with the Crow-Kimura with smooth obs (cyan line). Note that the Stratonovich interpretation \eqref{eq:zakaistratlingauss} coincides with the familar Zakai equation from filtering \eqref{eq:zakai}.} \label{fig:itovsstrat}
		\label{fig:itostratexample}
	\end{figure}
	
\end{example}

The following theorem establishes the convergence of the replicator-mutator equation (in the form \eqref{eq:repasKS}, as identified in Lemma \ref{lem:ckreform}) to a ``generalised'' form of the Kushner-Stratonovich equation from non-linear filtering as $d \rightarrow \infty$.  Convergence is studied via the unnormalised equations as this greatly simplifies the analysis but still yields the overall conclusion relating filtering and replicator-mutator equations due to the one to one correspondence between the unnormalised and normalised equations. Convergence to the standard filtering equations for the specific choice $r=1,s=0$; the benefits of the generalised form (ie. when $s \neq 0, r\neq 1$) will be further explored in the context of misspecified filtering in Section \ref{sec:misspec}.  
The proof of Theorem \ref{theo:limitcrowkimura} borrows many elements from the proof of Theorem 3.1 in \cite{hu_approximation_2002} and makes use of the representation formulae developed in \cite{kunita_stochastic_1982}.  We extend their work to consider unbounded observation drifts $h$ (where they had assumed uniform boundedness) and to the case of multivariate rather than scalar valued observations $Z_t$.  Due to the representation formula used here, we do not need to rely on strong convergence of piecewise smooth approximations with unbounded diffusion coefficients as developed in e.g. \cite{pathiraja_l_2024} (this aspect is discussed more specifically in the proof below).  The price paid is that we focus on pointwise convergence of the density functions, rather than stronger $L^p$ convergence, i.e. (i.e. $\mathbb{E}[\| \mu_t^d - q_t(x)\|_p^p] \rightarrow 0$ as $d \rightarrow \infty$ where $ \| f \|_p^p := \int_{\mathbb{R}^d} |f(x)|^p dx $.  The weaker mode of convergence is still useful, particularly given that it allows us to relax restrictive assumptions on $h$ which previously did not even cover the linear-Gaussian setting.  The proof of the following theorem can be found in section \ref{sec:prooftheolimitcrowkimura}. 

\begin{theorem}
	\label{theo:limitcrowkimura} 
	Assume that $g(x): \mathbb{R}^m \rightarrow \mathbb{R}^m$ and $h(x): \mathbb{R}^m \rightarrow \mathbb{R}^n$ are $C^2$, globally Lipschitz continuous functions satisfying linear growth conditions, i.e. there exists a constant $C > 0$ such that 
	\begin{align*}
		\| h(x) - h(y) \| + \|g(x) - g(y) \| &\leq C\|x - y\|, \enskip x, y, \in \mathbb{R}^m \\
		\| h(x)\| + \|g(x) \| &\leq C(1 + \|x \|), \enskip x \in \mathbb{R}^m
	\end{align*}
	Let $\Sigma$ and $\Xi$ be $m \times m$ and $n \times n$ positive definite matrices respectively. 
	Denote by $\mu_t^d(x)$ the solution to the unnormalised Crow-Kimura replicator-mutator equation with time-varying fitness landscape \eqref{eq:genquadfit} with $s < r, r >0$
	\begin{align}
		\nonumber 
		\partial_t \mu_t^d(x) &= \mathcal{L}^\ast \mu_t^d(x) +  \mu_t^d(x) \mathbb{E}_{z \sim \mu_t^d}[f_{t}(x,z)] \\
		\label{eq:repmutnonlinunnorm}
		& =  \mathcal{L}^\ast \mu_t^d(x) -\frac{r}{2}\|h(x) - \frac{dZ_t^d}{dt}\|_{\Xi}^2 \mu_t^d(x) + s \left \langle h(x) - \frac{dZ_t^d}{dt}, h(z) - \frac{dZ_t^d}{dt} \right \rangle_{\Xi} \mu_t^d(x)
	\end{align}
	Let $q_t(x)$ denote the solution of the modified Zakai equation (presented here in Ito form), 
	\begin{align}
		\label{eq:infzakai}
		d q_t = \mathcal{L}^\ast q_t(x) -\frac{s}{2}h(x)^\top \Xi^{-1} h(x) q_t(x)  + (r-s)q_t(x) h(x)^\top \Xi^{-1} dZ_t. 
	\end{align}
	Suppose also that $q_0(x) = \mu_0^d(x) = f(x)$ where $f$ is a uniformly bounded $C^\infty$ probability density function.  Then for any $T > 0$, 
	\begin{align}
		\lim_{d \rightarrow \infty} \mathbb{E} \left[ \sup_{0 \leq t \leq T} |\mu_t^d(x) - q_t(x) |^p\right] = 0, \quad \forall \enskip x \in \mathbb{R}^m 
	\end{align}
	for $p \geq 1$ and $r,s$ additionally satisfying $r-s < \frac{2}{ tC(1+ \mathbb{E}^\mathbb{Q}[|x|^2])p \lambda_\Xi^2 r_1r_2} $, where $r_1,r_2 >1$ and $1/r_1 + 1/r_2 = 1$ and $C$ depends on the linear growth constant of $h(x)$. 
	Importantly, when $r=1, s = 0$, \eqref{eq:repmutnonlin} converges to the standard Zakai equation \eqref{eq:zakai}.
	
\end{theorem}

\section{Replicator-mutator equations \& filtering with misspecified models: the Linear-Gaussian case}
\label{sec:misspec}

We now focus on the linear-Gaussian setting to demonstrate both analytically and numerically the benefits of the non-local replicator-mutator model for inference in the presence of model misspecification.  Throughout this section, we use the terminology \textit{local} and  \textit{non-local} replicator-mutator equation to refer to \eqref{eq:repasKS} with $s=0, r > 0$ and $s \neq 0, r >0$, respectively.  The terminology is motivated by the fact that the case $s \neq 0$ introduces a non-local or interaction-type term into the fitness function.  
Before detailing the main insights, we first present the fundamental equations in the linear-Gaussian setting and establish some useful findings on mean-field interacting particle systems corresponding to \eqref{eq:repasKS}. 
Firstly, consider the (normalised) linear-Gaussian Crow-Kimura replicator-mutator equation in the limit $\delta_d \rightarrow 0$ as derived in Theorem \ref{theo:limitcrowkimura} i.e. the normalised form of \eqref{eq:infzakai} with $h(x) = Hx, \enskip H \in \mathbb{R}^{n \times m}, \enskip g(x) = Gx,  \enskip G \in \mathbb{R}^{m \times m}, \enskip p_0 = \mathcal{N}(x; m_0,C_0)$,
\begin{align}
	\nonumber 
	dp_t(x) = \mathcal{L}^\ast p_t(x)dt \, -\,  &(r-s) p_t(x) (Hx - Hm_t)^\top\Xi^{-1}Hm_t dt + (r-s) p_t(x)(Hx - Hm_t)^\top\Xi^{-1}dZ_t \\
	\label{eq:modifiedkushner}
	&-\frac{s}{2} p_t(x) \left( \| Hx - Hm_t\|_\Xi^2  - \mathbb{E}_{p_t}[\| Hx - Hm_t\|_\Xi^2] \right)dt. 
\end{align}
The form with $\delta_d > 0$ is given by 
\begin{align}
	\nonumber 
	d p_t^d(x) = \mathcal{L}^\ast p_t^d(x)dt +  (r-s) p_t^d(x) &\left( - \frac12\left( x^\top H^\top \Xi^{-1} Hx   - \mathbb{E}_{p_t^d} \left[ x^\top H^\top \Xi^{-1} Hx \right]\right) + \left(Hx - Hm_t \right)^\top\Xi^{-1} \frac{dZ_t^d}{dt} \right)\\
	\label{eq:repmutlineargauss}
	-\frac{s}{2} p_t^d(x) &\left( \| Hx - Hm_t\|_\Xi^2  - \mathbb{E}_{p_t^d}[\| Hx - Hm_t\|_\Xi^2] \right) dt
\end{align}
where $m_t:= \mathbb{E}_{p_t}[x]$ and $C_t = \mathbb{E}_{p_t}[x x^\top] - \mathbb{E}_{p_t}[x]\mathbb{E}_{p_t}[x^\top]$ for the mean and covariance respectively.  It can be shown straightforwardly (see e.g. \cite{pathiraja_mathematical_2024}) that the evolution of the mean and covariance of the Crow-Kimura replicator-mutator equation \eqref{eq:modifiedkushner} is given by 
\begin{align}
	\label{eq:meanlinear}
	dm_t &= Gm_tdt \, + \, (r-s )K_t\left(dZ_t - Hm_t dt\right)\\
	\label{eq:covlinear}
	\frac{dC_t}{dt} &= GC_t \, + \, C_tG^\top \, + \, \Sigma -r C_t H^\top \Xi^{-1} H C_t 
\end{align}
where $K_t$ denotes the continuous time Kalman gain matrix, given by
\begin{align}
    \label{eq:kalmangain}
    K_t := C_t H^\top \Xi^{-1}
\end{align}
and also that the mean equation for the case $\delta_d >0$, i.e. \eqref{eq:repmutlineargauss} is given by
\begin{align}
	\label{eq:meanlinearsmoothobs}
	\frac{dm}{dt} &= Gm_t \, + \, (r-s )K_t\left(\frac{dZ_t^d}{dt} - Hm_t\right)
\end{align}
and the covariance equation coincides with \eqref{eq:covlinear}.  Notice that \eqref{eq:meanlinear} \& \eqref{eq:covlinear} coincide with the familiar Kalman-Bucy filter equations for $r=1, \, s=0$.  Also, while both $s$ and $r$ affect the evolution of the mean, only $r$ affects the evolution of the covariance.  This aspect will be discussed further in the context of misspecified filtering in section \ref{sec:biasvar}.  We discuss connections between the Crow-Kimura replicator-mutator equations and covariance inflated ensemble Kalman-Bucy filtering in the next section.  A nudging type filter similar to \eqref{eq:meanlinear}-\eqref{eq:covlinear} in the Linear-Gaussian case has appeared at around the same time as this work (see section 3.5 of \cite{gonzalez_nudging_2024}). 

\subsection{Replicator-mutator \& covariance inflated Kalman-Bucy filtering}
\label{sec:inflrepmut}

Though the replicator-mutator equations describe the macroscopic evolution of traits in a population, it is possible to derive microscopic or interacting particle based systems whose dynamics in the limit of infinite population sizes is described by \eqref{eq:modifiedkushner} and \eqref{eq:repmutlineargauss}.  Such mean-field processes arise in interacting particle implementations of Kalman Bucy filtering (the so-called ensemble Kalman-Bucy methods, in both the stochastic form \cite{VanLeeuwen2020,Houtekamer1997} and deterministic forms \cite{Bergemann2012, Taghvaei2017, bishop_mathematical_2023}), given by 
\begin{align}
\label{eq:stochenkbf}
d\bar{X}_t &= G \bar{X}_t + \Sigma^{1/2} dW_t + K_t \left( dZ_t - H \bar{X}_tdt + \sqrt{2} \, \Xi^{1/2}d\bar{B}_t \right) \quad \text{(stochastic enKBf)} \\
\label{eq:detenkbf}
d\bar{X}_t &= G \bar{X}_t dt + \Sigma^{1/2} dW_t + K_t \left(dZ_t - H\left( \frac{\bar{X}_t + m_t}{2}  \right)dt \right) \quad \text{(deterministic enKBf)}
\end{align}
where $\bar{B}$ is a Wiener process independent of $W,B$.  The variable $\bar{X}_t$ is a random variable representing the trait at time $t$ of an infinitely large population of individuals.  These are interacting processes due to the presence of $C_t$ in the driving equations.  It is well-known that the conditional probability density of $\bar{X}_t$ given $\mathcal{Z}_t$ (denoting the frequency of traits in the population) for both \eqref{eq:stochenkbf} and \eqref{eq:detenkbf} is a solution of \eqref{eq:modifiedkushner} with $r=1, s=0$.  These processes are fundamental in the field of particle based sequential Bayesian inference as they are more computationally efficient to simulate than discretising the corresponding PDEs.  The next lemma establishes mean-field processes more generally for \eqref{eq:modifiedkushner} and \eqref{eq:repmutlineargauss} with  $r >0, s < r$.  The proof of the lemma can be found in section \ref{sec:proofmeanfieldproc}.

\begin{lemma}
	\label{lem:meanfieldproc}
    Consider the multivariate linear-Gaussian setting, $g(x) := Gx, \enskip h(x):= Hx$ and a mean-field SDE of the form 
    \begin{align}
        \label{eq:basemf}
        d\bar{X}_t = \underbrace{G\bar{X}_t + \Sigma^{1/2}dW_t}_{\text{mutation}}  \quad + \underbrace{dI_t}_{\text{observation update}} 
    \end{align}
    with $\bar{X}_0 \sim \mathcal{N}(x; m_0, C_0)$ and the stochastic observation update term 
    \begin{align}
    \label{eq:mfinfla}
        dI_t = &-\frac{s}{2}K_t H (\bar{X}_t - m_t)dt + \underbrace{ (r-s) K_t \left(dZ_t - H\bar{X}_t dt  + \left(\frac{2}{r-s}\Xi\right)^{1/2} d\bar{B}_t\right)  }_{\text{Stochastic Kalman Innovation}} 
    \end{align}
    or the deterministic observation update term 
    \begin{align}
	\label{eq:mfinflabias}
	d I_t = \left( \frac{r}{2} - s \right) K_t H (\bar{X}_t - m_t)dt   + \underbrace{\enskip  (r-s) K_t\left( dZ_t - \frac{1}{2} \left(H\bar{X}_t + Hm_t  \right) dt \right)}_{\text{Deterministic Kalman Innovation}} 
\end{align}
where $\bar{B}$ is a scalar Wiener process independent of $W, B$. The time evolution of the conditional density of $\bar{X}_t$ given $\mathcal{Z}_t$ for \eqref{eq:basemf} with either \eqref{eq:mfinfla} or \eqref{eq:mfinflabias} observation update terms is given by the limiting non-local replicator-mutator equation or modified Kushner equation \eqref{eq:modifiedkushner}.  Similarly, consider observation update terms driven by piecewise smooth in time observations, 
\begin{align}
	\label{eq:mfdetode} 
	d I_t =    - \frac{s}{2} K_t H \left( \bar{X}_t - m_t \right)dt + \enskip (r-s) K_t \left(\frac{dZ_t^d}{dt} - \frac{1}{2}H\left( \bar{X}_t + m_t  \right)\right)dt 
\end{align}
or the stochastic version 
\begin{align}
	\label{eq:mfstochode}
	d I_t = -\frac{s}{2}K_t H (\bar{X}_t - m_t)dt + (r-s) K_t \left(\frac{dZ_t^d}{dt} - H\bar{X}_t  \right)dt + \sqrt{r-s} K_t \Xi^{1/2} d \bar{B}_t  
\end{align}
The time evolution of the conditional density of $\bar{X}_t$ given $\mathcal{Z}_t$ for \eqref{eq:basemf} with either \eqref{eq:mfstochode} and \eqref{eq:mfdetode} observation update terms is given by the non-local replicator-mutator equation \eqref{eq:repmutlineargauss}.  
\end{lemma}

Lemma \ref{lem:meanfieldproc} allows us to immediately relate the replicator-mutator equations with non-local fitness functions to the so-called covariance inflation methods in ensemble Kalman filtering \cite{anderson_adaptive_2007, mitchell_accounting_2015, hamill_accounting_2005, Mitchell2000, tong_nonlinear_2016, bishop_mathematical_2023}.  Covariance inflation is an important heuristic tool used to improve numerical stability of the ensemble Kalman filter when the number of particles is low \cite{bishop_perturbations_2018, Bishopstable2019} and also to account for model errors.  There are two widely used forms of covariance inflation, the first being multiplicative covariance inflation,  see e.g. \cite{bishop_perturbations_2018,bishop_mathematical_2023} which when applied to continuous time stochastic ensemble Kalman-Bucy filtering takes the form
\begin{align}
\label{eq:inflateenkbf}
d\bar{X}_t = G \bar{X}_t + \Sigma^{1/2} dW_t + (C_t + \epsilon T)H^\top \Xi^{-1} \left( dZ_t - H \bar{X}_tdt + \Xi^{1/2}d\bar{B}_t \right) 
\end{align}
where $\epsilon >0$ is a tuning parameter and $T \in \mathbb{R}^{m \times m}$ is a reference matrix guiding the inflation (most commonly, $T = C_t$).  A similar mean field ensemble Kalman-Bucy filter with the second form of inflation, so-called additive inflation in the spirit of \cite{hamill_accounting_2005} (see (7) in their paper), takes the form  
\begin{align}
\label{eq:inflatenkbfadd}
d\bar{X}_t = G \bar{X}_t + \Sigma^{1/2} dW_t + \epsilon T(\bar{X}_t - m_t)dt   + C_tH^\top \Xi^{-1} \left( dZ_t - H \bar{X}_tdt + \Xi^{1/2}d\bar{B}_t \right) 
\end{align}
where $\epsilon, T$ are as previously defined although most commonly $T = I$ is adopted in this form. Such inflation terms can similarly be applied to the deterministic ensemble Kalman-Bucy filter \eqref{eq:detenkbf}. It follows directly from Lemma \ref{lem:meanfieldproc} that choosing $r = 1+ \epsilon, \enskip s=0$ in \eqref{eq:modifiedkushner} yields the multiplicative covariance inflated stochastic enKBf \eqref{eq:inflateenkbf} with $T = C_t$.  Likewise, choosing $s = -2\epsilon, \enskip r = 1 -2\epsilon$ yields the additive covariance inflation stochastic enKBf \eqref{eq:inflatenkbfadd} with $T = K_tH$.  Connections similarly hold for the deterministic versions and also with piecewise smooth observations.  This implies that for an arbitrary choice of $r >0, s < r$, the Crow-Kimura replicator-mutator equation (and its limiting form) with non-local fitness function involves a combination of additive and multiplicative inflation.  We explore this observation in more detail in the next section.     

\subsection{Case study: non-local replicator-mutator vs inflation for misspecified model filtering}
\label{sec:biasvar}

 It is known empirically in the data assimilation literature that a combined used of additive and multiplicative inflation improves filtering estimates in the presence of model and finite-sample effects errors than using either additive or multiplicative inflation alone \cite{scheffler_dynamical_2022}.  This section focuses on  proving this empirical observation for the case of linear-Gaussian filtering in the presence of a constant bias term.  
 
It is well known that the Kalman-Bucy filter is the minimum variance unbiased estimator of $X_t, \enskip t > 0$.  This property also extends to the non-local replicator mutator equation when the system is perfectly known (see Lemma \ref{lem:unbiasedness}).  However, from inspection of \eqref{eq:covlinear} and also from Lemma \ref{lem:meanfieldproc}, it is clear that the minimum variance property is destroyed unless $r=1$.

We now turn out attention to analysing the performance of the replicator-mutator equation for a filtering problem with misspecified signal dynamics. This section focuses on the $\delta_d \rightarrow 0$ form of the replicator-mutator, \eqref{eq:modifiedkushner}, although the conclusions of this section are expected to similarly hold for $\delta_d >0$.   
The misspecified model filtering problem is as follows. Consider the following linear-gaussian problem where the hidden state $X_t \in \mathbb{R}^m$ evolves according to 
\begin{align}
\label{eq:realhidden}
dX_t = GX_tdt + bdt + \Sigma^{1/2}dW_t 
\end{align}
with $X_0 \sim \mathcal{N}(m_0, P_0)$ and $b$ a constant vector.  The hidden dynamics is known imperfectly and that the assumed model for the trait is instead
\begin{align}
\label{eq:asshidden}
dX_t = GX_tdt + \Sigma^{1/2}dW_t 
\end{align}
(i.e. \eqref{eq:realhidden} with $b=0$), and the observation process is as in \eqref{eq:obs}.  The replicator-mutator equation \eqref{eq:repmutlineargauss} is no longer expected to produce unbiased estimates due to the presence of the $b$ term.  The tracking error at time $t$ is denoted by $\varepsilon_t$, 
\begin{align*}
\varepsilon_t := m_t - x_t^\ast
\end{align*}
where $m_t := \int x p_t(x)dx$ and $p_t(x)$ is the solution of \eqref{eq:modifiedkushner} and $x_t^\ast$ the solution of \eqref{eq:realhidden} that generated the observation path, also known as the reference trajectory or true hidden state.  Denote also by $P_t$ the error covariance matrix,
\begin{align*}
P_t:=\mathbb{E}[(\varepsilon_t - \mathbb{E}[\varepsilon_t])(\varepsilon_t - \mathbb{E}[\varepsilon_t])^T]
\end{align*}
and 
\begin{align}
\label{eq:defnPtilde}
\tilde{P}_t := \mathbb{E}[ \varepsilon_t \varepsilon_t^\top].
\end{align}
Note well that $P_t$ is now distinct from $C_t$, the posterior/filtering covariance or covariance of the replicator-mutator equations.  Recall that in the perfect knowledge Kalman Bucy filtering setting (i.e. $b=0$), since $\mathbb{E}[\varepsilon_t] = 0$, we have that $C_t = \tilde{P}_t = P_t$.  When $b \neq 0$, even a standard Kalman-Bucy filter no longer produces an unbiased estimates of the hidden state, so that $\tilde{P}_t \neq P_t$ and also, $C_t \neq P_t$.  Finally, let 
\begin{align}
\label{eq:defnnut}
\nu_t &:= \| \mathbb{E}[\varepsilon_t] \|^2 = \text{Tr}[\mathbb{E}[\varepsilon_t]\mathbb{E}[\varepsilon_t^\top]] \\
\label{eq:defnEt}
\quad E_t &:= \mathbb{E}[ \| \varepsilon_t \|^2] =  \text{Tr}[\tilde{P}_t]
\end{align}
denote the squared expected error (or squared bias) and expected squared error (mean squared error, MSE), respectively.  Additionally, the standard bias-variance decomposition holds
\begin{align*}
\underbrace{E_t}_{\text{MSE}} = \underbrace{Tr(P_t)}_{\text{total error var.}} + \underbrace{\nu_t}_{\text{squared bias}}  
\end{align*}

\subsubsection{Mean squared error analysis}
\label{sec:mseanalysis}

Here we analytically determine values of $r,s$ at which the minimum MSE is obtained and evaluate the corresponding posterior covariance $C_t$.  The following quantity will be used throughout the rest of the section, 
\begin{align}
\label{eq:defnAt}
A_t(r,s) := G - (r-s)K_t H = G - (r-s)C_t H^\top \Xi^{-1} H.
\end{align}
It is well-known that the stability of the Kalman-Bucy filter ($r=1,s=0$) is intimately tied to the spectral properties of this quantity $A_t$.  In particular, we rely on the following assumptions for the remainder of this section for $r >0, s < r$ more generally. 
\begin{aspt}
\label{ass:Astable}
        The noise covariance matrices $\Sigma$ and $\Xi$ are positive definite square matrices. Furthermore, the signal drift and observation drift operators satisfy observability and controllability conditions, i.e. 
        \begin{align*}
         \left[   \begin{array}{cc}
                 \Sigma^{1/2}, G \Sigma ^{1/2}, \dots, G^{m-1} \Sigma^{1/2}&  
            \end{array} \right], \quad \text{and} \quad \left[\begin{array}{cc}
                H &  \\
                HG & \\
                \vdots \\
                HG^{m-1}
            \end{array} \right]
        \end{align*}
        have rank $m$.


\end{aspt}
Assumption \ref{ass:Astable} guarantees the existence of a unique $C_\infty$.  Note that in the case of Kalman filtering, assumption \ref{ass:Astable} is also sufficient to ensure $\alpha(A_\infty) < 0$ when $s= 0$ (see e.g. \cite{bishop_mathematical_2023}).  When $s \neq 0$, the presence of the $(r-s)$ term in \eqref{eq:defnAt} no longer guarantees this, therefore we rely on the following additional assumption. 

\begin{aspt}
\label{ass:eigsymmA}
    For $A_t$ as defined in \eqref{eq:defnAt}, the parameters $r,s$ are chosen such that $\alpha(A_\infty) < 0$ where $\alpha(A)$ is the spectral abscissa of an $\mathbb{R}^{m \times m}$ matrix $A$, i.e. $\max_i \{ \text{Re}(\lambda_i) \}$ where $\lambda_i$ is the $i$th eigenvalue of $A$.  
    
    Additionally, they are specified such that $\alpha(A_\infty (r,s) + A_\infty(r,s)^\top) < 0$, i.e. the largest eigenvalue of $A_\infty + A_\infty^\top$ is strictly negative. 
\end{aspt}

The following lemma characterises the time evolution of bias, variance and mean squared error.  As expected, the error variance $P_t$ evolves independently of the unknown term $b$.  The proof of the lemma can be found in section \ref{sec:proofcovreperror}

\begin{lemma}
\label{lem:covreperror}
Assume the system properties described by \eqref{eq:realhidden} for the true hidden state, \eqref{eq:asshidden} for the assumed hidden state and \eqref{eq:obs} for the observation model.  Given $r,s$ such that $r <s, r >0$, we have the following evolution equations for the  error covariance $P_t$ and expected squared error $\tilde{P}_t$,
\begin{align}
	\label{eq:odeforPt}
	\frac{dP_t}{dt} &= (G -(r-s)K_tH) P_t  + P_t (G -(r-s)K_tH)^\top  + \Sigma   + (r-s)^2 K_t \Xi K_t^\top   \\
	\label{eq:odeptilde}
	\frac{d\tilde{P}_t}{dt} &= (G -(r-s)K_tH)\tilde{P}_t  + \tilde{P}_t(G -(r-s)K_tH)^\top  + \Sigma  + (r-s)^2 K_t \Xi K_t^\top   - \mathbb{E}[\varepsilon_t]b^\top  - b \mathbb{E}[\varepsilon_t^\top]  
\end{align}
where $K_t$ is the continuous time Kalman gain matrix as defined in \eqref{eq:kalmangain}.  For the special case $s=0,r=1$, it holds that $P_t = C_t, \enskip \forall \enskip t > 0$ if $P_0 = C_0$ where $C_t$ is the solution of the covariance equation of the replicator-mutator.\\
\\
Furthermore, for any p.d. $C_0$,  the evolution of the mean squared error $E_t$ satisfies the following inequality 
\begin{align}
	\label{eq:odeforEt}
	\frac{dE_t}{dt} \leq \alpha(A_t(r,s) + A_t(r,s)^\top)E_t - 2\text{Tr}[\mathbb{E}[\varepsilon_t]b^\top] + \text{Tr}[\Sigma] + (r-s)^2 \alpha(H^\top \Xi^{-1}H) \|C_t\|_F
\end{align}
where $A_t$ is defined in \eqref{eq:defnAt} and $\alpha(A)$ denotes the largest eigenvalue of $A$. 
\end{lemma}

To obtain further insights on optimal choices of $r,s$, we consider a simplified setting where $C_0 = C_\infty$ (i.e. where the covariance is initialised at the steady state covariance matrix in \eqref{eq:covlinear}).  This setting is still rich enough to provide insights on the role of $r,s$ in the non-local replicator-mutator, particularly as we are primarily interested in the time asymptotic behaviour of mean squared error. The following lemma gives explicit representations of the time asymptotic squared bias $\nu_\infty$ and mean squared error $E_\infty$ in terms of the system parameters. The proof of the lemma can be found in section \ref{sec:proofasympbiasvar}.


\begin{lemma}
\label{lem:asympbiasvar}
\textbf{Steady state Bias-Variance}. Suppose assumptions \ref{ass:Astable} and \ref{ass:eigsymmA} holds along with $\mathbb{E}[\varepsilon_0] = 0$ and $C_0 = C_\infty$ where $C_\infty$ satisfies
\begin{align}
	\label{eq:steadystatecov}
	0 = GC_\infty + C_\infty G^\top + \Sigma -rC_\infty H^\top \Xi^{-1} H C_\infty.
\end{align}
For $\nu_t, E_t$ as defined in \eqref{eq:defnnut} and \eqref{eq:defnEt} respectively,  $\nu_t \rightarrow \nu_\infty$ and $E_t \rightarrow E_\infty$ as $t \rightarrow \infty$, where
\begin{align}
	\label{eq:asympbias}
	\nu_\infty &= \| A_\infty^{-1} b \|^2 \\
	\label{eq:asympmerror}
	E_\infty &=  \text{Tr}\left[ (\Sigma + (r-s) \left( G - A_\infty \right)C_\infty  - 2A_\infty^{-1} b b^\top) X_\infty \right] \\
    \label{eq:boundEinf}
    & \leq -\frac{1}{\alpha(A_\infty + A_\infty^\top)} \left( - 2 \text{Tr}[A_\infty^{-1}b b^\top] + \text{Tr}[\Sigma] + (r-s)^2 \alpha(H^\top \Xi^{-1}H) \|C_\infty\|_F\right)
\end{align}         
where $X_\infty$ is the unique solution of the Lyapunov equation  
\begin{align*}
	A_\infty^\top X_\infty + X_\infty A_\infty + I = 0.  
\end{align*}
In the scalar case $m=n=1$, $\nu_\infty$ and $E_\infty$ can be characterised entirely in terms of the system parameters $G, H, \Sigma, \Xi, b$, 
\begin{align}
	\label{eq:asympbiasscalar}
	\nu_\infty &= \left(\frac{b}{A_\infty(r,s)} \right)^2 \\
	\label{eq:asympmsescalar}
	E_\infty = \tilde{P}_\infty &=  - \frac{1}{2}\left(\Sigma + \left( \frac{G - A_\infty(r,s)}{H} \right)^2 \Xi  \right)  \frac{1}{A_\infty(r,s)}  + \left(\frac{b}{A_\infty(r,s)} \right)^2
\end{align}
where 
\begin{align}
	\label{eq:Ainfscalar}
	A_\infty(r,s) &=  \frac{s}{r} G + \frac{(r-s)}{r} A_\infty(r,0) \\
	A_\infty(r,0) &= -\sqrt{G^2 + rH^2 \Xi^{-1} \Sigma}.
\end{align}


\end{lemma}

For the sake of demonstration, the remainder of this section focuses on the fully scalar setting $m=1,n=1$ to characterise the optimal $r,s$ values minimising MSE and leave the multivariate setting to future work.  We expect many of the qualitative insights from the remaining lemmas to hold in the multivariate setting.  The next lemma gives explicit expressions for the optimal $r,s$ minimising time asymptotic MSE entirely in terms of the system parameters; firstly an expression for $s$ for any given $r$  \eqref{eq:optsr} and then for $r$ for given $s$ (see \eqref{eq:optr_cond2} - \eqref{eq:optr_cond1}). The proof of this lemma can be found in section \ref{sec:proofoptimalrs}. 

\begin{lemma}
\label{lem:optimalrs}
\textbf{Optimal $r,s$ minimising MSE.}    Adopt the same conditions as in Lemma  \ref{lem:asympbiasvar} and assume further that $m = n = 1$ (scalar setting).  Given a fixed $r > 0$, the optimal value of $s$ minimising $E_\infty$ is given by
\begin{align}
	\label{eq:optsr}
	s^{\text{opt}} = \frac{r(A_\infty^\ast +\sqrt{G^2 + rH^2 \Xi^{-1} \Sigma})}{G +\sqrt{G^2 + rH^2 \Xi^{-1} \Sigma}}
\end{align}
where $A_\infty^\ast$ denotes the optimal value of $A_\infty$ as defined in \eqref{eq:defnAt},  which has an explicit representation in terms of the system parameters, 
\begin{align}
	\label{eq:Ainfstarexp}
	A_\infty^\ast := 
	\begin{cases}
		\left(-\frac{q}{2} + \sqrt{\tau}\right)^{1/3} + \left(-\frac{q}{2} - \sqrt{\tau}\right)^{1/3}, \enskip \tau > 0\\
		2\,\sqrt{-\frac{p}{3}}\,\cos\left[\,\frac{1}{3} \cos^{-1}\left(\frac{3q}{2p}\sqrt{\frac{-3}{p}}\,\right) - \frac{4\pi}{3}\,\right] , \enskip \tau < 0
	\end{cases}
\end{align}
and $\tau := \frac{q^2}{4} + \frac{p^3}{27}$ with $ p := - (H^2 \Xi^{-1}\Sigma + G^2) $ and $q := 4b^2H^2\Xi^{-1}$.  Furthermore, the admissible values of $s$ guaranteeing $A_\infty < 0$ satisfy 
\begin{align}
	\label{eq:condsons}
	s < \min \left( r, \frac{ \sqrt{G^2 + rH^2 \Xi^{-1} \Sigma}}{G + \sqrt{G^2 + rH^2 \Xi^{-1} \Sigma}} r \right). 
\end{align}
Conversely, given a fixed $s$, there may be either one or two corresponding optimal $r$ values minimising $E_\infty$.  Specifically, given a fixed $s \geq s^u$ where $s^u:= \frac{(|G| - G)(A_\infty^\ast + G)}{H^2 \Xi^{-1} \Sigma}$, there is a unique optimal $r$ value given by 
\begin{align}
	\label{eq:optr_cond2}
	r^{\text{opt}} =  \frac{1}{4H^2 \Xi^{-1} \Sigma} \left( G - A_\infty^\ast + \sqrt{(G - A_\infty^\ast)^2 + 4(G A_\infty^\ast + sH^2 \Xi^{-1} \Sigma)} \right)^2 - \frac{G^2}{H^2 \Xi^{-1} \Sigma}.
\end{align}
Furthermore, if $G - 2|G| \geq A_\infty^\ast$ then given an $s$ such that $s^l < s < s^u$ where $s^l := -\frac{(G + A_\infty^\ast)^2}{4 H^2 \Xi \Sigma} $, there are two optimal $r$ values given by 
\begin{align}
	\label{eq:optr_cond1}
	r^{\text{opt}} = \frac{1}{4H^2 \Xi^{-1} \Sigma} \left( G - A_\infty^\ast \pm \sqrt{(G - A_\infty^\ast)^2 + 4(G A_\infty^\ast + sH^2 \Xi^{-1} \Sigma)} \right)^2 - \frac{G^2}{H^2 \Xi^{-1} \Sigma}
\end{align}

\end{lemma}
 
Importantly, the above lemma shows that there are infinitely many $(r,s)$ pairs satisfying some minimal conditions that minimise MSE $E_\infty$, although the covariance $C_\infty$ will differ due to its explicit dependence on $r$ (see \eqref{eq:asympCscalar}). Notably, these estimates hold true regardless of the stability characteristics of the hidden state and observation dynamics.  Additionally, it holds that for $s \in (s^l, s^u)$ where $s^l, s^u$ are as defined in the above lemma, there exists two possible $r$ values for a given $s$ that will yield the minimal asymptotic mean square error \eqref{eq:optr_cond1}. This is particularly beneficial in terms of allowing for a realistic $C_\infty$, as will be explored further in the following lemma.  Before stating the lemma, we note that larger values of $r$ give rise to smaller $C_\infty$.    
From an inference perspective, smaller $C_\infty$ indicates greater confidence in the estimator, which can be problematic especially when $C_\infty$ is smaller than the minimum covariance achievable in the perfect model setting.  To analyse this phenomenon further, define 
\begin{align}
\label{eq:covperfect}
\hat{C}_\infty := \frac{G + \sqrt{G^2 + H^2 \Xi^{-1} \Sigma}}{ H^2 \Xi^{-1}}
\end{align}
which is the steady state covariance coinciding with the Kalman-Bucy filter (the optimal filter) in the case of perfect knowledge, (i.e. $r=1,s=0,b=0$).  Any choice of $r$ for which $C_\infty < \hat{C}_\infty$ should be avoided when there is model misspecification, as we cannot expect to be more confident than when we have no misspecification.  Recall that in the perfect knowledge setting, we have that $\tilde{P}_\infty = \hat{C}_\infty$, in other words, the asymptotic covariance from the Kalman-Bucy filter coincides with the mean squared error.  The following lemma establishes a relationship between $r,s$ that ensures the specific pair satisfies $C_\infty = \tilde{P}_\infty$ in the misspecified model setting where $b \neq 0$.  Its proof can be found in section \ref{sec:proofcovimpact}.

\begin{lemma}
\label{lem:covimpact}
Adopt the same conditions as in Lemma \ref{lem:optimalrs}.  Then the following results hold: 
\begin{enumerate}
    \item The optimal $(r,s)$ such that $\tilde{P}_\infty = C_\infty$ satisfies \eqref{eq:optsr} and 
	\begin{align}
		\label{eq:rminuscond}
		r^{\text{opt}} - s^{\text{opt}} = \frac{(A_\infty^\ast)^2(G - A_\infty^\ast)}{-0.5A_\infty^\ast \left( H^2 \Xi^{-1} \Sigma + (G - A_\infty^\ast)^2   \right) + b^2H^2 \Xi^{-1}}
	\end{align}
	where $A_\infty^\ast$ is given by \eqref{eq:Ainfstarexp}.
    \item For $s = 0$, it holds that the value of $r$ minimising $E_\infty$, $r_0^{\text{opt}}$, satisfies 
\begin{align}
		\label{eq:ropt0}
		r_0^{\text{opt}} > \frac{2(G^2 + H^2\Xi^{-1}\Sigma)}{H^2\Xi^{-1}\Sigma} > 1
	\end{align}
independently of the unknown $b$.  Furthermore, for the case $\tau > 0$ where $\tau$ is as defined in Lemma \ref{lem:optimalrs} and $G >0$, it holds that 
\begin{align}
    \label{eq:boundCratio}
    \frac{C_\infty(r_0^{\text{opt}})}{\hat{C}_\infty} < \frac{1}{\sqrt{r_0^{\text{opt}}}} \left( \frac{2|G| + \sqrt{H^2\Xi^{-1} \Sigma}}{(G + \sqrt{G^2 + H^2\Xi^{-1} \Sigma})} \right) < 1
\end{align}
where $\hat{C}_\infty$ denotes the asymptotic covariance in the perfect model setting as defined in \eqref{eq:covperfect}. 
\end{enumerate}

\end{lemma}

This result highlights that the non-local replicator-mutator equation, unlike the multiplicative covariance inflated ensemble Kalman filter ($r >0, s = 0$), is capable of simultaneously minimising asymptotic mean squared error and providing realistic uncertainty estimates through $C_\infty$.  In particular, \eqref{eq:boundCratio} demonstrates that when $G > 0$ (unstable dynamics) the optimal $r$ minimising MSE for the multiplicative covariance inflated ensemble Kalman method typically yields a covariance $C_\infty < \hat{C}_\infty$, i.e. it is over-confident.  Adopting $s \neq 0$ is thus crucial to obtain more representative $C_\infty$.  Equivalently, combining additive and multiplicative covariance inflation along with a judicious choice of tuning parameters is necessary to minimise MSE while maintaining realistic uncertainty quantification.     
%
%
%
%
%
%

\subsubsection{Numerical experiments}
\label{Sec:numericsmissp}

The following experiment aims to provide further insights on the role of $s,r$ in \eqref{eq:repmutnonlin} with fitness landscape $f_t(x,z)$ given by \eqref{eq:genquadfit} for the misspecified model setting from the previous section (which we repeat here for convenience).  Although the analysis in the previous section has been done for the limiting case $\delta_d \rightarrow 0$, here we focus on the practically relevant discrete time observation case (i.e. piecewise smooth observations) and show that much of the analysis holds for $\delta_d$ small enough.  That is,  consider a partition of the time interval $[0,T]$, $0 < t_1 < t_2 \dots < t_d = T$ with time-step $t_{i+1} - t_i = \delta_d$.  Synthetic observations of the form
\begin{align}
\label{eq:sec5obs}
\frac{dZ_t^d}{dt} = Hx_{t_i}^\ast(\omega) + \Xi^{1/2} \frac{dB_t^d}{dt}, \quad t \in [t_i, t_{i+1})
\end{align}
(as in \eqref{eq:smoothobsnew}) are constructed, where $x_t^\ast(\omega)$ is a solution for a realisation $\omega$ at time $t$ of 
\begin{align}
\label{eq:sec5asshidden}
dX_t = GX_tdt + bdt + \Sigma^{1/2}dW_t 
\end{align}
with $X_0 \sim \mathcal{N}(m_0, P_0)$.  The process in \eqref{eq:realhidden} describes the evolution of the actual optimal trait and $\frac{dZ_t^d}{dt}$ corresponds to noisy observations of it.  Suppose the assumed model for the trait is instead
\begin{align}
\label{eq:sec5realhidden}
dX_t = GX_tdt + \Sigma^{1/2}dW_t 
\end{align}
(i.e. \eqref{eq:realhidden} with $b=0$), so that the corresponding replicator-mutator equation takes the form 
\begin{align}
\label{eq:repmutmiss}
\partial_t p_t^d(x) = -G\nabla \cdot p_t^d(x) + \frac{1}{2} \nabla \cdot (\Sigma \nabla p_t^d(x)) + p_t^d(x)\left( \mathbb{E}_{z \sim p_t^d}[f_t(x,z)] - \mathbb{E}_{x, z \sim p_t^d}[f_t(x,z)]\right)
\end{align}
with $f_t(x,z)$ given by \eqref{eq:genquadfit}.  
The remainder of this section will focus on the following experimental settings, which have been randomly generated. 
\begin{center}
\begin{tabularx}{0.6\textwidth} { 
		| >{\centering\arraybackslash}X 
		| >{\centering\arraybackslash}X 
	| >{\centering\arraybackslash}X | }
\hline
\textbf{Parameter} & \textbf{System 1} $(\tau > 0)$ & \textbf{System 2} $(\tau < 0)$ \\
\hline
$G$  & 0.5  & 2.5  \\
\hline
$H$  & 8.5  & 2.9  \\
\hline
$\Sigma$  & 0.8  & 18  \\
\hline
$\Xi$  & 6.3  & 26  \\
\hline
$b$  & 9.9  & 1.2  \\
\hline
\end{tabularx}
\end{center}

Note that throughout, we assume $C_0 = C_\infty$ given by \eqref{eq:asympCscalar}.  The settings in System 1 and 2 correspond to the case where $\tau > 0$  (i.e. \eqref{eq:depressedcubic} has one real root) and   $\tau < 0$  (i.e. \eqref{eq:depressedcubic} has three real roots), respectively.  Notice that in both systems, the hidden state evolves according to unstable dynamics and the Crow-Kimura replicator mutator is capable of tracking an unstable signal as $A_t < 0$.  We restrict the time domain to one where machine precision doesn't become an issue.  We adopt $\delta_d = 10^{-3}$ and use a simulation time step of $10^{-4}$ to construct the true hidden state as well as to discretise the mean and covariance equation \eqref{eq:meanlinear} and \eqref{eq:covlinear} using forward euler. We have the following main insights.

\begin{figure}[h!]
\begin{minipage}{\textwidth} 
\centering
\includegraphics[width=0.8\linewidth]{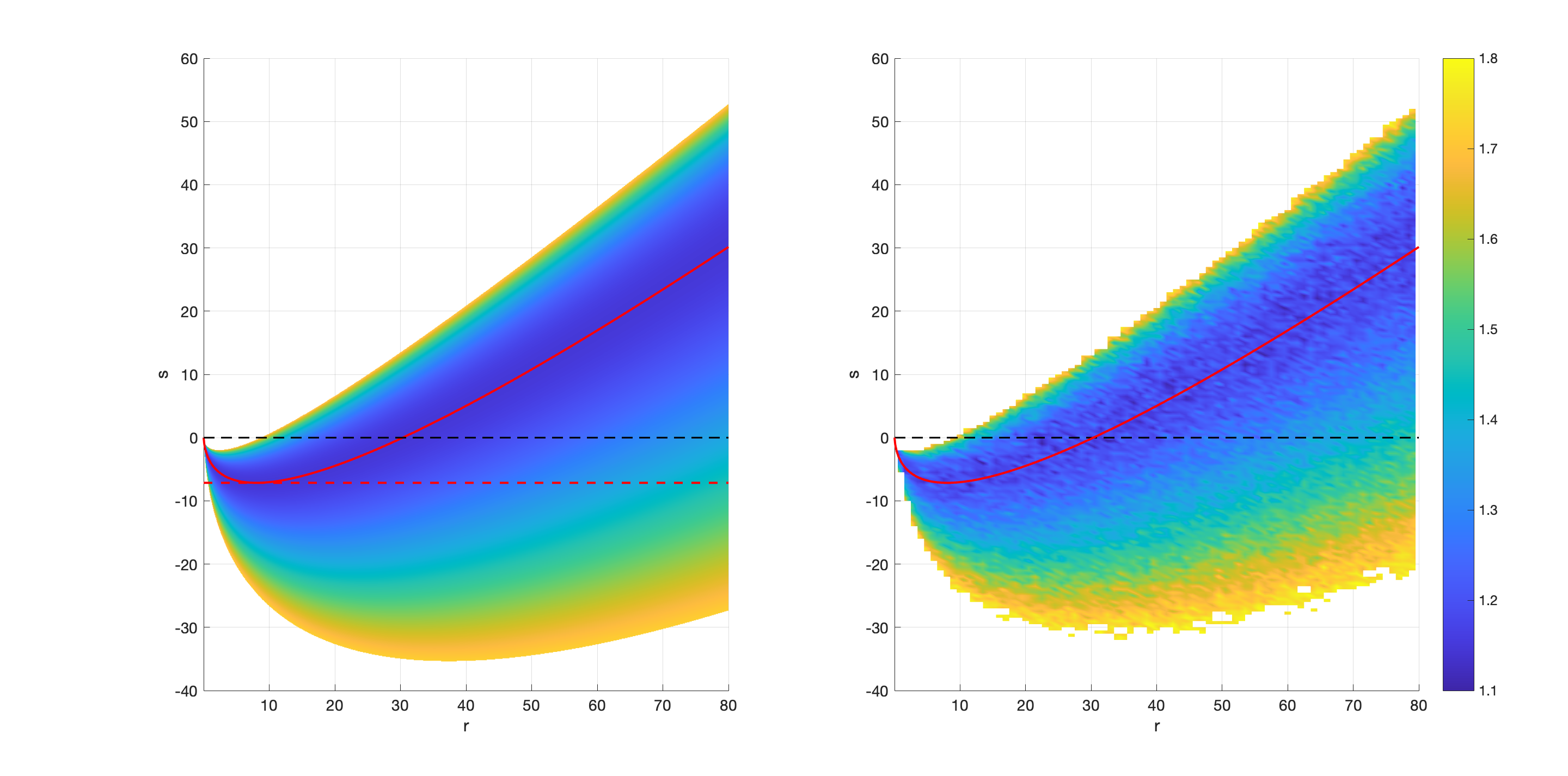}
\caption{Plot of asymptotic MSE $E_\infty$ for various s vs r values for system 2.  The optimal (in terms of asymptotic mse) values are indicated by the red line, calculated using \eqref{eq:optsr}.  The colourbar shows corresponding values of the asymptotic MSE $E_\infty$.  The dashed red line on the left plot shows the theoretical expression for $s^l$ as given in lemma \ref{lem:optimalrs}}. 
\label{fig:analhm1}
\end{minipage}
\vspace{0.5cm}
\begin{minipage}{\textwidth}
\centering
\includegraphics[width=0.8\linewidth]{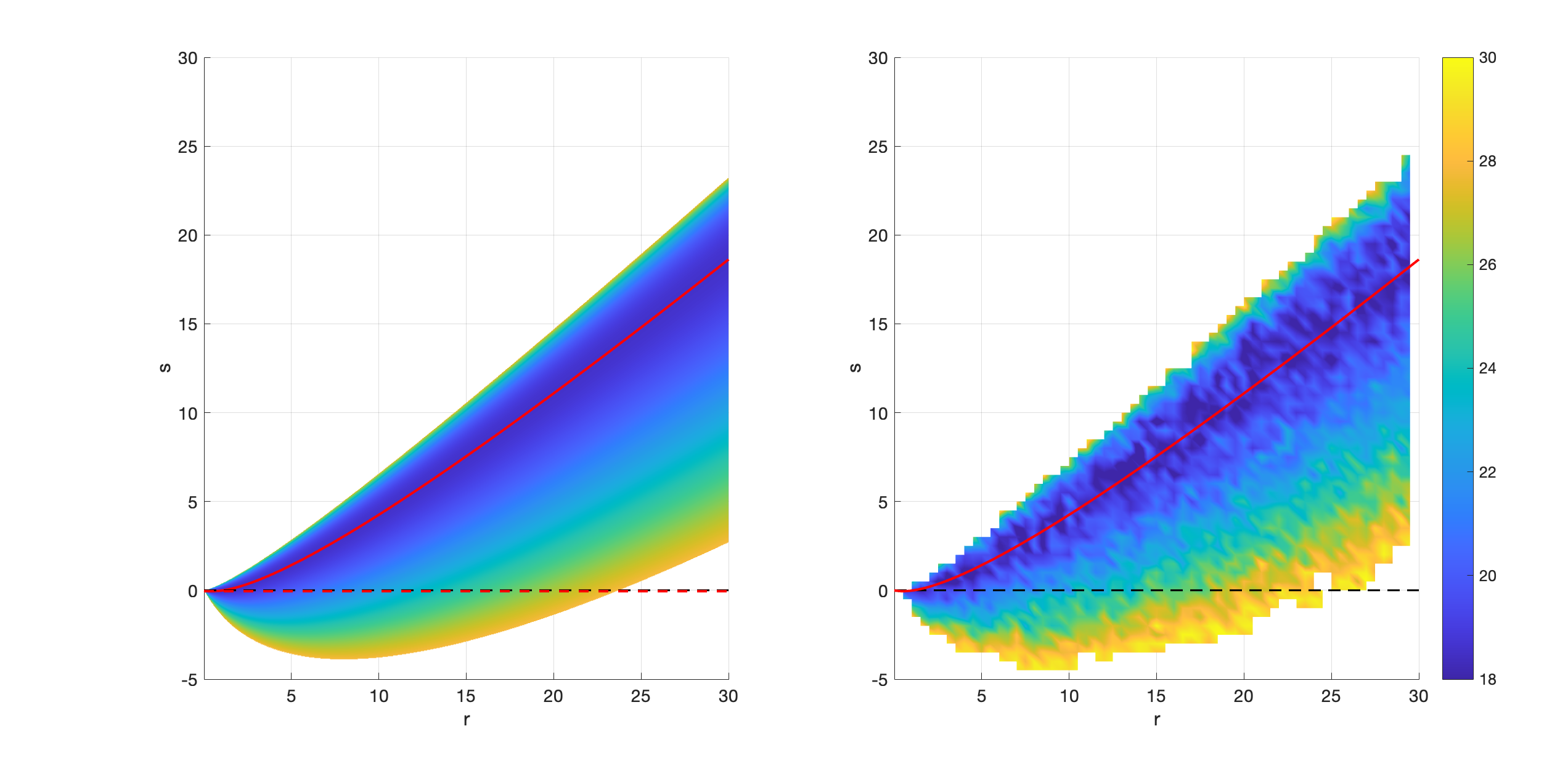}
\caption{Plot of asymptotic MSE $E_\infty$ for various s vs r values for system 2.  The optimal (in terms of asymptotic mse) values are indicated by the red line, calculated using \eqref{eq:optsr}.  The colourbar shows corresponding values of the asymptotic MSE $E_\infty$.}
\label{fig:analhm2}
\end{minipage}
\end{figure}

\paragraph{Verification of optimal $r,s$ (Lemma \ref{lem:optimalrs}).}   To verify the analytic results of Lemma \ref{lem:optimalrs}, we calculate the empirical MSE at time $t$, $\mathcal{E}_t$,
\begin{align}
\mathcal{E}_t =  \frac{1}{N_s} \sum_{j=1}^{N_s} \left( m_t^j - x_t^\ast \right)^2 
\end{align}
where $x_t^\ast$ is a single fixed realisation of \eqref{eq:realhidden}, $m_t^j$ is a solution of \eqref{eq:meanlinearsmoothobs} with $C_0 = C_\infty$ and the index $j$ refers to a single realisation of the smooth observation path $Z_t^d$ and $N_s = 5000$ in these experiments. Figures \ref{fig:analhm1} and \ref{fig:analhm2} show the empirical estimate of the asymptotic mean square error $E_\infty$ for different $(r,s)$ pairs for system 1 and 2 respectively.  They show that in both systems, the analytic expressions in \eqref{eq:optsr} and for the smallest optimal $s$ value, $s^l$ match quite closely.  The existence of two optimal $r$ values for every given $s^l < s < s^u$ is also confirmed numerically (this is more visible for system 1 than in system 2 where $s^l = -0.047$ and in both cases, $s^u \approx 0$).  
Finally, for system 1 in particular it is apparent that choosing $s < 0$ allows to choose smaller $r$ values that can simultaneously reduce mean square error and provide a realistic representation of uncertainty via $C_\infty$. 

\begin{figure}[h!]
	\begin{minipage}{0.45\textwidth} 
		\centering
		\includegraphics[width=1\linewidth]{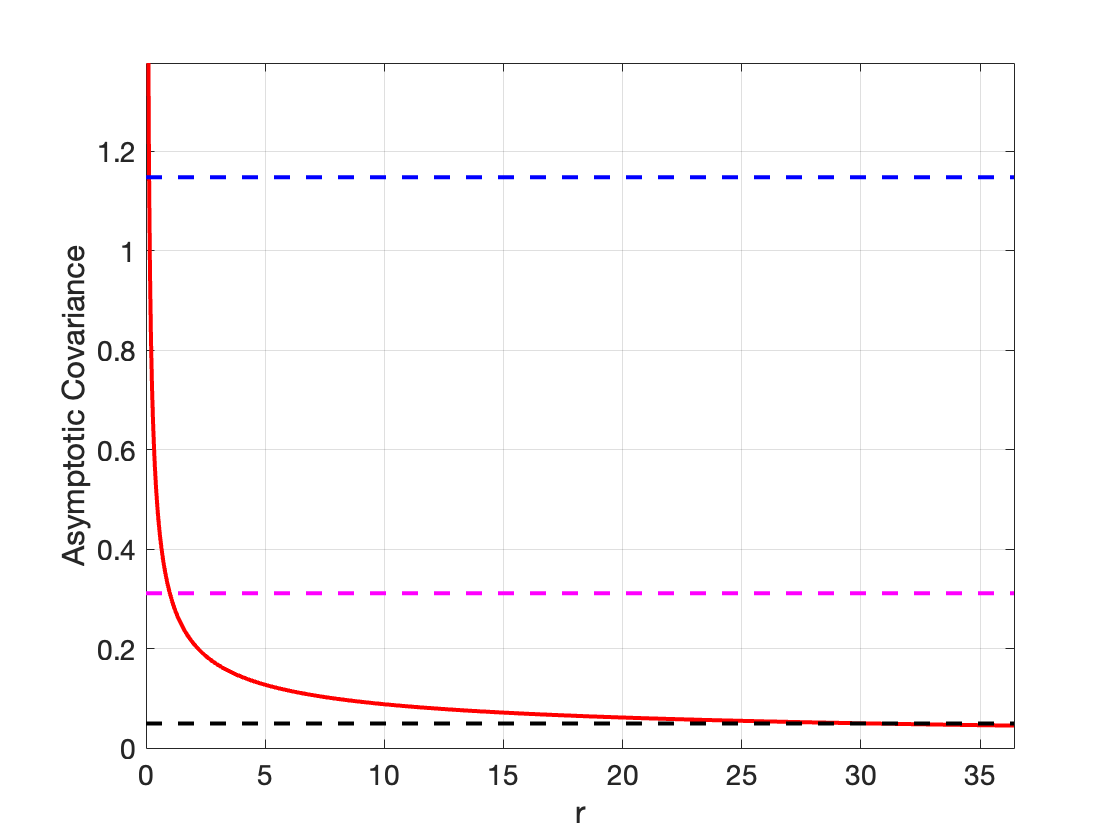}
		
	\end{minipage}
	\hfill 
	\begin{minipage}{0.45\textwidth}
		\centering
		\includegraphics[width=1\linewidth]{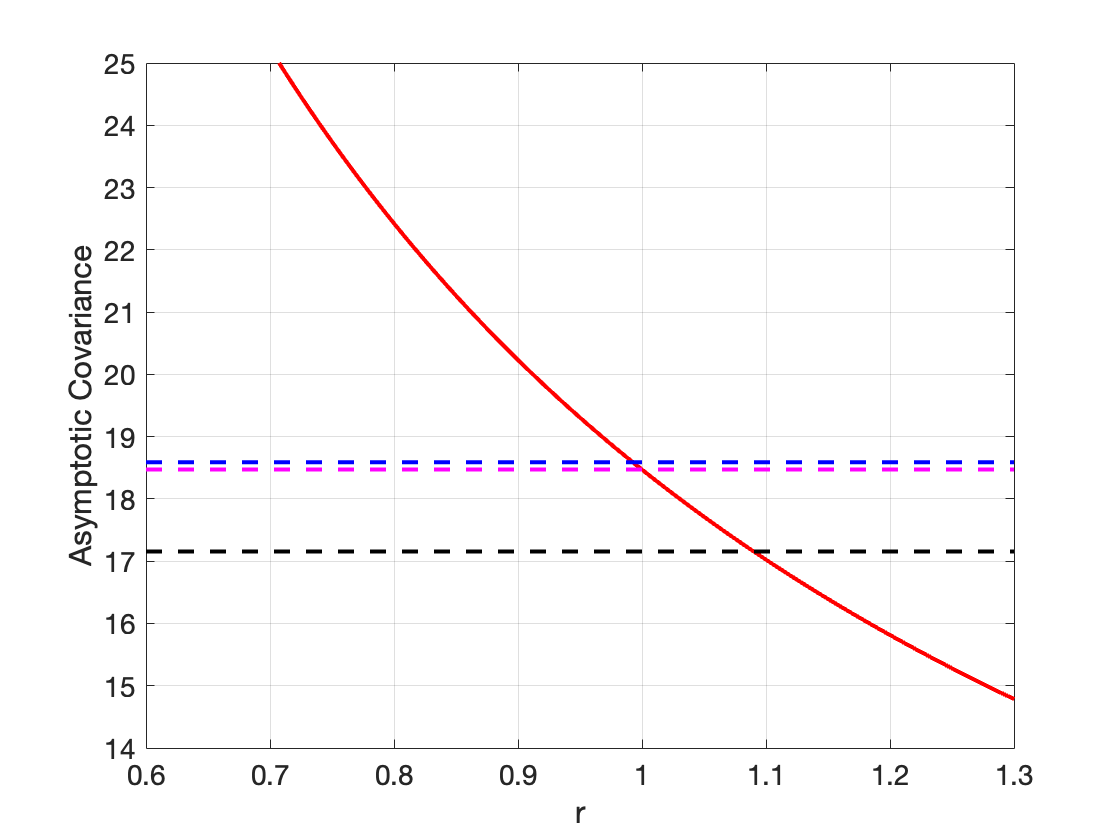}
	\end{minipage}
	\caption{Demonstration of more realistic/representative covariances that can be obtained with the non-local replicator mutator, i.e. with $s \neq 0$.  Left plot indicates system 1, right plot indicates system 2.  The red solid line indicates $C_\infty$ vs $r$ using \eqref{eq:asympCscalar}. The blue horizontal dashed line indicates the analytic $E_\infty^{\text{opt}}$, i.e. the minimum asymptotic MSE obtained by using $r^{\text{opt}}, s^{\text{opt}}$ in \eqref{eq:asympmsescalar}, the pink horizontal dashed line indicates $C_\infty$ for $r=1$ (i.e. $\hat{C}_\infty$ as defined in \eqref{eq:covperfect} , the covariance in the perfect model case).  Finally, the black dashed line indicates $C_\infty$ for $(s=0, r^\text{opt}_0)$, i.e. the optimal $r$ for multiplicative covariance inflation.} 
	\label{fig:covtest1}
\end{figure}

\begin{figure}[h!]
\begin{minipage}{0.45\textwidth} 
\centering
\includegraphics[width=1\linewidth]{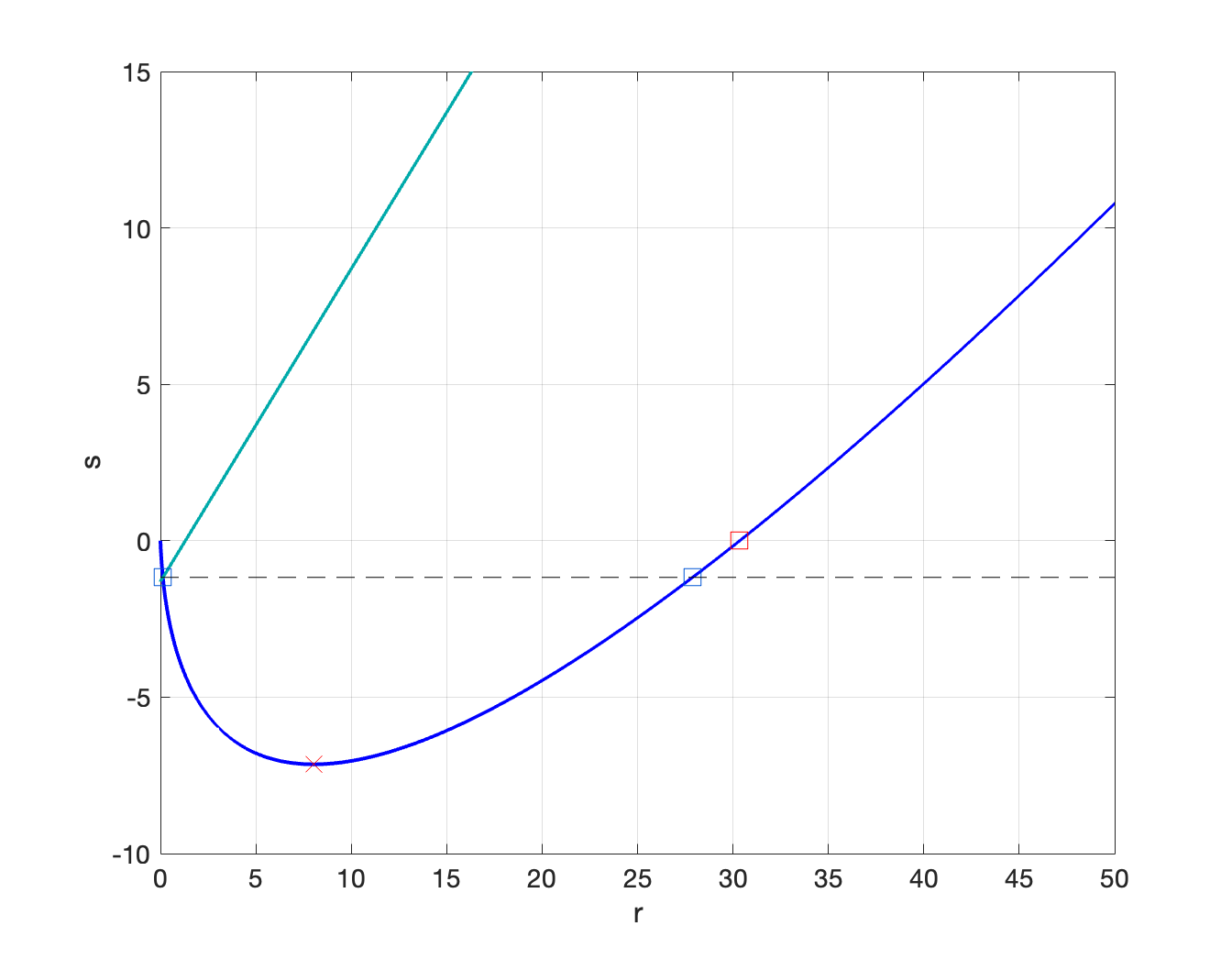}

\end{minipage}
\hfill 
\begin{minipage}{0.45\textwidth}
\centering
\includegraphics[width=1\linewidth]{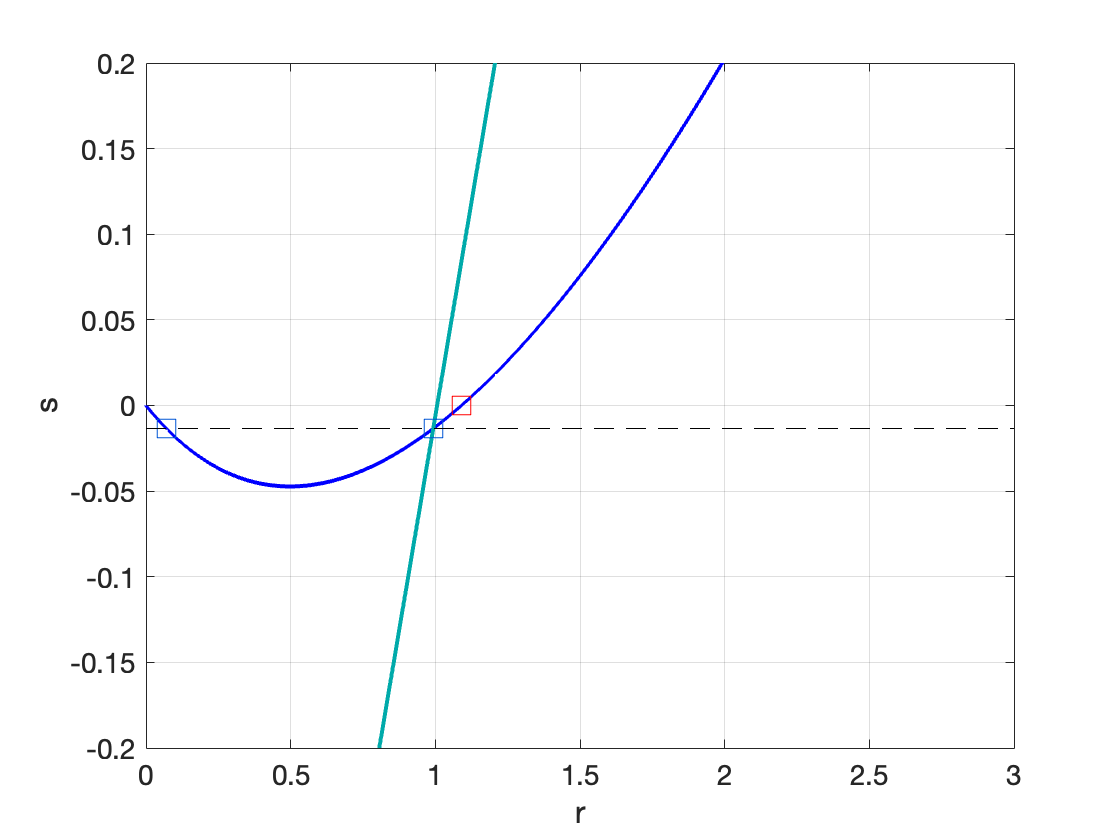}
\end{minipage}
\caption{Plot of optimal $(r,s)$ values for system 1 (left plot) and system 2 (right plot). 
The blue line indicates the $(r,s)$ pairs minimising mse only, obtained from \eqref{eq:optsr} and the cyan line indicates the $(r,s)$ pairs such that $C_\infty = E_\infty$, obtained from \eqref{eq:rminuscond}.  The point of intersection of the two lines indicates the $(r,s)$ pair that achieves both.  The red square indicates the optimal $r$ value corresponding to the multiplicative covariance inflation case $(s=0)$.  The blue squares indicate the possible $r$ values corresponding to the optimal $s$ in terms of $E_\infty$.  In system 1, multiplicative covariance inflation leads to overconfident estimates ($C_\infty$ too small), whereas in system 2, it is not so far off from the optimal choice of $s = -0.0135, r = 0.99$.}  
\label{fig:optrsplot}
\end{figure}

\paragraph{Uncertainty quantification with $C_\infty$.} Recall that we may further constrain the optimal $(r,s)$ values by enforcing the requirement that $C_\infty = \tilde{P}_\infty$.  This allows us to construct a filter/replicator-mutator dynamic for which the covariance of the estimate coincides with the actual mean squared error, as one obtains from the regular Kalman-Bucy filter in the perfect model setting. 
Figure \ref{fig:covtest1} shows the variation in asymptotic covariance $C_\infty$ for different $r$ values, obtained from \eqref{eq:asympCscalar}. Notice here that in the multiplicative covariance inflation case $s=0$, the asymptotic covariance for the corresponding optimal $r_0^{\text{opt}}$ is ($\tilde{C}_\infty = 0.05$), which is considerably smaller than the covariance in the perfect model setting  ($\tilde{C}_\infty = 0.31$), indicating overconfidence in the estimation (as also obtained in Lemma \ref{lem:covimpact}).  The non-local replicator-mutator on the other hand allows to obtain estimates that simultaneously minimise MSE and provide a realistic representation of uncertainty.  More specifically, the choice $r= 0.13$ (which coincides with $s^\text{opt} = -1.18$) yields a $C_\infty$ which coincides with $E_\infty$, so that the covariance produced by the estimation algorithm provides us with useful uncertainty quantification.  In particular, it represents an increase in uncertainty over the perfect knowledge case (pink line), which should be reflected given the unknown bias in the system.\\
\\
Finally, Figure \ref{fig:optrsplot} examines the optimal $r,s$ values minimising $E_\infty$ \eqref{eq:optsr} (blue line) and ensuring $C_\infty = E_\infty$,  \eqref{eq:rminuscond} in Lemma \ref{lem:covimpact} ( cyan line).  The intersection of these two lines indicates the single optimal $r,s$ pair that simultaneously achieves both.  Somewhat counter-intuitively, this pair does not always correspond to the smaller of the two possible $r$ values for the given $s$ value (see Figure \ref{fig:optrsplot} where in system 1 the smaller $r$ yields $C_\infty = E_\infty$ whereas the opposite is true for system 2).  In fact, for system 2, the optimal $s^{\text{opt}}$ has two corresponding optimal $r$ values that minimise $E_\infty$ only, $r^{\text{opt}} = 27.9$ or $0.13$.  The former $r$ value corresponds closely to the optimal $r$ for $s=0$ (pure multiplicative inflation) for which one obtains over-confident $C_\infty$.  However, the opposite is true for system 2, where the optimal $r,s$ where $C_\infty = E_\infty$ is relatively close to that of the pure multiplicative covariance inflation case.   We leave a study of specific system characteristics that would benefit most from the non-local replicator-mutator approach to future work.

\afterpage{\blankpage}

\newpage 

\section{Proofs}

\subsection{Proof of Lemma \ref{lem:fisherrao}}
\label{sec:prooffisherrao}

\begin{proof}
We now restrict $\mathcal{P}$ to be the set of smooth probability density functions 
\begin{align*}
\mathcal{P}(\mathbb{R}^d) := \Bigl\{ p \in C^\infty(\mathbb{R}^d) : \int p(x)dx = 1, \enskip p \geq 0 \Bigr\}
\end{align*}
whose tangent space is given by\footnote{Actually, if $p$ has a non-trivial support (e.g., there is an interval on which $p$ vanishes), then the tangent space needs to be replaced by a tangent cone $T_p \mathcal{P}(\mathbb{R}^d) = \Bigl\{ \sigma \in C_c^\infty(\supp(p)) : \int \sigma(x)dx = 0\Bigr\}$, see \cite{Mertikopoulos2018}, but we forgo the technical details here.}
\begin{align*}
T_p \mathcal{P}(\mathbb{R}^d) = \Bigl\{ \sigma \in C^\infty(\mathbb{R}^d) : \int \sigma(x)dx = 0\Bigr\}
\end{align*}
We start by computing the Fr\'echet derivative of \eqref{eq:expecfit} at $p$: Set $\eps > 0$ and $q\in \mathcal M$.
\begin{align*}
\frac{\mathcal F(q+\eps \tilde q) - \mathcal F(q)}{\eps} &= -\iint f(x,z)\frac{(q(x)+\eps \tilde q(x))(q(z)+\eps \tilde q(z)) - q(x)q(z)}{2\eps}\d z \d x\\
&= -\frac{1}{2}\iint f(x,z) \left[q(x)\tilde q(z) + q(z)\tilde q(x) + \eps \tilde q(x)\tilde q(z) \right]\d z \d x.
\end{align*}
This means that, because $f$ is symmetric in its components, 
\begin{align*}
D_q\mathcal F[\tilde q] &= -\iint f(x,z)q(z)\tilde q(x)\d z \d x = -\int \pi_q(x) \tilde q(x)\d x= -\langle \pi_q, \tilde q\rangle
\end{align*}
which shows that 
\begin{align*}
\mathcal{F}'(p) = -\pi_p 
\end{align*}
(this is to be understood as a linear operator). 
For the dissipation mechanism, consider the Fisher-Rao metric defined as 
\begin{align*}
g_p^{FR}(\sigma_1, \sigma_2) &=\int \frac{\sigma_1}{p(x)}\frac{\sigma_2}{p(x)} \d p(x) = \int \sigma_1 \frac{\sigma_2}{p(x)} \d x \\
& \equiv \int \sigma_1 \frac{\sigma_2}{p(x)} \d x - \int \sigma_1 dx \cdot \int \sigma_2 \d x
\end{align*}
since $\sigma \in T_p \mathcal{P}$ satisfies $\int \sigma dx = 0$.  Its corresponding isomorphism/dual action $\mathcal{G}^{FR}(p): T_p \mathcal{P} \rightarrow T_p^\ast \mathcal{P}$ is given (by inspection) as 
\begin{align}
\label{eq:gfr}
\mathcal{G}^{FR}(p)\sigma = \frac{\sigma}{p} - \int \sigma dx 
\end{align}
The inverse mapping is given by 
\begin{align}
\label{eq:gfrinv}
\mathcal{G}^{FR}(p)^{-1}\Phi = \left(\Phi - \int \Phi \, p(x)dx   \right)p.
\end{align}
Straightforwardly, we have
\begin{align*}
-\mathcal{G}^{FR}(p)^{-1} \mathcal{F}'(p) = \left(\pi_p(x) - \int \pi_p(x) p(x)dx   \right)p(x)
\end{align*}
which yields the right-hand side of the replicator equation.  

\end{proof}

\subsection{Proof of Lemma \ref{lem:ckreform}}
\label{sec:proofckreform}

\begin{proof}

The result follows from a simple re-arrangement of \eqref{eq:repmutnonlin}.  Throughout, we use the shorthand notation $y_t := \frac{dZ_t^d}{dt}$ and recall that $\|v(x) \|_{\Xi}^2 := v(x)^{\top} \Xi^{-1} v(x)$ for vector valued functions $v(x): \mathbb{R}^d \rightarrow \mathbb{R}^d$.  Firstly, 
\begin{align*}
\mathbb{E}_{z \sim \rho_t}[f_t(x,z)] -  \mathbb{E}_{\rho_t}[f_t]  &= \int \left(-\frac{r}{2}\|h(x) - y_t\|_{\Xi}^2 +  s \langle h(x) - y_t, h(z) - y_t \rangle_{\Xi} \right) \rho_t(z) dz \\
& - \int \int \left(  -\frac{r}{2}\|h(x) - y_t\|_{\Xi}^2 +  s \langle h(x) - y_t, h(z) - y_t \rangle_{\Xi}  \right)  \rho_t(x) \rho_t(z) dx dz  \\
& = -\frac{r}{2} \left( \|h(x) - y_t\|_{\Xi}^2 -  \mathbb{E}_{\rho_t}[\|h -y_t\|_{\Xi}^2] \right)  + s \langle h(x) - y_t, \mathbb{E}_{\rho_t}[h] - y_t \rangle_{\Xi} - s \| \mathbb{E}_{\rho_t}[h] - y_t \|_\Xi^2 \\
& =: I_1(x) + sI_2(x) + sI_3(x)
\end{align*}
For the first term, we have
\begin{align*}
I_1(x)  &  = -\frac{r}{2} \left( h(x)^\top\Xi^{-1} h(x) -h(x)^\top \Xi^{-1}y_t  - y_t^\top \Xi^{-1} h(x)   - \mathbb{E}_{\rho_t}[h^\top\Xi^{-1} h] + \mathbb{E}_{\rho_t}[h^\top] \Xi^{-1}y_t  +y_t^\top \Xi^{-1} \mathbb{E}_{\rho_t}[h] \right) \\
& = -\frac{r}{2} \left( h(x)^\top\Xi^{-1} h(x)   - \mathbb{E}_{\rho_t}[h(x)^\top\Xi^{-1} h(x)]\right)  +  r\left(h(x) - \mathbb{E}_{\rho_t}[h]\right)^\top\Xi^{-1}y_t
\end{align*}
For the remaining terms, we have that  
\begin{align*}
I_2(x) + I_3(x)  &=   \langle h - \mathbb{E}_{\rho_t}[h], \mathbb{E}_{\rho_t}[h] - y_t \rangle_\Xi \\
& = -\frac{1}{2} \| h - \mathbb{E}_{\rho_t}[h] \|_\Xi^2 - \frac{1}{2} \mathbb{E}_{\rho_t}[h^\top] \Xi^{-1} \mathbb{E}_{\rho_t}[h] + \frac{1}{2}h^T \Xi^{-1} h - \frac{1}{2} \mathbb{E}_{\rho_t}[h^T \Xi^{-1} h] \\
& + \frac{1}{2} \mathbb{E}_{\rho_t}[h^T \Xi^{-1} h] - (h - \mathbb{E}_{\rho_t}[h])^T \Xi^{-1} y_t \\
&  = -\frac{1}{2} \| h - \mathbb{E}_{\rho_t}[h] \|_\Xi^2  + \frac{1}{2} \left( \mathbb{E}_{\rho_t}[h^T \Xi^{-1} h] -  \mathbb{E}_{\rho_t}[h^\top] \Xi^{-1} \mathbb{E}_{\rho_t}[h] \right)\\
& + \frac{1}{2} \left( h^T \Xi^{-1} h - \mathbb{E}_{\rho_t}[h^T \Xi^{-1} h]  \right)  - (h - \mathbb{E}_{\rho_t}[h])^T \Xi^{-1} y_t. 
\end{align*}
Combining yields
\begin{align*}
I_1(x) + sI_2(x) + sI_3(x) &= -\frac{1}{2} (r-s) \left( h(x)^\top\Xi^{-1} h(x)   - \mathbb{E}_{\rho_t}[h(x)^\top\Xi^{-1} h(x)]\right)  + (r - s) \left(h(x) -  \mathbb{E}_{\rho_t}[h] \right)^\top\Xi^{-1}y_t \\
& -\frac{s}{2} \| h - \mathbb{E}_{\rho_t}[h] \|_\Xi^2  + \frac{s}{2} \left( \mathbb{E}_{\rho_t}[h^T \Xi^{-1} h] -  \mathbb{E}_{\rho_t}[h^\top] \Xi^{-1}\mathbb{E}_{\rho_t}[h] \right).
\end{align*}
Substituting the above into \eqref{eq:repmutnonlin} yields \eqref{eq:repasKS}.  

\end{proof}

\subsection{Proof of Theorem \ref{theo:limitcrowkimura}}
\label{sec:prooftheolimitcrowkimura}

\begin{proof} 

Start with the reformulation of \eqref{eq:repmutnonlin} as derived in Lemma \ref{lem:ckreform}, (repeating here for convenience) 
\begin{align}
\label{eq:thm4ck}
\partial_t \mu_t^d(x) = \mathcal{L}^\ast \mu_t^d(x) +  \left(   - \frac{r}{2} h(x)^\top\Xi^{-1} h(x) +  (r-s) h(x)^\top \Xi^{-1}\frac{dZ_t^d}{dt}  \right) \mu_t^d(x) 
\end{align}
and the Stratonovich form of \eqref{eq:infzakai}, 
\begin{align}
\label{eq:thm4zakai}
d q_t = \mathcal{L}^\ast q_t(x)dt -\frac{r}{2}h(x)^\top \Xi^{-1} h(x) q_t(x)dt  + (r-s)q_t(x) h(x)^\top \Xi^{-1} \circ dZ_t. 
\end{align}
We proceed with the following steps.  The proof below is inspired by the proof of Theorem 3.1 in \cite{hu_approximation_2002} except for the following important extensions: 1) we no longer assume $h$ is uniformly bounded; 2) $h$ is no longer a scalar valued function but may be vector valued; 3) we make use of the forward stochastic Feynman-Kac style representation formula in \cite{kunita_stochastic_1982} rather than the backward formula.  Existence and uniqueness of density valued solutions to the zakai equation in $ L^2(\mathbb{R}^m)$ with $g,h$ unbounded has been studied by a number of authors \cite{baras_existence_1983,bhatt_uniqueness_1995}, building on extensive works in the unbounded case (see e.g. the excellent summary in \cite{bain_fundamentals_2009}).\\ 
\\
\textbf{Step 1.} Use (stochastic) Feynman-Kac type formulae to obtain a probabilistic representation of solutions to the Zakai equation and replicator-mutator equation, as given in Theorem \ref{theo:pdeprobrep}.  Specifically, we make use of the formulae developed in \cite{kunita_stochastic_1982, kunita_cauchy_1981} as was done in \cite{hu_approximation_2002}  although  we rely on the forward rather than backward representation formulae.  Recall that $\mathcal{L}^\ast$ denotes the adjoint operator of the generator of the diffusion process 
\begin{align*}
dX_t = g(X_t)dt + \sigma(X_t)dW_t 
\end{align*}
By expanding the adjoint operator, we may express it as 
\begin{align*}
\mathcal{L}^\ast \mu_t^d(x) &= -\sum_{i=1}^m \frac{\partial}{\partial x^i} \left(g^i(x) \mu_t^d(x)  \right) + \frac{1}{2} \sum_{i=1}^m \sum_{j=1}^m \frac{\partial^2}{\partial x^i \partial x^j} \left( (\sigma \sigma^\top)^{ij} \mu_t^d(x) \right) \\
& = -\sum_{i=1}^m  \left( \mu_t^d(x) \frac{\partial g^i(x)}{\partial x^i} + g^i(x) \frac{\partial \mu_t^d(x)}{\partial x^i} \right)  + \frac{1}{2} \sum_{i=1}^m \sum_{j=1}^m \left( \mu_t^d \frac{\partial^2 (\sigma \sigma^\top)^{ij}}{\partial x^i \partial x^j} + (\sigma \sigma^\top)^{ij} \frac{\partial^2 \mu_t^d}{\partial x^i \partial x^j} \right) \\
&+ \sum_{i=1}^m \frac{\partial \mu_t^d(x)}{\partial x^i} \sum_{j=1}^m \frac{\partial (\sigma \sigma^\top)^{ij}}{\partial x^j} \\
& = \mathcal{G} \mu_t^d(x) + \left( -\sum_{i=1}^m \frac{\partial g^i(x)}{\partial x^i}  
+ \frac{1}{2} \sum_{i=1}^m \sum_{j=1}^m  \frac{\partial^2 (\sigma \sigma^\top)^{ij}}{\partial x^i \partial x^j} \right) \cdot \mu_t^d(x) 
\end{align*}
where in the second line we have used that $\sigma \sigma^\top$ is symmetric and $\mathcal{G}$ denotes the infinitesimal generator of the diffusion process
\begin{align}
\label{eq:diffproc}
dY_t = \left(b(Y_t) -g(Y_t) \right) dt + \sigma(Y_t)dW_t^y 
\end{align}
where $W_t^y$ is a Wiener process independent of $W_t$,  $b^i(x) = \nabla \cdot (\sigma(x) \sigma(x)^\top)^i$ and the superscript $i$ denotes the $i$th row of the matrix $\sigma \sigma^\top$.  That is, 
\begin{align*}
\mathcal{G} \mu_t^d(x) := -\sum_{i=1}^m  g^i(x) \frac{\partial \mu_t^d(x)}{\partial x^i}   + \frac{1}{2} \sum_{i=1}^m \sum_{j=1}^m  (\sigma \sigma^\top)^{ij} \frac{\partial^2 \mu_t^d}{\partial x^i \partial x^j} + \sum_{i=1}^m \frac{\partial \mu_t^d(x)}{\partial x^i} \sum_{j=1}^m \frac{\partial (\sigma \sigma^\top)^{ij}}{\partial x^j}
\end{align*}
This decomposition of $\mathcal{L}^\ast$ will be used in both \eqref{eq:thm4ck} and \eqref{eq:thm4zakai}; starting with \eqref{eq:thm4zakai} and using \eqref{eq:obs} yields,
\begin{align*}
d q_t(x) &= \mathcal{L}^\ast q_t(x) -\frac{r}{2} h(x)^\top \Xi^{-1} h(x) q_t(x) +  (r-s)q_t(x) h(x)^\top \Xi^{-1} \circ dZ_t \\
& = \mathcal{G}q_t(x) + \left( -\sum_{i=1}^m \frac{\partial g^i(x)}{\partial x^i}  + \frac{1}{2} \sum_{i=1}^m \sum_{j=1}^m  \frac{\partial^2 (\sigma \sigma^\top)^{ij}}{\partial x^i \partial x^j} -\frac{r}{2} h(x)^\top \Xi^{-1} h(x)  +  (r-s) h(x)^\top \Xi^{-1} h(x_t^\ast)dt \right) \cdot q_t(x) \\
& + (r-s)q_t(x) h(x)^\top \Xi^{-1/2} \circ dB_t
\end{align*}
This equation now takes the form of \eqref{eq:parapde2} in Theorem \ref{theo:pdeprobrep} with
\begin{align*}
l^{ik} &= 0 \\
c^k &= (r-s)(h(x)^\top \Xi^{-1/2})^k, \quad k = 1, 2, \dots, n \\
c^0 &=  -\sum_{i=1}^m \frac{\partial g^i(x)}{\partial x^i}  + \frac{1}{2} \sum_{i=1}^m \sum_{j=1}^m  \frac{\partial^2 (\sigma \sigma^\top)^{ij}}{\partial x^i \partial x^j}  -\frac{r}{2} h(x)^\top \Xi^{-1} h(x) + (r-s)h(x)^\top \Xi^{-1} h(x_t^\ast) \\
a &= \sigma \\
b^i &= -g^i(x) + \sum_{j=1}^m \sum_{k=1}^m \sigma^{ik} \frac{\partial \sigma^{kj}}{\partial x^j} + \frac{1}{2} \sigma^{kj} \frac{\partial \sigma^{ik}}{\partial x^j}
\end{align*}
recalling that $x_t^\ast$ is treated as a fixed realisation.  Then by Theorem \ref{theo:pdeprobrep}, there exists another probability space equipped with the measure $\mathbb{Q}$ (from here on we use the notation $\mathbb{E}^\mathbb{Q}$ denote the expectation with respect to this measure) such that the solution can be represented as
\begin{align*}
q_t(x) = \mathbb{E}^{\mathbb{Q}} \left[ f(\xi_t(x)) \exp \left(  (r-s)\int_0^t h(\xi_u(x))^\top \Xi^{-1} \circ dZ_u  - \frac{r}{2}\int_0^t h(\xi_u(x))^\top \Xi^{-1} h(\xi_u(x)) du + \int_0^t \tilde{c}(\xi_u(x))du \right) \right]
\end{align*}
where $\xi_s(x)$ is a vector-valued function of $x$ denoting the solution of an SDE in the form of \eqref{eq:sdefeynmankac} in Theorem \ref{theo:pdeprobrep} with $b, a, l$ as defined above, i.e. 
\begin{align}
\label{eq:sderepform}
d \xi_t(x) =  \sum_{i=1}^m b^{i}(\xi_t(x))  dt + \sum_{j=1}^m   \sum_{i=1}^m \sigma^{ij}(\xi_t(x)) \circ dW_t^j 
\end{align}
and $f(x) = \lim_{t \rightarrow 0} q_t(x)$ denotes the initial density and  
\begin{align}
\label{eq:feynctilde}
\tilde{c}(x) := -\sum_{i=1}^m \frac{\partial g^i(x)}{\partial x^i}  + \frac{1}{2} \sum_{i=1}^m \sum_{j=1}^m  \frac{\partial^2 (\sigma(x) \sigma^\top(x))^{ij}}{\partial x^i \partial x^j}
\end{align}
We can similarly apply Theorem \ref{theo:pdeprobrep} to \eqref{eq:thm4ck} with 
\begin{align*}
c^k &=  0 \\
c^0 &=  -\sum_{i=1}^m \frac{\partial g^i(x)}{\partial x^i}  + \frac{1}{2} \sum_{i=1}^m \sum_{j=1}^m  \frac{\partial^2 (\sigma \sigma^\top)^{ij}}{\partial x^i \partial x^j}  -\frac{r}{2} h(x)^\top \Xi^{-1} h(x) + (r-s)h(x)^\top \Xi^{-1} \frac{dZ_t^d}{dt} 
\end{align*}
and $l, a, b$ as defined previously, to obtain the representation 
\begin{align*}
\mu_t^d(x) = \mathbb{E}^{\mathbb{Q}} \left[ f(\xi_t(x)) \exp \left( (r-s)\int_0^t h(\xi_u(x))^\top \Xi^{-1} \dot{Z}_u^d du   - \frac{r}{2}\int_0^t h(\xi_u(x))^\top \Xi^{-1} h(\xi_u(x)) du +  \int_0^t \tilde{c}(\xi_u(x))du \right) \right]
\end{align*}
where $\tilde{c}(x)$ is given in \eqref{eq:feynctilde} and $\xi_s$ is as defined previously since $l^{ik} = 0$ in both cases and from now onwards we use the shorthand notation $\dot{Z}_t^d \equiv \frac{dZ_t^d}{dt}$. Recall also that both \eqref{eq:thm4ck} and \eqref{eq:thm4zakai} are assumed to be initialised by the same density $f(x)$.\\
\\
\textbf{Step 2.} We are now ready to prove pointwise convergence using the above representation formulae. Firstly, we have
\begin{align*}
\mu_t^d(x) - q_t(x) &= \mathbb{E}^\mathbb{Q} \left[ f(\xi_t(x)) \exp \left(   - \frac{r}{2}\int_0^t h(\xi_u(x))^\top \Xi^{-1} h(\xi_u(x)) du + (r-s)\int_0^t h(\xi_u(x))^\top \Xi^{-1} \dot{Z}_u^d du + \int_0^t \tilde{c}(\xi_u(x))du \right)  \right] \\
& - \mathbb{E}^{\mathbb{Q}} \left[ f(\xi_t(x)) \exp \left( - \frac{r}{2}\int_0^t h(\xi_u(x))^\top \Xi^{-1} h(\xi_u(x)) du +  (r-s)\int_0^t h(\xi_u(x))^\top\Xi^{-1} \circ dZ_u   + \int_0^t \tilde{c}(\xi_u(x))du \right) \right] \\
& =: \mathbb{E}^\mathbb{Q} \left[ f(\xi_t(x)) I_2(x) I_1(x; \mathcal{Z}) \right]
\end{align*}
where 
\begin{align*}
I_2(x) &=  \exp \left( \int_0^t \tilde{c}(\xi_s(x)) - \frac{r}{2} h(\xi_s(x))^\top \Xi^{-1} h(\xi_s(x)) ds \right) \\
I_1(x; \mathcal{Z}) &= \exp \left( (r-s) \int_0^t h(\xi_s(x))^\top \Xi^{-1} \dot{Z}_s^d ds  \right)  - \exp \left( (r-s) \int_0^t h(\xi_s(x))^\top\Xi^{-1} \circ dZ_s   \right)  \\
& =: \exp(I_3(x; \mathcal{Z})) - \exp(I_4(x; \mathcal{Z}))
\end{align*}
and the $\mathcal{Z}$ notation is used to denote the dependence of $I_1$ on the observation path.  In the below, let $\mathbb{E}$ refer to the expectation on the original space, i.e. wrt to the observation noise $W$.  Then for any fixed $x \in \mathbb{R}^m$, in other words, a realisation of the initialisation which has probability density $f(x)$, 
\begin{align*}
\mathbb{E} \left[| \mu_t^d(x) - q_t(x) \|^p \right]  &\leq \mathbb{E} \left[ \mathbb{E}^\mathbb{Q} \left[ | f(\xi_t(x)) I_2(x) I_1(x; \mathcal{Z}) |^p  \right] \right] \\
& =  \mathbb{E}^\mathbb{Q} \left[  |f(\xi_t(x)) I_2(x)|^p \mathbb{E} \left[|I_1(x; \mathcal{Z})|^p   \right]  \right]  \quad \text{(Fubini)} \\
& \leq C  \left( \mathbb{E}^\mathbb{Q} \left[ |f(\xi_t(x)) I_2(x)|^{p r_2}  \right]  \right)^{1/r_2} \cdot \left( \mathbb{E}^\mathbb{Q} \left[ (\mathbb{E} \left[|I_1(x; \mathcal{Z})|^p \right]   )^{r_1} \right] \right)^{1/r_1}  \quad \text{(Hoelder inequality)} \\
& =: CI_5^{1/r_2} \cdot I_6^{1/r_1}  
\end{align*}
for a constant $C>0$ independent of $t$ and with $1/r_1 + 1/r_2 = 1$ and $r_1, r_2 > 1$.  Starting with $I_2$, since $g(x)$ is assumed to be $C^2$ and globally Lipschitz continuous and $\Sigma$ is a constant, it holds that $\tilde{c}(x)$ as defined in \eqref{eq:feynctilde} is uniformly bounded.  Additionally, it holds that $h(x)^\top \Xi^{-1} h(x) > 0, \enskip \forall \enskip x \in \mathbb{R}^m$ since $\Xi$ is a positive definite matrix, so that there exists some $C_2> 0$ independent of $t$.  Combining, we have that 
\begin{align*}
|I_2(x)|^{pr_2} = \exp \left( pr_2\int_0^t \tilde{c}(\xi_u(x)) du \right) \cdot \exp \left( - p r_2\frac{r}{2} \int_0^t h(\xi_u(x))^\top \Xi^{-1} h(\xi_u(x)) du \right) \leq \exp(C_1 t) \cdot C_2  
\end{align*}
from which we obtain  
\begin{align*}
I_5 &\leq   \mathbb{E}^\mathbb{Q} \left[ |f(\xi_t(x))|^{p r_2}  |I_2(x)|^{p r_2}  \right]   \\
&\leq C_3(t) 
\end{align*}
using the fact that $f$ is uniformly bounded. Now turning to $I_6$, using the identity 
\begin{align*}
|\exp(x) - \exp(y)|^p &\leq (\exp(x) + \exp(y))^p|x-y|^p, \quad x,y \in \mathbb{R} \\
& \leq C(\exp(px) + \exp(py))|x-y|^p
\end{align*}
along with Hoelder inequality yields  
\begin{align*}
I_6 &\leq C  \mathbb{E}^\mathbb{Q} \left[ (\mathbb{E} \left[ (\exp(p I_3(x; \mathcal{Z})) + \exp(p I_4(x; \mathcal{Z}))| I_3(x; \mathcal{Z}) - I_4(x; \mathcal{Z})|^p \right]   )^{r_1} \right]   \\
& \leq C  \mathbb{E}^\mathbb{Q} \left[ (\mathbb{E} \left[ \exp(p I_3(x; \mathcal{Z}))| I_3(x; \mathcal{Z}) - I_4(x; \mathcal{Z})|^p \right]   )^{r_1} + (\mathbb{E} \left[  \exp(p I_4(x; \mathcal{Z}) | I_3(x; \mathcal{Z}) - I_4(x; \mathcal{Z})|^p \right])^{r_1} \right]  \\
& \leq C   \mathbb{E}^\mathbb{Q} \left[ \left( \mathbb{E} \left[ \exp(r_2 p I_3(x; \mathcal{Z})) \right]  \right)^{r_1/r_2} \cdot \mathbb{E}[|I_3(x;\mathcal{Z}) - I_4(x; \mathcal{Z})|^{pr_1}] \right]   \\
& + C  \mathbb{E}^\mathbb{Q} \left[ \left( \mathbb{E} \left[ \exp(r_2 p I_4(x; \mathcal{Z})) \right]  \right)^{r_1/r_2} \cdot \mathbb{E}[|I_3(x;\mathcal{Z}) - I_4(x; \mathcal{Z})|^{pr_1}] \right]  \\
& \leq C \left[\left( \mathbb{E}^\mathbb{Q}[I_9(x)^{r_1}] \right)^{1/r_2} + \left( \mathbb{E}^\mathbb{Q}[I_8(x)^{r_1}] \right)^{1/r_2} \right] \left(  \mathbb{E}^\mathbb{Q}[I_7(x)^{r_1}]\right)^{1/r_1}
\end{align*}
with 
\begin{align*}
I_7(x) := \mathbb{E}[|I_3(x;\mathcal{Z}) - I_4(x; \mathcal{Z})|^{pr_1}], \quad I_8(x):=  \mathbb{E} \left[ \exp(r_2 p I_4(x; \mathcal{Z})) \right], \quad I_9(x):=\mathbb{E} \left[ \exp(r_2 p I_3(x; \mathcal{Z})) \right] 
\end{align*}

\noindent \textbf{Step 3.} The remainder of the proof will focus on showing that the terms involving $I_8$ and $I_9$ can be bounded by constants (depending on $t$ only).  The term involving $I_7$ will the be shown to go to zero as $d \rightarrow \infty$, yielding the desired convergence result.  Starting with $I_8(x)$ and using the shorthand notation $p_r := p(r-s)$,
\begin{align}
\nonumber 
&\mathbb{E} \left[ \exp(r_2 p_r I_4(x; \mathcal{Z})) \right] = \mathbb{E} \left[ \exp \left(r_2 p_r \int_0^t  h(\xi_u(x))^\top \Xi^{-1} \circ dZ_u   \right) \right] \\
\nonumber 
& = \mathbb{E} \left[ \exp \left( r_2 p_r \int_0^t  h(\xi_u(x))^\top\Xi^{-1}h(x_u^\ast)du +   r_2 p_r \int_0^t h(\xi_u(x))^\top\Xi^{-1/2} dB_u  \right)  \right] \\   
\nonumber 
& \leq  \left( \mathbb{E} \left[ \exp \left( r_1 r_2 p_r \int_0^t  h(\xi_u(x))^\top\Xi^{-1}h(x_u^\ast)du \right) \right] \right)^{1/r_1} \cdot \left( \mathbb{E} \left[ \exp \left(r_2^2 p_r \int_0^t h(\xi_u(x))^\top\Xi^{-1/2} dB_u  \right)  \right] \right)^{1/r_2}  \\
\nonumber 
& \leq  \left( \mathbb{E} \left[ \exp \left( r_1 r_2 p_r \int_0^t  h(\xi_u(x))^\top\Xi^{-1}h(x_u^\ast)du \right) \right] \right)^{1/r_1} \cdot \left( \prod_{k=1}^n \mathbb{E} \left[ \exp \left(n r_2^2 p_r \int_0^t (h(\xi_u(x))^\top\Xi^{-1/2})^k dB_u^k  \right)  \right] \right)^{1/{(r_2n)}}  \\
\nonumber 
& \leq \left( \mathbb{E} \left[ \exp \left( r_1 r_2 p_r \int_0^t  h(\xi_u(x))^\top\Xi^{-1}h(x_u^\ast)du \right) \right] \right)^{1/r_1} \cdot \left( \exp \left( \frac{n r_2^2 p_r}{2} \sum_{k=1}^n  \int_0^t  ((h(\xi_u(x))^\top\Xi^{-1/2})^k)^2  du  \right) \right)^{1/{(r_2n)}} \\
\label{eq:i8bound}
& = \left( \mathbb{E} \left[ \exp \left( r_1 r_2 p_r \int_0^t  h(\xi_u(x))^\top\Xi^{-1}h(x_u^\ast)du \right) \right] \right)^{1/r_1} \cdot \left( \exp \left( \frac{n r_2^2 p_r}{2}   \int_0^t  h(\xi_u(x))^\top\Xi^{-1}h(\xi_u(x))  du  \right) \right)^{1/{(r_2n)}} 
\end{align}
where in the last line, we have used Lemma \eqref{lem:mgf} and the fact that $\xi_s$ is defined on a different probability space to the signal process. 
Then we have
\begin{align*}
\left( \mathbb{E}^\mathbb{Q}[I_8(x)^{r_1}] \right)^{1/r_2} &\leq \left( \mathbb{E}^\mathbb{Q} \left[ \mathbb{E} \left[ \exp \left( r_1 r_2 p_r \int_0^t  h(\xi_u(x))^\top\Xi^{-1}h(x_u^\ast)du \right) \right]  \cdot \exp \left( \frac{ r_1 r_2 p_r}{2}   \int_0^t  h(\xi_u(x))^\top\Xi^{-1}h(\xi_u(x))  du  \right)   \right] \right)^{1/r_2} \\
& = \left( \mathbb{E}^\mathbb{Q} \left[  \exp \left( r_1 r_2 p_r \int_0^t  h(\xi_u(x))^\top\Xi^{-1}h(x_u^\ast)du \right)  \cdot \exp \left( \frac{ r_1 r_2 p}{2}   \int_0^t  h(\xi_u(x))^\top\Xi^{-1}h(\xi_u(x))  du  \right)   \right] \right)^{1/r_2} \\
&  =  \left( \mathbb{E}^\mathbb{Q} \left[  \exp \left( r_1 r_2 p_r \int_0^t  h(\xi_u(x))^\top\Xi^{-1} \left( h(x_u^\ast) + \frac{1}{2} h(\xi_u(x)) \right)du  \right)   \right] \right)^{1/r_2} \\
&=: I_{10}(x)^{1/r_2}
\end{align*}
since the expectation $\mathbb{E}$ is with respect to the observation noise and  $x_s^\ast$ is taken as a fixed realisation. 

To proceed further, we will make use of Lemma \ref{lem:delmeyer} which allows us to bound exponential moments of a non-decreasing process by its raw moments.  Generally, it is not possible to bound exponential moments in terms of polynomial moments, as the exponential term grows faster.  The crucial point of this lemma is in a careful specification of the factor in the exponential (c.f. $L$ in Lemma \ref{lem:delmeyer}) which acts to ``dampen'' the growth relative to the growth of the raw moments.  We have using the shorthand notation $\tilde{L} := \frac{r_1 r_2 p_r}{2}$,
\begin{align*}
I_{10}(x) \leq \mathbb{E}^\mathbb{Q} \left[\exp \left( \tilde{L}  \int_0^t \tilde{h}_u(\xi_u(x))^\top \Xi^{-1} \tilde{h}_u(\xi_u(x)) du \right) \right] 
\end{align*}
where $\tilde{h}_u(y):= h(y) + h(x_u^\ast)$.
Furthermore, 
\begin{align*}
\tilde{h}_u(x)^\top \Xi^{-1} \tilde{h}_u(x) = |\Xi^{-1/2} \tilde{h}_u(x) |^2 \leq \lambda_{\Xi}^2 |\tilde{h}_u(x)|^2, \quad \forall \enskip x \in \mathbb{R}^m,  
\end{align*}
where $\lambda_{\Xi}$ is the smallest eigenvalue of $\Xi^{1/2}$.  
Combining, we have letting $L = \tilde{L} \lambda_{\Xi}^2$,
\begin{align*}
I_{10}(x) \leq  \mathbb{E}^\mathbb{Q} \left[\exp \left( L  \int_0^t |\tilde{h}_u(\xi_u(x))|^2  du \right) \right] 
\end{align*}
Let $Y_t(x):= \int_0^t |\tilde{h}_u(\xi_u(x))|^2 du$. Clearly this is a non-decreasing (and adapted) process, so that we may apply lemma \ref{lem:delmeyer}. In particular, 
\begin{align*}
\mathbb{E}^\mathbb{Q}[Y_t(x) - Y_\tau(x)] = \int_\tau^t  \mathbb{E}^\mathbb{Q} \left[|\tilde{h}_u(\xi_u(x))|^2 \right] du \leq tC(1+ \mathbb{E}^\mathbb{Q}[|x|^2]) =: K(t), \quad \forall \enskip \tau \in [0,t) 
\end{align*}
where $C$ is a constant depending on the growth properties of $h$ and for the last inequality, we have used that since $g,h$ satisfy lipschitz and linear growth assumptions, 
\begin{align*}
\mathbb{E}^\mathbb{Q}[|\tilde{h}_u(\xi_u(x))|^2] \leq C(1 + \mathbb{E}^\mathbb{Q}[|\xi_u(x)|^2]) \leq C(1+ \mathbb{E}^\mathbb{Q}[|x|^2])
\end{align*}
Therefore by Lemma \ref{lem:delmeyer}, whenever
\begin{align}
\label{eq:expgrowthcond}
\frac{r_1r_2 (r-s)p \lambda_\Xi^2}{2} < \frac{1}{K(t)}
\end{align}
where $K(t)$ is a constant depending on time, the second moment of the initial density $f(x)$ and the linear growth constant of $h$, we have
\begin{align*}
I_{10}(x) \leq  \mathbb{E}^\mathbb{Q} \left[\exp \left( L  \int_0^t |\tilde{h}_s(\xi_s(x))|^2  ds \right) \right] < \frac{1}{1 - L K(t)}.
\end{align*}
Notice that condition \eqref{eq:expgrowthcond} can be satisfied whenver $r,s$ are chosen such that $(r-s)$ is small enough.  As will be seen in Section \ref{sec:biasvar}, this is at least possible in the linear-Gaussian setting whilst maintaining optimality (in the mean squared error sense) even in the case of a misspecified model. 
Finally, we have
\begin{align*}
\left( \mathbb{E}^\mathbb{Q}[I_8(x)^{r_1}] \right)^{1/r_2} \leq (1 - L K(t))^{-\frac{1}{r_2}} 
\end{align*}

The term involving $I_9(x)$ can be analysed in much the same way as for $I_8(x)$, with the only difference being that the stochastic integral.  In particular, we have letting $j$ refer to the index such that $t \in (t_j, t_{j+1}]$, 
\begin{align*}
&\mathbb{E} \left[ \exp(r_2 p_r I_3(x; \mathcal{Z})) \right] = \mathbb{E} \left[ \exp \left(r_2 p_r \int_0^t  h(\xi_u(x))^\top \Xi^{-1} \dot{Z}_u^d du    \right) \right] \\
& \leq  \left( \mathbb{E} \left[ \exp \left( r_1 r_2 p_r \int_0^t  h(\xi_u(x))^\top\Xi^{-1}h(x_u^\ast)du \right) \right] \right)^{1/r_1} \cdot \left( \prod_{k=1}^n \mathbb{E} \left[ \exp \left(n r_2^2 p_r \int_0^t (h(\xi_u(x))^\top\Xi^{-1/2})^k \dot{B}_u^k du \right)  \right] \right)^{1/{(r_2n)}}. 
\end{align*}
where with a slight abuse of notation, we let $\dot{B}_t \equiv \dot{B}_t^d$ and $\dot{B}_t^k$ denotes the $k$th component of $\dot{B}_t$.  Furthermore, 
\begin{align*}
&\mathbb{E} \left[ \exp \left(n r_2^2 p_r \int_0^t (h(\xi_u(x))^\top\Xi^{-1/2})^k \dot{B}_\tau^k d\tau  \right)  \right]  \\
&= \mathbb{E} \left[ \exp \left(n r_2^2 p_r \sum_{i=1}^j \frac{1}{\delta_d} \int_{t_{i-1}}^{t_i}   \left( \int_{t_{i-1}}^{t_i}   (h(\xi_u(x))^\top\Xi^{-1/2})^k  du\right) dB_\tau^k  + nr_2^2 p_r \int_{t_j}^t (h(\xi_s(x))^\top\Xi^{-1/2})^k \dot{B}_\tau^k d\tau  \right)  \right] \\
&= \prod_{i=1}^j \mathbb{E} \left[ \exp \left(n r_2^2 p_r  \frac{1}{\delta_d}  \left( \int_{t_{i-1}}^{t_i}   (h(\xi_u(x))^\top\Xi^{-1/2})^k  du\right) \int_{t_{i-1}}^{t_i}  dB_\tau^k \right) \right] \mathbb{E} \left[ \exp \left( nr_2^2 p_r \left( \int_{t_{j}}^{t}   (h(\xi_u(x))^\top\Xi^{-1/2})^k  du\right) \int_{t_{j}}^{t_{j+1}}  dB_\tau^k   \right)  \right] \\
& \leq \prod_{i=1}^j \exp \left( \frac{n r_2^2 p_r}{2}  \int_{t_{i-1}}^{t_i}   (h(\xi_u(x))^\top\Xi^{-1/2})^k  du\right) \cdot \exp \left( \frac{n r_2^2 p_r}{2}   \int_{t_{j}}^{t}   (h(\xi_u(x))^\top\Xi^{-1/2})^k  du \right) \\
& = \exp \left( \frac{n r_2^2 p_r}{2}  \int_0^t (h(\xi_u(x))^\top\Xi^{-1/2})^k du \right) \\
& \leq  \exp \left( \frac{n r_2^2 p_r}{2}  \int_0^t |(h(\xi_u(x))^\top\Xi^{-1/2})^k|^2 du \right) \\
& \leq \exp \left( \frac{n r_2^2 p_r}{2}  \int_0^t |(h(\xi_u(x))^\top\Xi^{-1/2})^k|^2 du \right) 
\end{align*}
where the second equality holds due to independence of brownian increments and in the first inequality holds due to lemma \ref{lem:mgf}.  Combining this result with the same calculations as for $I_8$ yields an upper bound on $I_9$ which is identical to \eqref{eq:i8bound}.  Therefore, following the same reasoning as in for $I_8$, we have 
\begin{align*}
\left( \mathbb{E}^\mathbb{Q}[I_9(x)^{r_1}] \right)^{1/r_2} \leq (1 - L K(t))^{-\frac{1}{r_2}} 
\end{align*}
Finally, for $I_7$, first note that for $t \in [t_{j}, t_{j+1})$, 
\begin{align*}
\int_0^t h(\xi_u(x))^\top \Xi^{-1} \dot{Z}_u^d du &= \sum_{i=1}^{j} \int_{t_{i-1}}^{t_{i}} h(\xi_u(x))^\top \Xi^{-1}\left( h(x_{t_{i-1}}^\ast) + \frac{1}{\delta_d} \Xi^{1/2}(B_{t_i} - B_{t_{i-1}}) \right) du  \\
&+ \int_{t_{j}}^t h(\xi_u(x))^\top \Xi^{-1}\left( h(x_{t_{j}}^\ast) + \frac{1}{\delta_d} \Xi^{1/2}(B_{t_{j+1}} - B_{t_{j}}) \right)du 
\end{align*}
similarly, 
\begin{align*}
\int_0^t h(\xi_u(x))^\top \Xi^{-1} \circ dZ_u &= \int_0^t h(\xi_u(x))^\top \Xi^{-1}  dZ_u \\
&= \sum_{i=1}^j \int_{t_{i-1}}^{t_{i}} h(\xi_u(x))^\top \Xi^{-1} (h(x_u^\ast)du + \Xi^{1/2}dB_u) + \int_{t_j}^t h(\xi_u(x))^\top \Xi^{-1} (h(x_u^\ast)du + \Xi^{1/2}dB_u)
\end{align*}
Then
\begin{align*}
&|I_3 - I_4|^{pr_1} \leq |2(r-s)|^{pr_1}\left(  \left| \sum_{i=1}^{j} \int_{t_{i-1}}^{t_{i}} h(\xi_u(x))^\top \Xi^{-1} (h(x_{t_{i-1}}^\ast) - h(x_u^\ast)) du + \int_{t_{j}}^t h(\xi_u(x))^\top \Xi^{-1} (h(x_{t_{j}}^\ast) - h(x_u^\ast))du  \right|^{pr_1} \right.\\
& + \left. \left|  \int_0^{t} h(\xi_u(x))^\top \Xi^{-1/2} \dot{B}_u^d du -\int_0^{t} h(\xi_u(x))^\top \Xi^{-1/2} dB_u \right|^{pr_1^2} \right)  \\
& =: |2(r-s)|^{pr_1}\left( I_{11} + I_{12} \right)
\end{align*}
As we are required to upper bound $\mathbb{E}[|I_3 - I_4|^{pr_1} ]$, we obtain the following bound for $\mathbb{E}[I_{12}]$ following the same reasoning as in pg 41 of \cite{hu_approximation_2002} (bound on $I_5^\Pi(x)$ in their proof), 
\begin{align*}
\mathbb{E}[I_{12}] &= \mathbb{E}\left[ \left| \sum_{k=1}^n \int_0^{t} (h(\xi_u(x))^\top \Xi^{-1/2})^k \dot{B}_u^k du -\int_0^{t} (h(\xi_u(x))^\top \Xi^{-1/2})^k dB_u^k \right|^{pr_1} \right] \\
& \leq n^{pr_1} \sum_{k=1}^n \mathbb{E}\left[ \left|  \int_0^{t} (h(\xi_u(x))^\top \Xi^{-1/2})^k \dot{B}_u^k du -\int_0^{t} (h(\xi_u(x))^\top \Xi^{-1/2})^k dB_u^k \right|^{pr_1} \right] \\
& \leq C \frac{1}{\delta_d} \sum_{i=1}^j \int_{t_{i-1}}^{t_i} \int_{t_{i-1}}^{t_i} \mathbb{E}[ | (h(\xi_u(x)^\top \Xi^{-1/2})^k - (h(\xi_\tau(x))^\top\Xi^{-1/2})^k |^{pr_1} ]d\tau du\\
& + C\delta_d^{\frac{pr_1}{2}-1} \int_{t_j}^t \mathbb{E}[|(h(\xi_\tau(x))^\top\Xi^{-1/2})^k|^{pr_1}] d\tau 
\end{align*}
where $C$ is a constant depending on $n$.  It is worthwhile clarifing that classical Wong-Zakai/piecewise smooth convergence results are focused on stochastic integrals of the form $|\int h(X_s) \circ dB_s - \int h(X_s^d)\dot{B}_s^d ds|$ where the integrand $X_s^d$ is dependent on $B_s^d$.  As the coefficient here $h(\xi_s(x)$ evolves independently of $B$, we can resort to simpler convergence tools than used in e.g. \cite{pathiraja_l_2024}.  Starting with the second term, we have under the assumptions on $b,g,h$ using lemma \ref{lem:lipsde},
\begin{align*}
\mathbb{E}[|(h(\xi_\tau(x))^\top\Xi^{-1/2})^k|^{pr_1}] &\leq \mathbb{E}[|(h(\xi_\tau(x))^\top\Xi^{-1/2})^k - (h(x)^\top\Xi^{-1/2})^k|^{pr_1}] + \mathbb{E}[|(h(x)^\top\Xi^{-1/2})^k|^{pr_1}] \\
& \leq C(x) |\tau|^{pr_1/2} + C(x)
\end{align*}
so that 
\begin{align*}
C\delta_d^{\frac{pr_1}{2}-1} \int_{t_j}^t \mathbb{E}[|(h(\xi_\tau(x))^\top\Xi^{-1/2})^k|^{pr_1}] d\tau &\leq C\delta_d^{\frac{pr_1}{2}-1}  C(x)\int_{t_j}^t (|\tau|^{pr_1/2} +1) d\tau \\
& \leq  C\delta_d^{\frac{pr_1}{2}-1}  C(x)(\delta_d^{pr_1/2+1} + \delta_d) \\
& \leq C(x) \delta_d^{pr_1/2}
\end{align*}
Then for the first term on the rhs of the inequality, making use of lemma \ref{lem:lipsde} and that $h$ is lipschitz continuous and at most linear growth, we have
\begin{align*}
&C \frac{1}{\delta_d} \sum_{i=1}^j \int_{t_{i-1}}^{t_i} \int_{t_{i-1}}^{t_i} \mathbb{E}[ | (h(\xi_u(x)^\top \Xi^{-1/2})^k - (h(\xi_\tau(x))^\top\Xi^{-1/2})^k |^{pr_1} ]d\tau du \\
&\leq C(x) \frac{1}{\delta_d} \sum_{i=1}^j \int_{t_{i-1}}^{t_i} \int_{t_{i-1}}^{t_i} |\tau-u|^{pr_1/2}   d\tau du \\
&= \leq C(x) \frac{1}{\delta_d} \sum_{i=1}^j \delta^{pr_1/2} \int_{t_{i-1}}^{t_i} \int_{t_{i-1}}^{t_i}  d\tau du \\
& \leq C(x) \frac{1}{\delta_d} \sum_{i=1}^j  \delta^{pr_1/2 + 2} \\
& = C(x) \frac{1}{\delta_d} \frac{t_j}{\delta_d} \delta^{pr_1/2 + 2} = C(x) \delta_d^{pr_1/2}
\end{align*}
Combining the two yields 
\begin{align*}
\mathbb{E}[I_{12}(x)] \leq  C(x) \delta_d^{pr_1/2}
\end{align*}
Finally, for $I_{11}(x)$ we have using the notation $\lfloor x_t^\ast \rfloor = x_{t_i}^\ast, \enskip t \in [t_{i-1}, t_{i}]$ and Lemma \ref{lem:lipsde} and the linear growth assumption on $h$,
\begin{align*}
\mathbb{E}[I_{11}(x)] &= \mathbb{E}\left[\left| \int_0^{t_{j}} h(\xi_s(x))^\top \Xi^{-1} (h( \lfloor x_u^\ast \rfloor) - h(x_u^\ast)) du + \int_{t_{j}}^t h(\xi_u(x))^\top \Xi^{-1} (h(x_{t_{j}}^\ast) - h(x_u^\ast))du  \right|^{pr_1} \right] \\
& \leq C \int_0^{t_{j}} \mathbb{E}\left[ \left| h(\xi_u(x))^\top \Xi^{-1} (h( \lfloor x_u^\ast \rfloor) - h(x_u^\ast)) \right|^{pr_1} \right] du + C\int_{t_{j}}^t  \mathbb{E}\left[ \left| h(\xi_u(x))^\top \Xi^{-1} (h(x_{t_{j}}^\ast) - h(x_u^\ast)) \right|^{pr_1} \right] du \\
& \leq n^{pr_1}C(1+C_h(x)) \sum_{i=1}^{j} \int_{t_{i-1}}^{t_{i}} \mathbb{E}  \left[\left|  h(x_{t_{i-1}}^\ast) - h(x_u^\ast) \right|^{pr_1} \right] du + \int_{t_{j}}^{t} \mathbb{E} \left[ \left|  h(x_{t_{j}}^\ast) - h(x_u^\ast) \right|^{pr_1} \right] du \\
& \leq n^{pr_1}C(1+C_h(x)) \sum_{i=1}^{j+1} \int_{t_{i-1}}^{t_{i}} \delta_d^{pr_1/2} du \\
& \leq  n^{pr_1}C(1+C_h(x)) \frac{t_{j+1}}{\delta_d} \delta_d^{pr_1/2+1} \\
& \leq n^{pr_1}C(t)(1+C_h(x))  \delta_d^{pr_1/2}
\end{align*}
We are now ready to bound the remaining term involving $I_7$, i.e. 
\begin{align*}
\left(  \mathbb{E}^\mathbb{Q}[I_7(x)^{r_1}]\right)^{1/r_1} &\leq (2(r-s))^p \left(  \mathbb{E}^\mathbb{Q} [(\mathbb{E}[I_{11}(x)] + \mathbb{E}[I_{12}(x)])^{r_1}]
\right)^{1/r_1} \\
& \leq C (r-s)^p  \left( \mathbb{E}^\mathbb{Q} [C(x)^{r_1}] \delta_d^{pr_1^2/2}  \right)^{1/r_1} \\
& \leq C (r-s)^p \delta_d^{pr_1/2} 
\end{align*}
under the assumption of finite $p > 1$ moments of the initial density $f(x)$.  Choosing $r_1 >1$ large enough for a given $p$ gives us the required decay as $\delta_d \rightarrow 0$.

\end{proof}

\subsection{Proof of Lemma \ref{lem:meanfieldproc}}
\label{sec:proofmeanfieldproc}

\begin{proof}
We start with the proof of the first claim for the Stochastic EnKbf.  Note that  \eqref{eq:mfinfla} takes the form 
\begin{align*}
d\bar{X}_t = G\bar{X}_t dt + \mu_t(\bar{X}_t)dt + \Sigma^{1/2}dW_t + \sigma_t d\bar{B}_t + \sigma_t^z dZ_t 
\end{align*}
with 
\begin{align*}
\sigma_t &= \sqrt{2(r-s)} C_t H^\top \Xi^{-1/2}; \\
\mu_t(x) &=  -\frac{s}{2}C_t H^\top \Xi^{-1}H(x - m_t) - (r-s)C_t H^\top \Xi^{-1}Hx  \\
\sigma^z_t &= (r-s) C_t H^\top \Xi^{-1} 
\end{align*}
and denote by $D_t:= \frac{1}{2} \sigma \sigma^\top$.   The (observation conditioned) forward Kolmogorov equation (see e.g. \cite{pathiraja_mckean-vlasov_2021}) is then given by 
\begin{align}
\label{eq:forwardkolmrough} 
\partial_t \rho_t(x) =  \mathcal{L}^\ast \rho_t(x)  \underbrace{-\nabla \cdot (\rho_t(x)\mu_t(x)) - \nabla \rho_t(x)  \cdot \sigma_t^z dZ_t  + \nabla \cdot (D_t \nabla \rho_t(x))}_{=: \mathcal{M}^\ast \rho_t(x)}  
\end{align} 
where 
\begin{align}
\label{eq:gensig} 
\mathcal{L}^\ast \rho_t(x) = -\nabla \cdot (\rho_t(x) Gx) + \nabla \cdot (\Sigma \nabla \rho_t(x)) 
\end{align}
characterises the evolution under the signal process.   The remaining term can be simplified further using that  $\rho_t(x) = \mathcal{N}(x; m_t ,C_t) \enskip  \forall \enskip t \geq 0$ with $m_t, C_t$ to be determined, 
\begin{align*}
\mathcal{M}^\ast \rho_t(x) &= \rho_t(x) \text{Tr}\left[(x-m_t)^\top C_t^{-1} D_t C_t^{-1} (x - m_t)  - D_t C_t^{-1}  + (x-m_t)^\top C_t^{-1} \mu_t(x)   \right] - \rho_t(x)  \nabla \cdot \mu_t (x) \\
& +  (x-m_t)^\top C_t^{-1}\rho_t(x) \sigma_t^z dZ_t  \\
& =  \rho_t(x)  \frac{s}{2} \text{Tr} \left[  H^\top C_t H \Xi^{-1} \right]    + (r-s) (Hx-Hm_t)^\top \Xi^{-1} dZ_t \rho_t(x) \\
& - \frac{s}{2} \rho_t(x)  (Hx-Hm_t)^\top   \Xi^{-1}  (Hx - Hm_t)    - (r-s) \rho_t(x) (x-m_t)^\top H^\top \Xi^{-1} H m_t 
\end{align*}
using the cyclic property of the trace.  Substituting back into \eqref{eq:forwardkolmrough} yields an equation formally an equation formally equivalent to \eqref{eq:modifiedkushner}  since
\begin{align*}
\mathbb{E}_{p_t}[\| Hx - Hm_t\|_\Xi^2] =  \text{Tr} \left[  H^\top C_t H \Xi^{-1} \right].  
\end{align*}
The proof of the second claim follows from a similar line of reasoning with $\sigma_t = 0$ and $\mu_t(x) =   \left(\frac{r}{2} - s \right)C_t H^\top \Xi^{-1}H(x - m_t) - \frac{(r-s)}{2}C_t H^\top \Xi^{-1}H(x + m_t)$ and is therefore omitted.

The claim involving the deterministic update and piecewise smooth observation \eqref{eq:mfdetode} follows form a very similar line of reasoning. Nevertheless, we provide it here for completeness.  Note that \eqref{eq:mfdetode} can be written in the form, 
\begin{align*}
d\bar{X}_t = G\bar{X}_t dt + \mu_t(\bar{X}_t)dt + \Sigma^{1/2}dW_t 
\end{align*}
with 
\begin{align*}
\mu_t(x)  = -\left(\frac{s}{2} +\frac{(r-s)}{2} \right) C_t H^\top \Xi^{-1} H x + (r-s) C_t H^\top \Xi^{-1} \frac{dZ_t^d}{dt} + \left(\frac{s}{2} - \frac{r-s}{2} \right) C_t H^\top  \Xi^{-1} H m_t 
\end{align*}
The (conditional) forward Kolmgorov equation is then given by 
\begin{align}
\label{eq:forwardkolmdet}
\partial_t \rho_t(x) =  \mathcal{L}^\ast \rho_t(x)  \underbrace{-\nabla \cdot (\rho_t(x)\mu_t(x))}_{=: \mathcal{M}^\ast \rho_t}  
\end{align}
where $\mathcal{L}^\ast$ is as given in \eqref{eq:gensig}.  We have that  
\begin{align*}
\mathcal{M}^\ast \rho_t(x) &= \rho_t(x) \text{Tr}\left[  (x-m_t)^\top C_t^{-1} \mu_t(x)   \right] -  \rho_t(x) \nabla \cdot \mu_t (x) \\
& = \left(\frac{s}{2} +\frac{(r-s)}{2} \right) \rho_t(x) \mathbb{E}_{p_t}[\| Hx - Hm_t\|_\Xi^2] - \left(\frac{s}{2} +\frac{(r-s)}{2} \right) (Hx-Hm_t)^\top  \Xi^{-1} H x \rho_t(x)  \\
&  + \left(   \left(\frac{s}{2} - \frac{r-s}{2} \right) C_t H^\top   \Xi^{-1} H m_t  \right)^\top C_t^{-1} (x-m_t) \rho_t(x) + (r-s) (Hx-Hm_t)^\top \Xi^{-1} \frac{dZ_t^d}{dt} \rho_t(x) \\
& = (r-s)\rho_t(x) \left( -\frac{1}{2} x^\top H^\top \Xi^{-1} H x  + (Hx - Hm_t)^\top \Xi^{-1} \frac{dZ_t^d}{dt}  \right) -\frac{s}{2} \rho_t(x) \left( x^\top H^\top \Xi^{-1} H x - \mathbb{E}_{p_t}[\| Hx - Hm_t\|_\Xi^2]  \right) \\
& -\frac{(r-s)}{2} \mathbb{E}_{p_t}[\| Hx - Hm_t\|_\Xi^2]  +\left(\frac{s}{2} +\frac{(r-s)}{2} \right) (Hm_t)^\top  \Xi^{-1} H x \rho_t(x) \\
& + \left(   \left(\frac{s}{2} - \frac{r-s}{2} \right) C_t H^\top   \Xi^{-1} H m_t  \right)^\top C_t^{-1} (x-m_t) \rho_t(x) \\
& = (r-s)\rho_t(x) \left( -\frac{1}{2} x^\top H^\top \Xi^{-1} H x  + (Hx - Hm_t)^\top \Xi^{-1} \frac{dZ_t^d}{dt}  \right) -\frac{s}{2} \rho_t(x) \left( x^\top H^\top \Xi^{-1} H x - \mathbb{E}_{p_t}[\| Hx - Hm_t\|_\Xi^2]  \right) \\
& -\frac{(r-s)}{2} \mathbb{E}_{p_t^d} \left[ x^\top H^\top \Xi^{-1} Hx \right] \rho_t(x) + s \rho_t(x) (Hx)^\top \Xi^{-1} H m_t - \frac{s}{2} \rho_t(x) (Hm_t)^\top \Xi^{-1} Hm_t 
\end{align*}
once again using that 
\begin{align*}
\mathbb{E}_{\rho_t}[\| Hx - Hm_t\|_\Xi^2] =  \text{Tr} \left[  H^\top C_t \Xi^{-1}H \right].  
\end{align*}
Substituting back into  \eqref{eq:forwardkolmdet} yields an equation formally equivalent to \eqref{eq:repmutlineargauss}. 

For the proof of the final claim, i.e. for stochastic update with piecewise smooth observations $\{Z_t^d\}$, note that \eqref{eq:mfstochode} takes the form 
\begin{align*}
d\bar{X}_t = G\bar{X}_t dt + \mu_t(\bar{X}_t)dt + \Sigma^{1/2}dW_t + \sigma_t d\bar{B}_t 
\end{align*}
with 
\begin{align*}
\sigma_t &= \sqrt{r-s} C_t H^\top \Xi^{-1/2}; \\
\mu_t(x) &=   -\frac{s}{2} C_t H^\top \Xi^{-1}H(x - m_t) + (r-s) C_t H^\top \Xi^{-1} \left(\frac{dZ_t^d}{dt} - Hx  \right)
\end{align*}
and denote by $D_t := \frac{1}{2}  \sigma_t \sigma_t^\top $.   The (conditional) forward Kolmgorov equation is then given by 
\begin{align}
\label{eq:forwardkolm}
\partial_t \rho_t(x) =  \mathcal{L}^\ast \rho_t(x)  \underbrace{-\nabla \cdot (\rho_t(x)\mu_t(x)) + \nabla \cdot (D_t \nabla \rho_t(x))}_{=: \mathcal{M}^\ast \rho_t(x)}  
\end{align}
where $\mathcal{L}^\ast$ is as given in \eqref{eq:gensig}.
The remaining term can be simplified further using that $\rho_t(x) = \mathcal{N}(x; m_t ,C_t) \enskip  \forall \enskip t \geq 0$ with $m_t, C_t$ to be determined, 
\begin{align*}
\mathcal{M}^\ast \rho_t(x) &= \rho_t(x) \text{Tr}\left[(x-m_t)^\top C_t^{-1} D_t C_t^{-1} (x - m_t)  - D_t C_t^{-1}  + (x-m_t)^\top C_t^{-1} \mu_t(x)   \right] -  \rho_t(x) \nabla \cdot  \mu_t (x) \\
& =  \rho_t(x)  \frac{r}{2} \text{Tr} \left[  H^\top C_t \Xi^{-1}H \right]    + (r-s) (Hx-Hm_t)^\top \Xi^{-1} \frac{dZ_t^d}{dt} \rho_t(x) \\
& + \left(\frac{r}{2} - s \right) \rho_t(x)  (Hx-Hm_t)^\top   \Xi^{-1}  (Hx - Hm_t)    - (r-s) \rho_t(x) (x-m_t)^\top H^\top \Xi^{-1} H x \\
& = (r-s) (Hx-Hm_t)^\top \Xi^{-1} \frac{dZ_t^d}{dt} \rho_t(x) - \frac{(r-s)}{2} \text{Tr}[Hx \Xi^{-1} (Hx)^\top - \mathbb{E}_{\rho_t}[Hx \Xi^{-1}(Hx)^\top ] ]  \rho_t(x)  \\
& -\frac{s}{2} \text{Tr}[  (Hx - Hm_t)^\top \Xi^{-1} (Hx - Hm_t)  -H^\top C_t \Xi^{-1} H] \rho_t(x)
\end{align*}
using the cyclic property of the trace and the identity 
\begin{align*}
\text{Tr}[\mathbb{E}_{\rho_t}[x x^\top] - x x^\top  ] = \text{Tr}[(x-m_t)(x-m_t)^\top -2 x (x - m_t)^\top + C_t ] . 
\end{align*}
Inserting back into \eqref{eq:forwardkolm} yields a formal equivalence to \eqref{eq:repmutlineargauss}, as desired.

\end{proof}

\subsection{Proof of Lemma \ref{lem:covreperror}}
\label{sec:proofcovreperror}

\begin{proof}
Starting from the evolution equation for the mean \eqref{eq:meanlinear}, we have 
\begin{align}
\nonumber 
d \varepsilon_t &= dm_t - dX_t^\ast \\
\label{eq:sdevareps}
& = (G -(r-s)K_tH)\varepsilon_t dt  - b dt + \Sigma^{1/2}dW_t + (r-s)K_t \Xi^{1/2}dB_t 
\end{align}
where $x_t^\ast$ denotes the reference process that generated the observation path, from which we obtain the evolution of the expected error 
\begin{align}
\label{eq:odeerror}
d\mathbb{E}[\varepsilon_t] = (G -(r-s)K_tH)\mathbb{E}[\varepsilon_t] dt  - b dt. 
\end{align}
Also, 
\begin{align*}
d\mathbb{E}[\varepsilon_t]\mathbb{E}[\varepsilon_t^\top] &=  \mathbb{E}[\varepsilon_t] d \mathbb{E}[\varepsilon_t^\top] + d\mathbb{E}[\varepsilon_t]  \mathbb{E}[\varepsilon_t^\top] \\
& = \mathbb{E}[\varepsilon_t] \mathbb{E}[\varepsilon_t^\top](G -(r-s)K_tH)^\top dt +  (G -(r-s)K_tH)\mathbb{E}[\varepsilon_t] \mathbb{E}[\varepsilon_t^\top] dt  -  \mathbb{E}[\varepsilon_t] b^\top dt - b\mathbb{E}[\varepsilon_t^\top]dt 
\end{align*}
Since  $\tilde{P}_t = \mathbb{E}[\varepsilon_t \varepsilon_t^\top]$ and using Ito formula, we have   
\begin{align*}
d\tilde{P}_t &=  \mathbb{E}[\varepsilon_t (d\varepsilon_t)^\top] + \mathbb{E}[(d\varepsilon_t) \varepsilon_t^\top] + \mathbb{E}[d \varepsilon_t (d \varepsilon_t)^\top]\\
& =  \mathbb{E}[\varepsilon_t \varepsilon_t^\top (G -(r-s)K_tH)^\top dt  - \varepsilon_t b^\top dt] + \mathbb{E}[ (G -(r-s)K_tH)\varepsilon_t \varepsilon_t^\top dt - b \varepsilon_t^\top dt] + \Sigma dt + (r-s)^2 K_t \Xi K_t^\top dt \\
& = \tilde{P}_t(G -(r-s)K_tH)^\top dt + (G -(r-s)K_tH)\tilde{P}_t dt - \mathbb{E}[\varepsilon_t]b^\top dt - b \mathbb{E}[\varepsilon_t^\top] dt + \Sigma dt + (r-s)^2 K_t \Xi K_t^\top dt
\end{align*}
which yields \eqref{eq:odeptilde}.  
Using, $P_t = \tilde{P}_t - \mathbb{E}[\varepsilon_t]\mathbb{E}[\varepsilon_t^\top]$ and combining the above expressions yields  \eqref{eq:odeforPt} and \eqref{eq:odeptilde}. 

For the second claim, consider that when $s=0$, we have 
\begin{align*}
dP_t &= (G -rK_tH) P_t dt + P_t (G -rK_tH)^\top dt + \Sigma dt  + r^2 K_t \Xi K_t^\top dt 
\end{align*}
whose solution at any given time $t$ is formally equivalent to $C_t$ (whose time evolution is given by \ref{eq:covlinear}) if $P_0 = C_0$ and $r=1$, since $K_t \Xi K_t^\top = K_t H C_t$.  Notice that due to the $r^2K_t \Xi K_t^\top$ rather than $rK_t \Xi K_t^\top$ term, this equivalence only holds when $r=1$ in addition to $s=0$.\\     
\\
For the final claim, since $Tr(\tilde{P}_t) = E_t$, using the cyclic property of the trace yields 
\begin{align}
\label{eq:eqntildePt}
\frac{d Tr(\tilde{P}_t)}{dt} = \text{Tr}[(A_t + A_t^\top) \tilde{P}_t]   - 2\text{Tr} \left[ \mathbb{E}[\varepsilon_t]b^\top \right]   + \text{Tr}[\Sigma]  + (r-s)^2 \text{Tr}[K_t \Xi K_t^\top]
\end{align}
where $A_t(r,s):= G - (r-s)K_tH$ and we have used the cyclic property of the trace.  Then using $\text{Tr}[AB] \leq \alpha(A) \text{Tr}[B]$ for any real symmetric matrix $A$ and p.d. matrix $B$, where $\alpha(A)$ denotes the largest eigenvalue of real matrix $A$, we have 
\begin{align*}
\text{Tr}[(A_t + A_t^\top)P_t] &\leq \alpha(A_t + A_t^\top) \text{Tr}[\tilde{P}_t] \\
\text{Tr}[K_t \Xi K_t^\top]  = \text{Tr}[H^\top \Xi^{-1} H C_t^\top C_t] &\leq \alpha(H^\top \Xi^{-1} H) \|C_t\|_F 
\end{align*}
where $\|A \|_F = \text{Tr}[A^\top A]$ is the Frobenius norm.  Substituting into \eqref{eq:eqntildePt} yields \eqref{eq:odeforEt}.


\end{proof}

\subsection{Proof of Lemma \ref{lem:asympbiasvar}}
\label{sec:proofasympbiasvar}

\begin{proof}
Under assumption \ref{ass:Astable}, $C_\infty$ exists and is unique (see e.g. \cite{bishop_mathematical_2023}).  When $C_0 = C_\infty$ we have that $K_t = K_\infty = C_\infty H^\top \Xi^{-1}$ for all $t \geq 0$.  Using \eqref{eq:expsolnerror} and the shorthand notation $A_\infty(r,s) := G - (r-s)K_\infty H$, 
\begin{align}
\label{eq:expsolnerror}
\mathbb{E}[\varepsilon_t] = \exp(t A_\infty(s,r) )\mathbb{E}[\varepsilon_0] + A_\infty^{-1}(s,r) [I - \exp\left(t A_\infty(s,r)\right)]b 
\end{align}
When $\alpha(A_\infty(r,s)) < 0$, for a given $r,s$, it holds that as $t \rightarrow \infty$, $\enskip \mathbb{E}[\varepsilon_t] \rightarrow -A_\infty(r,s)^{-1} b$, from which we obtain \eqref{eq:asympbias}.\\
\\
For brevity, we drop the dependence on $r,s$ in $A_\infty(r,s)$ in the remainder of the proof.  For the second claim, rewrite  \eqref{eq:odeptilde} as 

\begin{align*}
\frac{d\tilde{P}_t}{dt} &= A_\infty \tilde{P}_t  + \tilde{P}_t A_\infty^\top + \mathcal{B}_t \\
\mathcal{B}_t &:= \Sigma  + (r-s)^2 K_\infty \Xi K_\infty^\top   - \mathbb{E}[\varepsilon_t]b^\top  - b \mathbb{E}[\varepsilon_t^\top] 
\end{align*}
which has explicit solution 
\begin{equation} \label{eq:expptilde1}
\begin{split}
\tilde{P}_t &= \exp\left( t  A_\infty\right) \tilde{P}_0 \exp\left( t  A_\infty^\top \right)  + \int_0^t \exp\left( (t-u)A_\infty \right) (\Sigma + (r-s)^2 K_\infty \Xi K_\infty^\top ) \exp\left( (t-u)A_\infty^\top \right)  \d u \\
&- \int_0^t \exp\left( (t-u)A_\infty \right) \left(\mathbb{E}[\varepsilon_u] b^\top + b \mathbb{E}[\varepsilon_u^\top]\right) \exp\left( (t-u)A_\infty^\top \right)  \d u. 
\end{split}
\end{equation}
The final term can be simplified further, using \eqref{eq:expsolnerror} and the assumption $\mathbb{E}[\varepsilon_0] = 0$,
\begin{align*}
&-\int_0^t \exp\left( (t-u)A_\infty \right) \mathbb{E}[\varepsilon_u] b^\top \exp\left( (t-u)A_\infty^\top \right)  \d u \\
& = -\int_0^t \exp\left( (t-u)A_\infty \right)  A_\infty^{-1}bb^\top \exp\left( (t-u)A_\infty^\top \right)  \d u + \exp\left( tA_\infty \right)  A_\infty^{-1} bb^\top \int_0^t  \exp\left( (t-u)A_\infty^\top \right)  \d u.
\end{align*}
and the remaining term involving $b \mathbb{E}[\varepsilon_u^\top]$ is the transpose of the above.  Evaluating the last integral in the above and substituting into \eqref{eq:expptilde1} and taking trace of both sides yields 
\begin{align}
\nonumber 
\text{Tr}[\tilde{P}_t] &= \text{Tr}[\exp\left( t  A_\infty\right) \tilde{P}_0 \exp\left( t  A_\infty^\top \right)] +2\text{Tr}[\exp(tA_\infty)A_\infty^{-1} b b^\top \exp(tA_\infty^\top) (I - \exp(-tA_\infty^\top)) (A_\infty^{-1})^\top]  \\
\label{eq:traceptilde}
&+ \int_0^t \text{Tr}[\exp\left( (t-u)A_\infty \right) \tilde{\Sigma} \exp\left( (t-u)A_\infty^\top \right)]  \d u 
\end{align}
where $\tilde{\Sigma} := \Sigma + (r-s)^2 K_\infty \Xi K_\infty^\top  - A_\infty^{-1} b b^\top - b b^\top (A_\infty^{-1})^\top$.  It holds that the all terms except for the integral term in \eqref{eq:traceptilde} go to zero as $t \rightarrow \infty$ when $\alpha(A_\infty) < 0$. 
For the integral term, by the cyclic property of the trace,  
\begin{align}
\int_0^t \text{Tr}[\exp\left( (t-u)A_\infty \right) \tilde{\Sigma} \exp\left( (t-u)A_\infty^\top \right)]  \d u & = \text{Tr}\left[\tilde{\Sigma}\int_0^t  \exp\left( vA_\infty \right) \exp\left( vA_\infty^\top   dv\right) \right]
\end{align}
Furthermore, when $\alpha(A_\infty) < 0$, 
\begin{align*}
\int_0^\infty \exp(v A_\infty) \exp(vA_\infty^\top)dv = X_\infty
\end{align*}
where $X_\infty$ is a symmetric matrix that is the unique solution of the Lyapunov equation 
\begin{align}
\label{eq:lyapunov}
A_\infty^\top X_\infty + X_\infty A_\infty + I = 0 
\end{align}
Combining yields the result \eqref{eq:asympmerror}. 
\\
For the upper bound on $E_\infty$, consider the differential inequality in \eqref{eq:odeforEt}, recall (and using assumption $\mathbb{E}[\varepsilon_0] = 0$), 
\begin{align*}
\frac{dE_t}{dt} \leq \alpha(A_\infty + A_\infty^\top) E_t + \beta_t + \gamma 
\end{align*}
where 
\begin{align*}
    \beta_t &:= 2\text{Tr}[A_\infty^{-1}  \exp(tA_\infty)b b^\top]  \\
    \gamma &:= - 2 \text{Tr}[A_\infty^{-1}b b^\top] + \text{Tr}[\Sigma] + (r-s)^2 \lambda_{max}(H^\top \Xi^{-1}H) \|C_\infty\|_F. 
\end{align*}
By Gronwall lemma, 
\begin{align*}
E_t \leq  e^{\alpha(A_\infty + A_\infty^\top) t} \left[ E_0 + \int_0^t \beta_s e^{-\alpha(A_\infty + A_\infty^\top) s} ds \right]  + \frac{\gamma}{\alpha(A_\infty + A_\infty^\top) } \left( e^{\alpha (A_\infty + A_\infty^\top) t}  - 1 \right)  
\end{align*}
using the cyclic property of the trace, 
\begin{align*}
    \int_0^t \text{Tr}[A_\infty^{-1} \exp(s A_\infty)b b^\top] e^{-\alpha(A_\infty + A_\infty^\top) s} ds & =  \text{Tr} \left[ b b^\top A_\infty^{-1} \int_0^t \exp(s A_\infty) e^{-\alpha(A_\infty + A_\infty^\top) s} ds \right] \\
    & = \text{Tr} \left[ b b^\top A_\infty^{-1} (A_\infty - \alpha(A_\infty + A_\infty^\top) I)^{-1} (\exp(t (A_\infty - \alpha(A_\infty + A_\infty^\top) I))) \right] 
\end{align*}
where the final equality holds since $A_\infty - \alpha(A_\infty + A_\infty^\top) I$ is invertible. Combining yields 
\begin{align*}
    E_t & \leq  e^{\alpha(A_\infty + A_\infty^\top) t} E_0 +  2\text{Tr} \left[ b b^\top A_\infty^{-1} (A_\infty - \alpha (A_\infty + A_\infty^\top)I)^{-1}  \exp(t A_\infty) \right]   + \frac{\gamma}{\alpha(A_\infty + A_\infty^\top) } \left( e^{\alpha (A_\infty + A_\infty^\top)t}  - 1 \right)
\end{align*}
Taking the limit $t \rightarrow \infty$ yields the result \eqref{eq:boundEinf} due to assumptions \ref{ass:Astable} and \ref{ass:eigsymmA}.

For the final claim, in the scalar case $m=n=1$, 
\eqref{eq:steadystatecov} is explicitly solvable with 
\begin{align}
\label{eq:asympCscalar}
C_\infty = \frac{G + \sqrt{G^2 + rH^2 \Xi^{-1} \Sigma}}{r H^2 \Xi^{-1}}
\end{align}
which when substituted into $K_\infty = C_\infty H \Sigma^{-1}$ yields  
\begin{align}
\label{eq:Ainfexpress}
A_\infty &= G - \frac{(r-s)}{r} \left(  G + \sqrt{G^2 + rH^2 \Xi^{-1} \Sigma} \right)   \\
\nonumber 
& = \frac{s}{r}G - \frac{(r-s)}{r} \sqrt{G^2 + rH^2 \Xi^{-1} \Sigma}. 
\end{align}
We can also trivially solve the Lyapunov equation \eqref{eq:lyapunov} to obtain 
\begin{align*}
X_\infty = -\frac{1}{2}A_\infty^{-1}. 
\end{align*}
Substituting the above expressions into \eqref{eq:asympmerror} yields \eqref{eq:asympmsescalar}.

\end{proof}

\subsection{Proof of Lemma \ref{lem:optimalrs} }    
\label{sec:proofoptimalrs}

\begin{proof}

In the below we drop the $(s,r)$ dependency in $A_\infty(s,r)$ where there is no ambiguity. Recall also that in the scalar case $E_t = \tilde{P}_t$, and we will use the notation $\tilde{P}_\infty$ in the below to refer to the asymptotic MSE.

We start with the proof of the first claim, i.e. the optimal $s$ for any given $r >0$ minimising  \eqref{eq:asympmsescalar} by solving 
\begin{align*}
\frac{\partial \tilde{P}_\infty}{\partial s}  = \frac{\partial \tilde{P}_\infty}{\partial A_\infty} \frac{\partial A_\infty}{\partial s} = 0 
\end{align*}
where     
\begin{align*}
\frac{\partial \tilde{P}_\infty}{\partial A_\infty} &= -2b^2 A_\infty^{-3} + A_\infty^{-1} \left( \frac{G-A_\infty}{H} \right) \frac{\Xi}{H} + \frac{1}{2}A_{\infty}^{-2} \left( \Sigma + \left( \frac{G - A_\infty}{H} \right)^2 \Xi \right) \\
\frac{\partial A_\infty}{\partial s} &=   \frac{G}{r} + \frac{\sqrt{G^2 + rH^2 \Xi^{-1} \Sigma}}{r}. 
\end{align*}
Since $\frac{\partial A_\infty}{\partial s}$ is clearly never zero, we need only solve 
\begin{align}
\label{eq:depressedcubic}
g(A_\infty) := A_\infty^3 - \frac{(H^2 \Sigma + \Xi G^2)}{\Xi} A_\infty + \frac{4b^2H^2}{\Xi} = 0
\end{align}
which takes the form of a depressed cubic $A_\infty^3 + pA_\infty + q = 0$ with $p = -(H^2 \Xi^{-1}\Sigma +  G^2), \enskip q = 4b^2H^2 \Xi^{-1}$, and we have that $q >0, p < 0$ always.  The nature of the roots can be characterised in the usual way via the discriminant $\tau$, i.e. when $\tau > 0$, \eqref{eq:depressedcubic} has one real root and two complex roots, whilst when $\tau < 0$, it has three real roots. We deal with these two cases below separately.\\
\\
\textbf{Case 1: $\tau > 0$.}   First we show that the real root here is strictly negative.  From inspection, we have that $g(0) > 0$.  Also, \eqref{eq:depressedcubic} always has two real extremal (turning) points, one of which is strictly negative and the other is strictly positive, since $g'(A_\infty) = 3A_\infty^2 + p = 0$ implies the extremal points occur at $A_\infty = \pm \sqrt{\frac{-p}{3}}$.  Furthermore, $g'(0) = p < 0$.  Combining these properties implies that \eqref{eq:depressedcubic} has one negative real root when $\tau > 0$.  
When $\tau > 0$, we can use Cardano's formula to obtain an expression for the real root, 
\begin{align*}
A_\infty^\ast = \left(-\frac{q}{2} + \sqrt{\tau}\right)^{1/3} + \left(-\frac{q}{2} - \sqrt{\tau}\right)^{1/3}
\end{align*}
from which we then obtain \eqref{eq:Ainfstarexp}, 
and substituting the expression for $A_\infty^\ast$ into \eqref{eq:Ainfscalar} and re-arranging yields \eqref{eq:optsr}.  It can be verified straightforwardly that this this is indeed a minimum since $\frac{\partial^2 \tilde{P}_\infty}{\partial A_\infty^2}(A_\infty^\ast) > 0$.  \\ 
\\
\textbf{Case 2: $\tau < 0$.}  In this case, \eqref{eq:depressedcubic} has three real roots due to the usual condition on the discriminant.  First we characterise the signs of the roots.  We have directly from \eqref{eq:depressedcubic} that $g(0) > 0$.  Also, \eqref{eq:depressedcubic} always has two real extremal (turning) points, one of which is strictly negative and the other is strictly positive, since $g'(A_\infty) = 3A_\infty^2 + p = 0$ implies the extremal points occur at $A_\infty = \pm \sqrt{\frac{-p}{3}}$.  Furthermore, $g'(0) = p < 0$.  Combining these properties implies that \eqref{eq:depressedcubic} has two positive real roots and one negative real root when $\tau < 0$.\\
\\
When $\tau < 0$, the following trigonometric formula holds for the characterisation of the three real roots, $A_\infty^k, \enskip k = 0, 1, 2$
\begin{align*}
A_\infty^k = 2\,\sqrt{-\frac{p}{3}}\,\cos\left[\,\frac{1}{3} \cos^{-1}\left(\frac{3q}{2p}\sqrt{\frac{-3}{p}}\,\right) - \frac{2\pi k}{3}\,\right] \qquad \text{for } k=0,1,2
\end{align*} 
To determine the negative root, first notice that $-\frac{2\pi k}{3} < \frac{1}{3}\cos^{-1}(y) -\frac{2\pi k}{3} < -\frac{2\pi k}{3} + \frac{\pi}{3}$ for all $-1 < y < 1$. Furthermore, $\cos(\theta) < 0$ for all $-\frac{3 \pi}{2} < \theta < -\frac{\pi}{2}$ and it holds that for $\theta := \frac{1}{3}\cos^{-1}(y) -\frac{2\pi k}{3}$ and $k=2$, $-\frac{4}{3}\pi  < \theta < -\pi$, so that $\cos(\theta) < 0$ for $k = 2$.  This yields the remaining case in \eqref{eq:Ainfstarexp}.\\
\\
For the remaining claims, we work with the following change of variable
\begin{align}
	\label{eq:covyr}
	y:= \sqrt{G^2 +r H^2 \Xi^{-1} \Sigma}
\end{align}
which when substituted into \eqref{eq:optsr} with $r = \frac{y^2 - G^2}{H^2\Xi^{-1} \Sigma}$ yields
\begin{align}
    \label{eq:sopty}
	s^{\text{opt}}(y) =  \frac{1}{H^2\Xi^{-1} \Sigma} (y-G)(A_\infty^\ast +y).
\end{align}
We may now characterise the admissible values of $s$ for a fixed $r$.  Since we require $A_\infty < 0$ to ensure the existence of $\tilde{P}_\infty$, any choice of $s$ must satisfy 
\begin{align*}
	\frac{s}{r}G - \frac{(r-s)}{r} \sqrt{G^2 + rH^2 \Xi^{-1} \Sigma} < 0
\end{align*}
from which the upper bound in \eqref{eq:condsons} immediately follows (along with the requirement that $s < r$).The minimiser of $s^{\text{opt}}(y)$ wrt $y$, which we donte by $y^\ast$, can be found straightforwardly by solving $\frac{d s^{\text{opt}}(y)}{dy} = 0$, giving 
\begin{align}
    \label{eq:yopt}
	y^\ast = \frac{G -A_\infty^\ast}{2}
\end{align}
which implies that for any $r >0$,
\begin{align}
	\label{eq:sopty}
	s^{\text{opt}} \geq -\frac{(G + A_\infty^\ast)^2}{4 H^2 \Xi \Sigma} =: s^{l}
\end{align}
since $y$ is monotonically related to $r$ (notice also that $s^l$ is independent of $r$). 
 Note however that by definition, $y > |G|$ since $r > 0$.  Therefore the lower bound on admissible $s$ values corresponds to $s^l$ if $G -2|G| < A_\infty^\ast$, otherwise it corresponds to $s^{\text{opt}}(y=|G|)$, which we denote by $s^u$ in the lemma.  \\  
\\
Finally, we can obtain expressions for $r^{\text{opt}}$ for a given $s$ satisfying the aforementioned constraints by a straightforward rearrangement of \eqref{eq:sopty}, 
\begin{align}
	\label{eq:quady}
	y^2 + (A_\infty^\ast - G)y - G A_\infty^\ast - sH^2 \Xi^{-1}\Sigma = 0,
\end{align}
yielding 
\begin{align}
	\label{eq:solnopty} 
	y = \frac{(G - A_\infty^\ast) \pm \sqrt{(G - A_\infty^\ast)^2 + 4(G A_\infty^\ast + s H^2 \Xi^{-1} \Sigma)} }{2}.
\end{align}
Note that $s^{\text{opt}}(y)$ is a convex quadratic in $y$, with roots $y = G$ and $y = -A_\infty^\ast$ and minimum attained at $y = y^\ast$ as defined in \eqref{eq:yopt}. Also, it holds that $y > |G|$ which follows directly from the definition \eqref{eq:covyr} and that $r > 0$. Therefore, whenever $y^\ast < |G|$, there is a unique admissible $y$ value for any given $s$ satisfying the aforementioned constraints.  More specifically, this value is given by 
\begin{align}
    \label{eq:ypos}
    y = \frac{(G - A_\infty^\ast) + \sqrt{(G - A_\infty^\ast)^2 + 4(G A_\infty^\ast + s H^2 \Xi^{-1} \Sigma)} }{2}
\end{align}
Finally, whenever $y^\ast \geq |G|$, there are two admissible $y$ values for any given $s^l < s < s^{\text{opt}}(y = |G|)$, and a unique $y$ value for $s > s^{\text{opt}}(y = |G|)$ given by \eqref{eq:ypos}.  Substituting the definition of $y$ into the above expressions for $y$ and rearranging yields the expressions for $r^{\text{opt}}$.  The final claim is now proved.

\end{proof}

\subsection{Proof of Lemma \ref{lem:covimpact}}
\label{sec:proofcovimpact}

\begin{proof}
The first claim follows from a simple re-arrangement of the identity $C_\infty = \tilde{P}_\infty$ with $A_\infty = A_\infty^\ast$, i.e. 
\begin{align*}
- \frac{1}{2}\left(\Sigma + \left( \frac{G - A_\infty^\ast}{H} \right)^2 \Xi  \right)  \frac{1}{A_\infty^\ast}  + \left(\frac{b}{A_\infty^\ast} \right)^2 = \frac{G + \sqrt{G^2 + rH^2 \Xi^{-1} \Sigma}}{rH^2 \Xi^{-1}}
\end{align*}
additionally, rearranging \eqref{eq:Ainfexpress} yields  
\begin{align*}
G + \sqrt{G^2 + rH^2 \Xi^{-1} \Sigma} = \frac{r(G - A_\infty^\ast)}{(r-s)} 
\end{align*}
substituting into the expression for $\tilde{P}_\infty = C_\infty$ yields the result. 
Additionally, we have 
\begin{align*}
r = \frac{1}{4H^2 \Xi^{-1} \Sigma} \left( G - A_\infty^\ast \pm \sqrt{(G - A_\infty^\ast)^2 + 4(G A_\infty^\ast + sH^2 \Xi^{-1} \Sigma)} \right)^2 - \frac{G^2}{H^2 \Xi^{-1} \Sigma}
\end{align*}
so that 
\begin{align*}
r-s &= \frac{1}{4H^2 \Xi^{-1} \Sigma} \left( G - A_\infty^\ast \pm \sqrt{(G - A_\infty^\ast)^2 + 4(G A_\infty^\ast + sH^2 \Xi^{-1} \Sigma)} \right)^2 - \frac{G^2}{H^2 \Xi^{-1} \Sigma} - s \\
& =  \frac{(A_\infty^\ast)^2(G - A_\infty^\ast)}{-0.5A_\infty^\ast \left( H^2 \Xi^{-1} \Sigma + (G - A_\infty^\ast)^2   \right) + b^2H^2 \Xi^{-1}}
\end{align*}
For the claim \eqref{eq:ropt0}, when $s = 0$, it follows directly from \eqref{eq:optsr} that
\begin{align}
\label{eq:ropts0}
r^{\text{opt}}_0 = \frac{(A_\infty^\ast)^2 - G^2}{H^2 \Xi^{-1} \Sigma }
\end{align}
First establish the following bound on $A_\infty^\ast$, 
\begin{align*}
A_\infty^\ast &= \left(-\frac{q}{2} + \sqrt{\tau}\right)^{1/3} + \left(-\frac{q}{2} - \sqrt{\tau}\right)^{1/3} \\
& <  \left(-\frac{q}{2} + \sqrt{\frac{q^2}{4} + \frac{p^3}{27}}\right)^{1/3} + \left(-\frac{q}{2} - \sqrt{\tau}\right)^{1/3} \\
& < \left(-\left(-\frac{p^3}{27} \right)^{1/2} + \sqrt{\frac{q^2}{4} + \frac{p^3}{27}}\right)^{1/3} + \left(-\frac{q}{2} - \sqrt{\tau}\right)^{1/3}\\
& <  \left(-\left(-\frac{p^3}{27} \right)^{1/2} \right)^{1/3} + \left(-\frac{q}{2} - \sqrt{\tau}\right)^{1/3}
\end{align*}
recalling that $p = -(H^2 \Xi^{-1} \Sigma+G^2), \enskip q = 4b^2 H^2 \Xi^{-1}$.  The condition $\tau > 0$ also implies 
\begin{align*}
-\frac{q}{2} < -\left(-\frac{p^3}{27} \right)^{1/2}
\end{align*}
Then using the inequality $\sqrt{2}(a_1 + a_2)^{1/2} \geq \sqrt{a_1} + \sqrt{a_2}$, 
\begin{align*}
-\frac{q}{2} -\sqrt{\tau} < -\frac{q}{2} -\frac{1}{\sqrt{2}} \left( \frac{q}{2} + \left| \frac{p^3}{27} \right|^{1/2}  \right) = -\frac{\sqrt{2} + 1}{\sqrt{2}} \frac{q}{2} - \frac{1}{\sqrt{2}} \left| \frac{p^3}{27} \right|^{1/2} < -\left(\frac{\sqrt{2} + 2}{\sqrt{2}} \right) \left| \frac{p^3}{27} \right|^{1/2}
\end{align*}
Putting altogether, we have
\begin{align*}
A_\infty^\ast &<  -\frac{2\sqrt{2} +2}{\sqrt{6}}\sqrt{|p|} 
\end{align*}
Then letting $c := \frac{2\sqrt{2} +2}{\sqrt{6}}$ and substituting into \eqref{eq:ropts0} yields 
\begin{align*}
(A_\infty^\ast )^2 > c^2(H^2 \Xi^{-1} \Sigma + G^2)
\end{align*}
and 
\begin{align*}
r_0^{\text{opt}} > \frac{(c^2 - 1)G^2 + c^2H^2\Xi^{-1}\Sigma}{H^2\Xi^{-1}\Sigma} > \frac{(c^2-1)(G^2 + H^2\Xi^{-1}\Sigma)}{H^2\Xi^{-1}\Sigma} > 1
\end{align*}
since $c^2 > 1$, which yields \eqref{eq:ropt0}.  Then for the remainder of the claim, we have (using $r$ in place of $r_0^{\text{opt}}$ for brevity),
\begin{align*}
\frac{C_\infty}{\hat{C}_\infty} &= \frac{1}{r}  \left( \frac{G + \sqrt{G^2 + rH^2 \Xi^{-1} \Sigma}}{G + \sqrt{G^2 + H^2 \Xi^{-1} \Sigma}} \right) \\
& \leq  \frac{1}{r}  \left( \frac{2|G| + \sqrt{rH^2 \Xi^{-1} \Sigma}}{G + \sqrt{G^2 + H^2 \Xi^{-1} \Sigma}} \right) \\
&= \frac{1}{r} \frac{2|G|}{(G + \sqrt{G^2 + H^2\Xi^{-1} \Sigma})} + \frac{1}{\sqrt{r}} \frac{\sqrt{H^2\Xi^{-1} \Sigma}}{(G + \sqrt{G^2 + H^2\Xi^{-1} \Sigma})} \\
& \leq \frac{1}{\sqrt{r}} \left( \frac{2|G| + \sqrt{H^2\Xi^{-1} \Sigma}}{(G + \sqrt{G^2 + H^2\Xi^{-1} \Sigma})} \right) 
\end{align*}
whenever $r > 1$.  Furthermore, when $G > 0$ we have that
\begin{align*}
    \frac{2G + \sqrt{H^2 \Xi^{-1} \Sigma}}{G + \sqrt{G^2 + H^2 \Xi^{-1} \Sigma} } \leq \sqrt{2} 
\end{align*}
Finally, due to \eqref{eq:ropt0}
\begin{align*}
    \frac{1}{\sqrt{r}} \leq \frac{\sqrt{H^2 \Xi^{-1} \Sigma}}{\sqrt{2} \sqrt{G^2 + H^2 \Xi^{-1} \Sigma}} < \frac{1}{\sqrt{2}}
\end{align*}
Multiplying the two inequalities yields the final claim \eqref{eq:boundCratio}.


\end{proof}

\section{Conclusion}
We presented a detailed investigation of connections between continuous time, continuous trait Crow-Kimura replicator-mutator dynamics \cite{kimura_stochastic_1965} and the fundamental equation of non-linear filtering, the Kushner-Stratonovich partial differential equation.  Inspired by a non-local fitness functional presented in the mathematical biology literature \cite{cressman2006stability}, we extended this connection to obtain a ``modified'' Kushner-Stratonovich equation.  This equation was shown to beneficial for filtering with misspecified models and a specific choice of parameters in the fitness functional was shown to coincide with covariance inflated Kalman Bucy filtering, in the linear-Gaussian setting.  Additionally, we considered the misspecified model filtering problem, with linear-Gaussian dynamics and where the misspecification arises through an unknown constant bias in the signal dynamics.  We proved that through a judicious choice of parameters in the fitness functional, mean squared error and uncertainty quantification (through the covariance) could be improved via this modified Kushner-Stratonovich equation.  Estimation is improved over traditional covariance inflation techniques, as well as over the standard filtering setup (assuming perfect model knowledge).\\  
\\
There are several avenues for further work, most notably, the analysis on misspecified models in Section \ref{sec:misspec} has primarily focused on the scalar setting which has simplified the analysis.  In future works, the multivariate setting, as well as extensions to nonlinear dynamics should be explored.  Additionally, it would be worthwhile to extend the mode of convergence in Theorem \ref{theo:limitcrowkimura} to $L^p$ convergence rather than pointwise convergence.  

\newpage 

\section*{Acknowledgements}
SP gratefully acknowledges funding from the Eva Mayr Stihl foundation through membership in the Young Academy for Sustainability research.  The authors also thank Aimee Maurais for bringing the works \cite{belhadji_weighted_2025, zhu_kernel_2024} to our attention and pointing out the connection to Fisher-Rao gradient flows of the MMD.

\begin{appendix}
\section{Technical lemmata}

\begin{theorem}
\label{theo:pdeprobrep}
\textbf{Forward representation formula} (Theorem 3.1 in \cite{kunita_stochastic_1982}.) Let $u_t(x, \omega)$ for $t \in [0, T]$, $x \in \mathbb{R}^d$ denote a measurable stochastic process on a probability space $(\Omega, \mathcal{F}, \mathbb{P})$ with $\omega \in \Omega$ (from now on we suppress the $\omega$ notation).  Suppose its time evolution is given by 
\begin{align}
\label{eq:parapde2}
\partial_t u_t(x) = \mathcal{L}u_t(x) + \sum_{k=1}^n \mathcal{M}^k u_t(x) \circ dB_t^k, \quad u_0(x) = f(x)
\end{align}
where $B_t^k$ is a Brownian motion wrt $\mathbb{P}$ and $f$ a bounded $C^2$ function with bounded derivatives and 
\begin{align*}
\mathcal{L} &:= \frac{1}{2} \sum_{j=1}^m \left( \sum_{i=1}^m a^{ij}(x) \frac{\partial }{\partial x^i} \right)^2 + \sum_{i=1}^m b^{i}(x) \frac{\partial }{\partial x^i} + c^0 \\
& = \frac{1}{2} \sum_{i=1}^m \sum_{j=1}^m \left( \sum_k a^{ik}a^{jk} \right)  \frac{\partial^2}{\partial x^i \partial x^j} + \frac{1}{2} \sum_{i,j} a^{ij} \sum_k \left( \frac{\partial a^{kj} }{\partial x^i}\right) \frac{\partial}{\partial x^k} + \sum_{i=1}^m b^{i}(x) \frac{\partial }{\partial x^i} + c^0 \\
\mathcal{M}^k &:= \sum_{i=1}^m l^{ik}(x) \frac{\partial }{\partial x^i} + c^k
\end{align*}
with $c^k(x)$ being uniformly bounded $C^2$ functions with bounded derivatives in $x$ and $a^{ij}(x), b^i(x), m^{ik}(x)$ being uniformly bounded $C^4$ functions with bounded first derivatives in $x$.\\  
\\
Then there exists another probability space $(\tilde{\Omega}, \mathcal{B}, \mathbb{Q})$ on which the Brownian motion $W_t = [W_t^1,\dots, W_t^m]^\top$ is defined and an SDE on the product space $(\Omega \otimes \tilde{\Omega}, \mathcal{F} \otimes \mathcal{B}, \mathbb{P} \otimes \mathbb{Q})$ given by 
\begin{align}
\label{eq:sdefeynmankac}
d \xi_t(x) =  \sum_{i=1}^d b^{i}(\xi_t(x))  dt + \sum_{j=1}^m   \sum_{i=1}^m a^{ij}(\xi_t(x))  \circ dW_t^j + \sum_{k=1}^n  \sum_{i=1}^m l^{ik}(\xi_t(x)) \frac{\partial \xi_t(x)}{\partial x^i}  \circ dB_t^k 
\end{align}
where the notation $\xi_t(x)$ is used to denote the solution of the SDE with initial condition $\xi_0 = x$ and $t > 0$. The solution of \eqref{eq:parapde2} has the representation 
\begin{align*}
u_t(x) = \mathbb{E}^{\mathbb{Q}} \left[ f(\xi_t(x)) \exp \left( \sum_{k=1}^n \int_0^t c^k(\xi_s(x)) \circ dB_s^k  + \int_0^t c^0(\xi_s(x)) ds\right) \right], \quad x \in \mathbb{R}^m
\end{align*}
\end{theorem}

The following form of Ito-Stratonovich correction will also be useful.  Consider the following Stratonovich SDE
\begin{align*}
d\xi_t = X_0(t, \xi_t)dt + \sum_{j=1}^m X_j(t, \xi_t) \circ dW_t^j 
\end{align*}
where $X_j$ is a Lipschitz cts function with derivatives in the second argument up to second order. Let $\xi_{0,r}(x)$ denote the solution of the above Stratonovich SDE at time $r$ with initial condition $\xi_0 = x$.  Then the following relation between the Ito and Stratonovich integral holds (backward form!)
\begin{align}
\label{eq:moditostrat}
\int_s^t X^k( f(\xi_{r,t}(x)) \circ dB_r^k = \int_s^t X^k (f(\xi_{r,t}(x)))dW_r^k + \frac{1}{2} \int_s^t (X^k( f(\xi_{r,t}(x))))^2dr 
\end{align}
where $f:\mathbb{R}^d \rightarrow \mathbb{R}$ is a $C^3$ function.\\

\noindent Some well-known results from stochastic analysis are presented below. 
\begin{lemma}
\label{lem:mgf}
\textbf{Moment generating function}.  Consider 
\begin{align*}
Y_t = \int_0^t \vartheta_s dW_s 
\end{align*}
where $\vartheta_s$ is a deterministic, scalar-valued continuous function of $s$ and $W_t$ is a real-valued Brownian motion wrt a probability measure $\mathbb{P}$.  Then 
\begin{align*}
\mathbb{E}^\mathbb{P}[\exp(\lambda Y_t)] = \exp \left( \frac{\lambda}{2} \int_0^t \vartheta_s^2 ds \right) 
\end{align*}
\end{lemma}

\begin{lemma}
Consider a probability space $(\Omega, \mathcal{F}, \mathbb{P})$ on which a scalar valued Wiener process is defined.  Suppose $f(s,\omega)$ for $\omega \in \Omega$ is an $\mathcal{F}_t$-adapted process.  Then the following inequality holds for $q > 1$
\label{lem:lpineq}
\begin{align}
\mathbb{E} \left[ \left| \int_0^t f(s,\omega)dW_s \right|^{2q} \right] \leq t^{q-1}(q(2q-1))^q \mathbb{E} \left[ \int_0^t  | f(s,\omega) |^{2q} ds  \right]
\end{align}
\end{lemma}

The following well-known lemma will also be used 
\begin{lemma}
\label{lem:lipsde}
\textbf{Continuity of solutions of SDEs} 
Let $g, \sigma$ be lipschitz continuous functions satisfying linear growth conditions. Denote by $\xi_t(x)$ the unique solution to 
\begin{align*}
d\xi_t =  g(\xi_s)ds +  \sigma(\xi_s)dB_s 
\end{align*}
with $\xi_0 = x$.  Suppose $h$ is a globally lipschitz continuous real vector valued function.  Then it holds that 
\begin{align*}
\mathbb{E} \left[ |h(\xi_{s}(x)) - h(\xi_r(x))|^q  \right] \leq C(x)|s-r|^{q/2}
\end{align*}
\end{lemma}

\begin{lemma}
\label{lem:delmeyer}
\textbf{Exponential moment bounds.}  Let $Y_s$ for $s \in [0,t]$ denote a non-decreasing scalar-valued adapted process such that $\mathbb{E}[Y_t - Y_s|\mathcal{F}_s] \leq K$ for all $s \in [0,t]$.  Then for any $L < \frac{1}{K}$, 
\begin{align*}
\mathbb{E}[\exp \left( L Y_t \right)] < \frac{1}{1 - LK}
\end{align*}
\end{lemma}

\begin{lemma}
\label{lem:unbiasedness} 
\textbf{Unbiasedness of the non-local replicator-mutator with perfect system.} Consider the signal-observation pair \eqref{eq:sig}-\eqref{eq:obs} with $g(x) = Gx, \enskip h(x) = Hx$ and Gaussian initial conditions. If $m_0= \mathbb{E}[X_0]$ then the generalised Kalman-Bucy filter with $s < r, r > 0$ is unbiased, i.e. $\mathbb{E}[m_t] = \mathbb{E}[X_t], \enskip t > 0$. 
\begin{proof}
Begin with the mean equation from the generalised Kalman-Bucy filter
\begin{align*}
	m_t &= m_0 + \int_0^t Gm_u du + (r-s)\int_0^t K_u dZ_u -  (r-s)\int_0^t K_u H m_u du \\
	& = m_0 + \int_0^t Gm_u du + (r-s)\int_0^t K_u H (X_u^{\ast} - m_u)du + (r-s)\int_0^t K_u\Xi^{1/2} dB_u
\end{align*}
where $K_t:= C_t H^\top \Xi^{-1}$ corresponds to the Kalman gain at time $t$ and we have substituted in the form of the observations in the second equality.  Then 
\begin{align*}
	\mathbb{E}[m_t] &= m_0 + \int_0^t G \mathbb{E}[m_u]du + (r-s)\int_0^t K_u H (\mathbb{E}[X_u^{\ast}] - \mathbb{E}[m_u])du 
\end{align*}
Also recall that 
\begin{align*}
	\mathbb{E}[X_t^\ast] = \mathbb{E}[X_0^\ast] + \int_0^t G \mathbb{E}[X_u^\ast]du
\end{align*}
then 
\begin{align*}
	\mathbb{E}[m_t] - \mathbb{E}[X_t^\ast] = m_0 - \mathbb{E}[X_0^\ast] + \int_0^t (G - (r-s)K_sH)(\mathbb{E}[m_s] - \mathbb{E}[X_s^\ast])ds
\end{align*}
which corresponds to an ODE of the form $\dot{y}_t =A_ty_t$ with $y_t := \mathbb{E}[m_t] - \mathbb{E}[X_t^\ast]$, which has solutions $y_t = 0, \enskip t > 0$ whenever $y_0 = 0$. 
\end{proof}
\end{lemma}

\end{appendix}

\printbibliography


\end{document}